\documentclass[12pt]{amsart}
\usepackage{latexsym}
\usepackage{amssymb}
\usepackage{amscd}
\renewcommand\a{\alpha}
\renewcommand\b{\beta}

\newcommand\la{\lambda}

\newcommand\e{\eta}
\renewcommand\th{\theta}
\newcommand\io{\iota}
\newcommand\m{\mu}

\newcommand\f{\phi}
\newcommand\vf{\varphi}
\newcommand\p{\psi}

\renewcommand\r{\rho}
\newcommand\Om{\Omega}
\newcommand\w{\omega}

\newcommand\vS{\varSigma}

\newcommand\vD{\varDelta}

\newcommand\vL{\varLambda}

\newcommand{\vT}{\varTheta}

\newcommand\ve{\varepsilon}

\newcommand{\ZZ}{\mathbb Z}

\newcommand\BP{\mathbf P}

\newcommand\Ba{\mathbf a}

\newcommand\Bp{\mathbf p}
\newcommand\Bm{\mathbf m}
\newcommand\Bb{\mathbf b}

\newcommand\CB{\mathcal{B}}
\newcommand\ZC{\mathcal{C}}
\newcommand\CH{\mathcal{H}}

\newcommand\CS{\mathcal{S}}
\newcommand\CM{\mathcal{M}}

\newcommand\CK{\mathcal{K}}
\newcommand\CP{\mathcal{P}}

\newcommand\CT{ \mathcal{T}}

\newcommand\FS{\mathfrak S}
\newcommand\Fu{\mathfrak u}
\newcommand\Fv{\mathfrak v}

\newcommand\Fs{\mathfrak s}

\newcommand\Ft{\mathfrak t}

\newcommand\iv{^{-1}}
\newcommand\wh{\widehat}
\newcommand\wt{\widetilde}

\newcommand\ol{\overline}

\newcommand\trreq{\trianglerighteq}
\newcommand\dia{\diamondsuit}
\newcommand\hra{\hookrightarrow}
\newcommand\lra{\leftrightarrow}

\newcommand\Hom{\operatorname{Hom}}
\newcommand\End{\operatorname{End}}

\newcommand\Id{\operatorname{Id}}
\newcommand\lp{\operatorname{\!\langle\!}}
\newcommand\rp{\operatorname{\!\rangle\!}}
\renewcommand\Im{\operatorname{Im}}

\newcommand\nat{^{\natural}}

\newcommand\Std{\operatorname{Std}}
\newcommand\rStd{\operatorname{r-Std}}

\newcommand\rad{\operatorname{rad}}

\newcommand{\isom}{\,\raise2pt\hbox{$\underrightarrow{\sim}$}\,}
\numberwithin{equation}{section}

\newtheorem{thm}{Theorem}[section]
\newtheorem{lem}[thm]{Lemma}
\newtheorem{cor}[thm]{Corollary}
\newtheorem{prop}[thm]{Proposition}

\def \para#1{\par\medskip\textbf{#1}
              \addtocounter{thm}{1}}

\def \remark#1{\par\medskip\noindent
                \textbf{Remark #1}
                \addtocounter{thm}{1}}

\def \remarks#1{\par\medskip\noindent
                \textbf{Remarks #1}
                \addtocounter{thm}{1}}

\begin{document}
\setlength{\baselineskip}{4.9mm}
\setlength{\abovedisplayskip}{4.5mm}
\setlength{\belowdisplayskip}{4.5mm}
%%%
%%%
\renewcommand{\theenumi}{\roman{enumi}}
\renewcommand{\labelenumi}{(\theenumi)}
\renewcommand{\thefootnote}{\fnsymbol{footnote}}
%%%
\renewcommand{\thefootnote}{\fnsymbol{footnote}}
%%%
\allowdisplaybreaks[2]
%\NoBlackBoxes
\parindent=20pt
%%%%%%%%%%%%%%%%%%%%
%%%%%%%%%%%%%%%%%%%%%%%%%%%%%%%%%%%
\medskip
\begin{center}
{\bf Cyclotomic $q$-Schur algebras associated to the Ariki-Koike algebra} 
\\
\vspace{1cm}
Toshiaki Shoji and Kentaro Wada 
\\ 
\vspace{0.5cm}
Graduate School of Mathematics \\
Nagoya University  \\
Chikusa-ku, Nagoya 464-8602,  Japan
\end{center}
\title{}
\maketitle
\begin{abstract}
Let $\CH_{n,r}$ be the Ariki-Koike algebra associated to 
the complex reflection group $\FS_n\ltimes (\ZZ/r\ZZ)^n$, 
and $\CS(\vL)$ be the cyclotomic $q$-Schur algebra associated to
$\CH_{n,r}$, introduced by Dipper-James-Mathas.  For each 
$\Bp = (r_1, \dots, r_g) \in \ZZ_{>0}^g$ such that 
$r_1 +\cdots + r_g = r$, we define a subalgebra
$\CS^{\Bp}$ of $\CS(\vL)$ and its quotient algebra $\ol\CS^{\Bp}$.
It is shown that $\CS^{\Bp}$ is a standardly based algebra and 
$\ol\CS^{\Bp}$ is a cellular algebra.  By making use of these algebras, 
we prove a product formula for decomposition numbers of $\CS(\vL)$, 
which asserts that certain decomposition numbers are expressed 
as a product of 
decomposition numbers for various cyclotomic $q$-Schur algebras 
associated to Ariki-Koike algebras $\CH_{n_i,r_i}$ of smaller rank. 
This is a generalization of the result of N. Sawada.  
We also define a modified Ariki-Koike algebra $\ol\CH^{\Bp}$ of type
$\Bp$, and prove the Schur-Weyl duality between $\ol\CH^{\Bp}$ and 
$\ol\CS^{\Bp}$. 
\end{abstract}

\pagestyle{myheadings}
\markboth{SHOJI}{CYCLOTOMIC $q$-SCHUR ALGEBRAS}

%%%%%%%%%%%%%%%%%%%%
%%%%%%%%%%%%%%%%%%%%%%%%%%%%%%%%%%%
\bigskip
\medskip
\addtocounter{section}{-1}
\section{Introduction}
Let $\CH = \CH_{n,r}$ be the Ariki-Koike algebra over an integral 
domain $R$ associated to 
the complex reflection group $W_{n,r} = \FS_n\ltimes (\ZZ/r\ZZ)^n$
with parameters $q, Q_1, \dots, Q_r \in R$ such that $q$ is a unit 
in $R$.
Let $\wt\CP_{n,r}$ (resp. $\CP_{n,r}$) be the set of 
$r$-compositions (resp. $r$-partitions) of $n$.
The cyclotomic $q$-Schur algebra $\CS(\vL)$ associated to
$\CH$ was introduced by Dipper-James-Mathas [DJM], which is 
the endomorphism algebra of a certain $\CH$-module 
$M = \bigoplus_{\mu \in \vL}M^{\mu}$, where 
$\vL$ is a saturated subset of $\wt\CP_{n,r}$.
They showed that $\CS(\vL)$ is a cellular algebra in the sense of
Graham-Lehrer [GL], and that
the Schur-Weyl duality (i.e., the double centralizer property)
holds between $\CH$ and $\CS(\vL)$ in the case where $\vL = \wt\CP_{n,r}$.
\par
On the other hand, the modified Ariki-Koike algebra $\ol\CH$ was 
introduced in [SawS], under the condition that (*) ``$Q_i - Q_j$ are units 
in $R$ for each $i \ne j$'', based on the study of the Schur-Weyl
duality between $\CH$ and a certain subalgebra of the quantum group 
of type $A$ ([SakS], [Sh]).
By using the cellular structure of $\ol\CH$, a cyclotomic
$q$-Schur algebra associated to $\ol\CH$ was constructed,   
in analogy to $\CS(\vL)$. 
It was shown in [SawS] that this cyclotomic $q$-Schur algebra is 
isomorphic to the quotient algebra $\ol\CS^0$ of a certain subalgebra
$\CS^0$ of $\CS(\vL)$, and that the Schur-Weyl duality holds between 
$\ol\CH$ and $\ol\CS^0$.
Moreover, the structure theorem for $\ol\CS^0$was proved, which asserts that
$\ol\CS^0$ is a direct sum of tensor products of various $q$-Schur 
algebras $\CS(\wt\CP_{n_i,1})$ associated to the Iwahori-Hecke algebra 
of type $A_{n_i-1}$.
\par
In [Sa], N. Sawada reconstructed the subalgebra $\CS^0$ of $\CS(\vL)$ 
and its quotient $\ol\CS^0$ based on the cellular structure on 
$\CS(\vL)$, which works without the assumption (*).  
He proved that $\CS^0$ is a standardly based algebra in the sense of 
Du and Rui [DR], and showed, in the case where $R$ is a filed, 
 that the decomposition number $d_{\la\mu}$ 
between the Weyl module
$W^{\la}$ and the irreducible module $L^{\mu}$ of $\CS(\vL)$ 
($\la, \mu \in \CP_{n,r}$) coincides
with the corresponding decomposition number for $\ol\CS^0$ whenever
$|\la^{(i)}| = |\mu^{(i)}|$ for $i = 1, \dots, r$.
This implies in the case where $\vL = \wt\CP_{n,r}$, under 
the condition (*), that $d_{\la\mu}$ can be written 
as a product of $d_{\la^{(i)}\mu^{(i)}}$ for $i = 1, \dots, r$, 
where $d_{\la^{(i)}\mu^{(i)}}$ is the decomposition number of the 
$q$-Schur algebra $\CS(\wt\CP_{n_i,1})$ with 
$|\la^{(i)}| = |\mu^{(i)}| = n_i$.
\par
The subalgebra $\CS^0$ is regarded, in some sense, as
a Borel type subalgebra of $\CS(\vL)$. 
For example, we have $\CS(\vL) = \CS^0\cdot(\CS^0)^*$, where $(\CS^0)^*$
is the image of $\CS^0$ under the involution $*$ of $\CS(\vL)$, and 
$\ol\CS^0$ is a quotient of both $\CS^0$ and $(\CS^0)^*$.  Thus 
$\ol\CS^0$ corresponds to a Cartan subalgebra.  In this paper, we consider 
a parabolic analogue of $\CS^0$ and $\ol\CS^0$.
We fix $\Bp = (r_1, \dots, r_g) \in \ZZ_{>0}^g$ such that 
$r_1 + \cdots r_g = r$.  According to $\Bp$, we regard an
$r$-partition $\la = (\la^{(1)}, \dots, \la^{(r)})$  as a $g$-tuple
of multi-partitions $\la = (\la^{[1]}, \dots, \la^{[g]})$, 
where $\la^{[1]} = (\la^{(1)}, \dots, \la^{(r_1)}), 
\la^{[2]} = (\la^{(r_1+1)}, \dots, \la^{(r_1+r_2)})$, and so on. 
For example $\la^{[i]} = \la^{(i)}$ for $i = 1, \dots, r$ if 
$\Bp = (1^r)$ with $g = r$, and $\la^{[1]} = \la$ if $\Bp = (r)$ 
with $g = 1$.  For each $\Bp$, we define a subalgebra $\CS^{\Bp}$ of 
$\CS(\vL)$, and its quotient algebra $\ol\CS^{\Bp}$.  
The algebra $\CS^{\Bp}$ coincides with $\CS^0$ if
$\Bp = (1^r)$, and coincides with $\CS(\vL)$ if $\Bp = (r)$.
Thus $\CS^{\Bp}$ is a generalization of $\CS^0$, and is regarded 
as an intermediate object between $\CS(\vL)$ and $\CS^0$. 
\par
All the results in [Sa] can be generalized to our cases;
$\CS^{\Bp}$ is a standardly based algebra and 
$\ol\CS^{\Bp}$ is a cellular algebra. 
Assume that $R$ is a field. For $\la = (\la^{[1]}, \dots, \la^{[g]}), 
\mu = (\mu^{[1]}, \dots, \mu^{[g]}) \in \CP_{n,r}$ such that 
$|\la^{[i]}| = |\mu^{[i]}|$ for $i = 1, \dots, g$, one can show 
(Theorem 3.13) that
the decomposition number $d_{\la\mu}$ coincides with the corresponding 
decomposition number in the algebra $\ol\CS^{\Bp}$.
In the case where $\vL = \wt\CP_{n,r}$, we prove the structure 
theorem (Theorem 4.15) for $\ol\CS^{\Bp}$, which asserts 
that $\ol\CS^{\Bp}$ is a 
direct sum of tensor products of various $\CS(\wt\CP_{n_i,r_i})$.
We remark, contrast to the argument in [SawS], that no assumptions 
on parameters are required in this proof.
Combining with the previous results, we obtain the product formula
for decomposition numbers, namely, $d_{\la\mu}$ 
coincides with the product of $d_{\la^{[i]}\mu^{[i]}}$ for 
$i =1, \dots, g$, where $d_{\la^{[i]}\mu^{[i]}}$ is the decomposition 
number for $\CS(\wt\CP_{n_i,r_i})$ with 
$|\la^{[i]}| = |\mu^{[i]}| = n_i$, which holds without 
any restriction on parameters (Theorem 4.17).
\par
By making use of the Schur functors on $\CS(\vL)$, 
one can define a modified
Ariki-Koike algebra $\ol\CH^{\Bp}$ of type $\Bp$ as a 
certain subalgebra of $\CS^{\Bp}$.  The algebra $\ol\CH^{\Bp}$ 
is isomorphic to  
$\ol\CH$ if $\Bp = (1^r)$, and coincides with $\CH$ if $\Bp = (r)$.
Put $Q_i^{\Bp} = Q_{r_1 + \cdots + r_i}$ for $i = 1, \dots, g$.
Under the assumption (**) ``$Q_i^{\Bp} - Q_j^{\Bp}$ are units in $R$
for each $i \ne j$'',  we give a presentation of $\ol\CH^{\Bp}$ which
is a generalization of the presentation of $\ol\CH$ given in [SawS].
We show that $\ol\CS^{\Bp}$ is realized as an endomorphism algebra
of a certain $\ol\CH^{\Bp}$-module 
$\ol M^{\Bp} = \bigoplus_{\mu \in \vL}\ol M^{\mu}$,   
and prove the Schur-Weyl duality between $\ol\CS^{\Bp}$ and 
$\ol\CH^{\Bp}$.  
In the case where the parameters $q, Q_1, \dots, Q_r$ satisfy the 
separation condition in the sense of [A] (see (8.3.1)), it is shown 
that all the $\ol\CH^{\Bp}$ are isomorphic to $\CH$, and so the above
results give new presentations of $\CH$.  
\par
By using the Jantzen filtration, $v$-decomposition numbers
$d_{\la\mu}(v)$ for $\CS(\vL)$ can be defined, which is a polynomial 
analogue of $d_{\la\mu}$.  The results in this paper concerning the 
decomposition numbers for $\CS(\vL), \CS^{\Bp}, \ol\CS^{\Bp}$ 
are generalized to $v$-decomposition numbers.
In particular, the product formula for $v$-decomposition numbers is 
obtained, which is discussed in [W].
\par\bigskip\noindent
{\bf Notation}
\par
Let $R$ be an integral domain and $M$ a free $R$-module
of finite rank.  We denote by $\End M$ the endomorphism 
algebra of $M$, where the composition is defined by 
$(f\circ g)(m) = f(g(m))$ for $f,g \in \End M, m \in M$. 
Thus $\End M$ acts on $M$ from the left by $(f,m) \mapsto f(m)$.
We denote by $\End^0M$ the opposite algebra of $\End M$. 
If an $R$-algebra $A$ (resp. $B$) acts on $M$ from the left
(resp. from the right), then we have a natural homomorphism of 
$R$-algebras $A \to \End M$ (resp. $B \to \End^0M$).
If an $R$-algebra $X$ acts on $M$ from the right or left, we denote
by $\End_XM$ the subalgebra of $\End M$ consisting of endomorphisms 
commuting with $X$.  The subalgebra $\End^0_XM$ of $\End^0M$ is 
defined similarly. 
\par\bigskip\noindent
{\bf Table of contents}
\par\medskip\noindent 
0. Introduction \\
1. Recollection of cyclotomic $q$-Schur algebras \\ 
2. Parabolic type subalgebras of $\CS(\vL)$  \\
3. Decomposition numbers for $\CS(\vL), \CS^{\Bp}$ and 
   $\ol\CS^{\Bp}$  \\
4. Structure theorem for $\ol\CS^{\Bp}$  \\
5. Modified Ariki-Koike algebra of type $\Bp$  \\
6. Presentation for $\ol\CH^{\Bp}$ \\
7. Schur-Weyl duality  \\
8. Comparison of $\ol\CH^{\Bp}$ for various $\Bp$

\section{Recollection of cyclotomic $q$-Schur algebras}
\para{1.1.}
Let $R$ be an integral domain, $q, Q_1, \dots, Q_r$ be elements in 
$R$ with $q$ invertible.  The Ariki-Koike algebra $\CH = \CH_{n,r}$
associated to the complex reflection group $\FS_n\ltimes (\ZZ/r\ZZ)^n$
is an associative
algebra over $R$ with generators $T_0, T_1, \dots, T_{n-1}$ subject to 
the condition
\begin{equation*}
\begin{aligned}
&(T_0 -Q_1)\cdots (T_0-Q_r) = 0, &   &\\
&(T_i-q)(T_i+q\iv) = 0 &\quad (&i \ge 1), \\ 
&T_0T_1T_0T_1 = T_1T_0T_1T_0,  &   &\\
&T_iT_j = T_jT_i &\quad (&|i-j| \ge 2),  \\
&T_iT_{i+1}T_i = T_{i+1}T_iT_{i+1} &\quad (&1 \le i \le n-2). 
\end{aligned} 
\end{equation*}
It is known that $\CH$ is a free $R$-module with rank $n!r^n$.
We denote by $\CH_n$ the subalgebra of $\CH$ generated by 
$T_1, \dots, T_{n-1}$, which is isomorphic to the Iwahori-Hecke
algebra associated to the symmetric group $\FS_n$ of degree $n$.
\para{1.2.}
It is known by [DJM] that $\CH$ has a structure of 
the cellular algebra.  In order to describe the cellular basis
of $\CH$, we prepare some notation.
An element $\mu = (\mu_1, \dots, \mu_m) \in \ZZ_{\ge 0}^m$ 
is called a composition of length $\le m$, 
and $|\mu| = \sum \mu_i$ is called the size of $\mu$.
An $r$-composition $\la = (\la^{(1)}, \dots, \la^{(r)})$ 
is an $r$-tuple of compositions 
$\la^{(i)} = (\la^{(i)}_1, \dots, \la^{(i)}_{m_i})$.  The size 
$|\la|$ of $\la$ is defined by $|\la| = \sum_{i=1}^r |\la^{(i)}|$.
We denote by $\la$ by $\la = (\la^{(i)}_j)$.
A composition $\mu = (\mu_1, \dots, \mu_m)$ is called a partition
if $\mu_1 \ge \cdots \ge \mu_m \ge 0$.  An $r$-composition $\la$ is called 
an $r$-partition if $\la^{(i)}$ is a partition for all $i$.
We fix $\Bm = (m_1, \dots, m_r) \in \ZZ_{> 0}^r$ once and for all, 
and denote by 
$\wt\CP_{n,r} = \wt\CP_{n,r}(\Bm)$ the set of $r$-compositions 
$\la = (\la^{(1)}, \dots, \la^{(r)})$ of size $n$ such that 
$\la^{(i)}$ is a composition of length $\le m_i$.  
Similarly, we define the set $\CP_{n,r} = \CP_{n,r}(\Bm)$ of $r$-partitions.
If $m_i \ge n$ for any $i$, $\CP_{n,r}(\Bm)$ are mutually
identified with for all $\Bm$.  However even in that case, 
$\wt\CP_{n,r}(\Bm)$ depends on the choice of $\Bm$.
\par
For $r$-compositions $\la = (\la^{(i)}_j)$ and 
$\mu = (\mu^{(i)}_j)$, we define a dominance order 
$\la \trianglerighteq \mu$ by the condition 
\begin{equation*}
\sum_{c=1}^{k-1}|\la^{(c)}| + \sum_{j=1}^{i}|\la^{(k)}_j|
   \ge \sum_{c=1}^{k-1}|\mu^{(c)}| +
             \sum_{j=1}^i|\mu^{(k)}_j|
\end{equation*}
for any $1 \le k \le r$ and $1 \le i \le m_k$.
If $\la \trianglerighteq \mu$ and $\la \ne \mu$, we write
it as $\la \triangleright \mu$.
\par
Let $\la$ be an $r$-partition of $n$.  We identify $\la$ with 
the $r$-tuple of Young diagrams, and refer it as the Young diagram 
of $\la$. We denote by $\Std(\la)$ the 
set of standard tableau $\Ft = (\Ft^{(1)}, \dots, \Ft^{(r)})$ of 
shape $\la$, i.e., $\Ft$ is a Young diagram of $\la$ with  
letters $1, \dots, n$ attached to the nodes of the diagram, under the
condition that $\Ft^{(i)}$ is a standard tableau in the usual sense
for each $i$.
We define $\Ft^{\la} \in \Std(\la)$ by attaching the letters 
$1, 2, \dots, n$ to the nodes of
the Young diagram $\la$ in this order, from left to right, 
and from top to down for $\Ft^{(1)}$, and then for $\Ft^{(2)}$, and so
on.  
$\FS_n$ acts naturally on $\Std(\la)$ from the right, 
and we denote by $d(\Ft)$ 
the element in $\FS_n$ such that $\Ft = \Ft^{\la}d(\Ft)$ 
for each $\Ft \in \Std(\la)$.
More generally, the set $\rStd(\mu)$ of row-standard tableaux
of shape $\mu$ is defined for $\mu \in \wt\CP_{n,r}$, 
by replacing a standard tableau $\Ft^{(i)}$ by 
a row-standard tableau.  Then $\Ft^{\mu}$ is defined similarly, 
and $d(\Ft) \in \FS_n$ is defined 
also for $\Ft \in \rStd(\mu)$.
\par
For $\mu \in \wt\CP_{n,r}$, we define 
$r$-tuples of integers
\begin{equation*} 
\a(\mu) = (\a_1, \dots, \a_r), \qquad 
   \Ba(\mu) = (a_1, \dots, a_r)
\end{equation*}
 by 
$\a_i  = |\mu^{(i)}|$, and $a_i = \sum_{j=1}^{i-1}|\mu^{(j)}|$
for $i = 1, \dots, r$. (Note that $a_1 = 0$.)
\par
We define $L_k \in \CH$ by 
$L_k = T_{k-1}\cdots T_1T_0T_1\cdots T_{k-1}$ for 
$k = 1, \dots, n$.  Then $L_1, \dots, L_n$ commute with each other.
For $\Ba = (a_1, \dots, a_r) \in \ZZ_{\ge 0}^r$, we define
$u_{\Ba}^+ \in \CH$ by  
$u_{\Ba}^+ = u_{\Ba,1}u_{\Ba,2}\cdots u_{\Ba,r}$, where
\begin{equation*}
u_{\Ba,k} = \prod_{i=1}^{a_k}(L_i - Q_k).
\end{equation*}
For $\mu \in \wt\CP_{n,r}$, let   
$\FS_{\mu} = \FS_{\mu^{(1)}}\times\cdots\times \FS_{\mu^{(r)}}$
be the Young subgroup of $\FS_n$.
we define $x_{\la} \in \CH_n$ by 
$x_{\mu} = \sum_{w \in \FS_{\mu}}q^{l(w)}T_w$, where 
$l(w)$ is the length of $w \in \FS_n$, and $T_w$ is a basis element 
of $\CH_n$ corresponding to $w \in \FS_n$.
We define $m_{\mu} \in \CH$ by 
$m_{\mu} = u_{\Ba}^+x_{\mu} = x_{\mu}u_{\Ba}^+$. 
For $\Fs, \Ft \in \rStd(\mu)$, we define $m_{\Fs\Ft} \in \CH$
by $m_{\Fs\Ft} = T_{d(\Fs)}^*m_{\mu}T_{d(\Ft)}$, where 
$x \mapsto x^*$ is an anti-automorphism on $\CH_n$ defined by 
$T_i^* = T_i$ for $i = 1, \dots, n-1$. 
Then it is known by [DJM, Theorem 3.26] that 
the set 
\begin{equation*}
\tag{1.2.1}
\{ m_{\Fs\Ft} \mid \Fs, \Ft \in \Std(\la) 
         \text{ for some } \la \in \CP_{n,r}\}
\end{equation*}
gives a cellular basis of $\CH$ with respect to 
the dominance order on $\CP_{n,r}$ in the sense of [GL].
In particular, if we denote by $h \mapsto h^*$ the anti-automorphism
on $\CH$ defined by $T_i^*= T_i$ for $i = 0, \dots, n-1$, 
we have $m_{\Fs\Ft}^* = m_{\Ft\Fs}$.
\para{1.3.}
Here we recall the concept of semistandard tableau in the
case of multi-partitions due to [DJM]. We consider the set $X$ of 
pairs $(i,k)$ with $1 \le i \le n, 1 \le k \le r$, and define 
a total order on this set by $(i_1,k_1) < (i_2,k_2)$ if 
$k_1 < k_2$, or if $k_1 = k_2$ and $i_1 < i_2$.
For an $r$-partition $\la$ of $n$, 
a Tableau $T$ of shape $\la$ is defined as  
a Young diagram $\la$ with an element of $X$ attached to each node
of $\la$. 
For each $(i,k) \in X$, let
$\mu^{(k)}_i$ be the number of entries of $T$ containing $(i,k)$.
Then $\mu = (\mu^{(k)}_j)$ is an $r$-composition of $n$.
The Tableau $T$ is called a $\la$-tableau of type $\mu$. 
A Tableau $T = (T^{(1)}, \dots, T^{(r)})$ of shape $\la$ is called 
a semistandard tableau if it satisfies the properties;  
the entries of $T^{(i)}$ are weakly increasing along the rows, 
strictly increasing along the columns with respect to $X$, 
and furthermore the entries of $T^{(k)}$ consist of 
$(i,k')$ with $k' \ge k$.  
We denote by $\CT_0(\la, \mu)$ the set of semi standard tableau 
of shape $\la$ and type $\mu$ for $\la \in \CP_{n,r}$ and 
$\mu \in \wt\CP_{n,r}$.  Note that $\CT_0(\la,\mu)$ is empty
unless $\la \trianglerighteq \mu$.
\par
Let $\Ft$ be a standard tableau of shape $\la$.  For 
$\mu \in \wt\CP_{n,r}$, we construct a Tableau $\mu(\Ft)$ 
from $\Ft$ as follows; replace 
the entry $j$ in $\Ft$ by $(i,k)$ if $j$ appears in 
the $i$-th row of the $k$-th component $(\Ft^{\mu})^{(k)}$ of 
$\Ft^{\mu}$. $\mu(\Ft)$ is a $\la$-tableau of type $\mu$, 
but it is not necessarily semistandard.
\para{1.4.}
For each $\mu \in \wt\CP_{n,r}$, we define a right $\CH$-module
$M^{\mu}$ by $M^{\mu} = m_{\mu}\CH$.
It is known by [DJM, Theorem 4.14] that $M^{\mu}$ is a free $R$-module
with basis
\begin{equation*}
\tag{1.4.1}
\{ m_{S\Ft} \mid S \in \CT_0(\la,\mu), \Ft \in \Std(\la)
                      \text{ for some } \la \in \CP_{n,r}\},
\end{equation*}
where
\begin{equation*}
\tag{1.4.2}
m_{S\Ft} = \sum_{\substack{\Fs \in \Std(\la) \\
             \mu(\Fs) = S}}q^{l(d(\Fs)) + l(d(\Ft))}m_{\Fs\Ft}.
\end{equation*}
\par
A subset $\vL$ of $\wt\CP_{n,r}(\Bm)$ is called a saturated 
set if any partition $\la$ such that  $\la \trianglerighteq \mu$
for some $\mu \in \vL$ is contained in $\vL$. 
We denote by $\vL^+$ the set of $r$-partitions of $n$ contained in $\vL$.
Put $M = \bigoplus_{\mu \in \vL}M^{\mu}$.  The cyclotomic
$q$-Schur algebra $\CS(\vL)$ associated to $\CH$ (and to $\vL$) is defined by 
\begin{equation*}
\CS(\vL) = \End_{\CH}(M) = 
   \bigoplus_{\nu,\mu \in \vL}\Hom_{\CH}(M^{\nu}, M^{\mu}).
\end{equation*}
\par
We consider the structure of $\Hom_{\CH}(M^{\nu}, M^{\mu})$ for 
$\mu, \nu \in \vL$. 
Let $M^{\nu*} = (M^{\nu})^*$ be the image of $M^{\nu}$ under $*$.
We have $M^{\nu*} = \CH m_{\nu}$.  It is easy to see that 
for any $m \in M^{\nu*} \cap M^{\mu} = \CH m_{\nu} \cap m_{\mu}\CH$, 
the map $m_{\nu}h \mapsto mh$ ($h \in \CH$) gives rise to 
an $\CH$-module homomorphism $\vf_m : M^{\nu} \to M^{\mu}$.
It is known by [DJM, Corollary 5.17] that the map
$\vf \mapsto \vf(m_{\nu})$ gives an isomorphism of $R$-modules
\begin{equation*}
\tag{1.4.3}
\Hom_{\CH}(M^{\nu}, M^{\mu}) \to M^{\nu*} \cap M^{\mu}.
\end{equation*}
Suppose that $\mu, \nu \in \vL$, and $\la \in \vL^+$.
We define for $S \in \CT_0(\la,\mu), T \in \CT_0(\la,\nu)$
\begin{equation*}
\tag{1.4.4}
m_{ST} = \sum_{\substack{ \Fs,\Ft \in \Std(\la) \\
       \mu(\Fs) = S, \nu(\Ft) = T}}q^{l(d(\Fs)) + l(d(\Ft))}m_{\Fs\Ft}.
\end{equation*}
Then it is known by [DJM, Proposition 6.3] that the set
\begin{equation*}
\{ m_{ST} \mid S \in \CT_0(\la, \mu), T \in \CT_0(\la, \nu) 
                 \text{ for some } \la \in \vL^+\}
\end{equation*}
gives rise to a basis of $M^{\nu*} \cap M^{\mu}$.
We denote by $\vf_{ST}$ the element of 
$\Hom_{\CH}(M^{\nu}, M^{\mu})$ corresponding to $m_{ST}$ via 
the isomorphism (1.4.3).  Thus $\vf_{ST}$ is a map 
$M^{\nu} \to M^{\mu}$  defined by $\vf_{ST}(m_{\nu}h) = m_{ST}h$
for any $h \in \CH$.  
For each $\la \in \vL^+$, put 
$\CT_0(\la) = \bigcup_{\mu \in \vL}\CT_0(\la,\mu)$.
The fundamental result of Dipper-James-Mathas
is the following theorem.
%%%
\begin{thm}[{[DJM]}] %%% Theorem 1.5.
The cyclotomic $q$-Schur algebra $\CS(\vL)$ is a cellular algebra 
with a cellular basis
\begin{equation*}
\ZC(\vL) = \{ \vf_{ST} \mid S, T \in \CT_0(\la) \text{ for some }  
              \la \in \vL^+ \}.
\end{equation*}
\end{thm}
\para{1.6.}
For $\la \in \vL^+$, let $T^{\la} = \la(\Ft^{\la})$.
Then $T^{\la}$ is a semistandard tableau obtained from $\Ft^{\la}$
by replacing the entries $j$ in $(\Ft^{\la})^{(k)}$ by $(j,k)$.
Then $\Ft = \Ft^{\la}$ is the unique standard tableau such that 
$\la(\Ft) = T^{\la}$.  It follows that 
$m_{T^{\la}T^{\la}} = m_{\Ft^{\la}\Ft^{\la}} = m_{\la}$, and 
$\vf_{T^{\la}T^{\la}}$ is the identity element in 
$\Hom_{\CH}(M^{\la}, M^{\la})$.  We put 
$\vf_{\la} = \vf_{T^{\la}T^{\la}}$.
\par
For each $\la \in \vL^+$, we define $\CS^{\vee\la}$ as the $R$-submodule
of $\CS(\vL)$ spanned by $\vf_{ST}$, where  
$S, T \in \CT_0(\la', \mu)$ for various $\la' \in \vL^+$ such that 
$\la' \triangleright \la$,  and for various $\mu \in \vL$.
Then $\CS^{\vee\la}$ is a two-sided ideal of $\CS(\vL)$, and 
we define the Weyl module $W^{\la}$ as the right $\CS(\vL)$-submodule 
of $\CS(\vL)/\CS^{\vee\la}$ generated by the image of 
$\vf_{\la} \in \CS(\vL)$. 
For each $T \in \CT_0(\la)$, let $\vf_T$ be the image of 
$\vf_{T^{\la}T}$ in $W^{\la}$.  Then the following holds;  
$W^{\la}$ is an $R$-free module
with basis $\{ \vf_T \mid T \in \CT_0(\la)\}$.  There exists 
a canonical bilinear form $\lp\ , \ \rp$ on $W^{\la}$ determined by
\begin{equation*}
\vf_{T^{\la}S}\vf_{TT^{\la}} \equiv 
       \lp\vf_S,\vf_T\rp\vf_{T^{\la}T^{\la}} \mod \CS^{\vee\la}.
\end{equation*}
Let $\rad W^{\la} = \{ x \in W^{\la} \mid \lp x, y\rp = 0 
      \text{ for any } y \in W^{\la}\}$.
Then $\rad W^{\la}$ is an $\CS(\vL)$-submodule of $W^{\la}$. 
Put $L^{\la} = W^{\la}/\rad W^{\la}$.  Assume that $R$ is a field. 
Then it is known by [DJM] that  
$L^{\la}$ is a (non-zero) absolutely irreducible module, and that
the set $\{ L^{\la} \mid \la \in \vL^+\}$ gives a complete set
of non-isomorphic $\CS(\vL)$-modules.      
%%%
%%%
\section{Parabolic type subalgebras of $\CS(\vL)$}
\para{2.1.}
In [Sa] Sawada constructed a  subalgebra 
$\CS^0$ of $\CS(\vL)$
and showed that its quotient algebra $\ol\CS^0$ coincides with the cyclotomic 
$q$-Schur algebra associated to the modified Ariki-Koike algebra  
discussed in [SawS] under some condition on parameters (see
Introduction). $\CS^0$ is 
regarded, in some sense, a Borel type subalgebra of $\CS(\vL)$.  In this 
section, we extend his result to a more general situation, 
i.e, to the parabolic type subalgebras.
\para{2.2.}
Let $\vL$ and $\vL^+$ be as in Section 1.
We fix $\Bp = (r_1, \dots, r_g) \in \ZZ_{>0}^g$ such that 
$r = r_1 + \cdots + r_g$ for some $g$, and put 
$p_k = \sum_{i=1}^{k-1}r_i$ for $k = 1, \dots, g$ with $p_1 = 0$.
For each $\mu = (\mu^{(1)}, \dots, \mu^{(r)}) \in \vL$, 
put 
\begin{equation*}
\a_{\Bp}(\mu) = (n_1, \dots, n_g), \quad  
\Ba_{\Bp}(\mu) = (a_1, \dots, a_g),
\end{equation*}
where  
$n_k = \sum_{i=1}^{r_k}|\mu^{(p_k + i)}|$ and  
$a_k = \sum_{i=1}^{k-1}n_i$ for $k = 1, \dots, g$ with 
$a_1 = 0$.
By making use of $\Bp$, we express the $r$-compositions 
as the $g$-tuples of multi-compositions as follows; let 
$\mu = (\mu^{(1)}, \dots, \mu^{(r)}) \in \wt\CP_{n,r}$.   
We write $\mu$ as $(\mu^{[1]}, \dots, \mu^{[g]})$, where 
$\mu^{[k]} = (\mu^{(p_k+1)}, \dots, \mu^{(p_k + r_k)})$ 
is an $r_k$-composition of $n_k$.
Note that 
$\Ba_{\Bp}(\mu)$ (resp. $\a_{\Bp}(\mu)$) coincides with 
$\Ba(\mu)$ (resp. $\a(\mu)$) in 1.2 in the special case where
$\Bp = (1^r)$. 
\par
We define a partial order on $\ZZ_{\ge 0}^g$ by 
$\Ba = (a_1, \dots, a_g) \ge \Bb = (b_1, \dots, b_g)$ if
$a_k \ge b_k$ for $k = 1,\dots, g$.  We write $\Ba > \Bb$ if 
$\Ba \ge \Bb$ and $\Ba \ne \Bb$.  The following properties 
are easily verified.
\par\medskip\noindent
(2.2.1) \ Let $\mu, \nu \in \vL, \la \in \vL^+$.  Then we have 
\begin{enumerate}  
\item
$\Ba_{\Bp}(\mu) = \Ba_{\Bp}(\nu)$
if and oliy if $\a_{\Bp}(\mu) = \a_{\Bp}(\nu)$.
\item
If $\nu \trianglerighteq \mu$, then 
$\Ba_{\Bp}(\nu) \ge \Ba_{\Bp}(\mu)$.
In particular if $\CT_0(\la,\mu) \ne \emptyset$, 
then $\la \trianglerighteq \mu$ (cf. 1.3), and so
$\Ba_{\Bp}(\la) \ge \Ba_{\Bp}(\mu)$.  
\end{enumerate}
\par\medskip
For each $\la \in \vL^+, \mu \in \vL$, we define a set
$\CT_0^{\Bp}(\la, \mu)$ by $\CT_0(\la,\mu)$ if 
$\Ba_{\Bp}(\la) = \Ba_{\Bp}(\mu)$ and by the empty set otherwise.
Put $\CT_0^{\Bp}(\la) = \bigcup_{\mu \in \vL}\CT_0^{\Bp}(\la,\mu)$. 
\par\medskip\noindent
{\bf Example 2.3.}  
\addtocounter{thm}{1}
Let $n = 20, r = 5$ and take
$\mu = (21; 121; 32; 1^3; 41) \in \wt\CP_{20,5}$.
Let $\Bp = (2,2,1)$.  Then $\a_{\Bp}(\mu) = (7, 8, 5)$ and
$\a_{\Bp}(\mu) = (0, 7, 15)$.  
We have $\mu = (\mu^{[1]},\mu^{[2]}, \mu^{[3]})$ with 
$\mu^{[1]} = (21; 121), \mu^{[2]} = (32; 1^3), \mu^{[3]} = (41)$.
\para{2.4.}
Let $\ZC^{\Bp} = \ZC^{\Bp}(\vL)$ be the set of $\vf_{ST} \in \ZC(\vL)$ 
for $S \in \CT_0(\la,\mu), T \in \CT_0(\la,\nu)$, 
where $\mu,\nu \in \vL, \la \in \vL^+$ are taken subject to the condition
that $\Ba_{\Bp}(\la) > \Ba_{\Bp}(\mu)$ if  
  $\a_{\Bp}(\mu) \ne \a_{\Bp}(\nu)$.
We define an $R$-submodule $\CS^{\Bp} = \CS^{\Bp}(\vL)$ of 
$\CS(\vL)$ as the $R$-span of $\ZC^{\Bp}$.  We will see
that $\CS^{\Bp}$ is a subalgebra of $\CS(\vL)$ 
and that $\CS^{\Bp}$ turns out 
to be a standardly based algebra in the sense of Du-Rui [DR].
  Note that in the case where
$\Bp = (1^r)$, $\CS^{\Bp}$ coincides with $\CS^0$.
\par
First we note that the identity element $1_{\CS(\vL)}$ is 
contained in $\CS^{\Bp}$.
In fact, one can write $1_{\CS(\vL)} = \sum_{\mu \in \vL}\vf_{\mu}$, 
where $\vf_{\mu} \in \CS(\vL)$ is the identity map on $M^{\mu}$.  
Since $\vf_{\mu}$ is written as a linear combination of 
$\vf_{ST}$ with $S, T \in \CT_0(\la,\mu)$, we see that  
$1_{\CS(\vL)} \in \CS^{\Bp}$. 
\par
In order to relate $\CS^{\Bp}$ with the standardly based algebra,
we introduce a different kind of labeling for $\ZC^{\BP}$, following 
the idea of [Sa].  Let us define a subset $\vS^{\Bp}$ of 
$\vL^+\times \{ 0,1\}$ by 
\begin{equation*}
\begin{split}
\vS^{\Bp} = (\vL^+\times \{ 0,1\}) \backslash 
            &\{ (\la, 1) \mid \CT_0(\la,\mu) = \emptyset \\
               &\text{ for any } \mu \in \vL 
               \text{ such that } \Ba_{\Bp}(\la) > \Ba_{\Bp}(\mu)\}.
\end{split}
\end{equation*}
We define a partial order $\ge$ on $\vS^{\Bp}$ by 
$(\la_1,\ve_1) > (\la_2,\ve_2)$ if $\la_1 \triangleright \la_2$ 
or $\la_1 = \la_2$ and $\ve_1 > \ve_2$.
\par 
For each $\e = (\la,\ve) \in \vS^{\Bp}$, put
\begin{align*}
I^{\Bp}(\e) &= \begin{cases}
                 \CT_0^{\Bp}(\la) &\quad\text{ if } \ve = 0,  \\
                  \\                 
\displaystyle\bigcup_{\substack{\mu \in \vL \\
                  \Ba_{\Bp}(\la) > \Ba_{\Bp}(\mu)}}\CT_0(\la,\mu)
                                  &\quad\text{ if } \ve = 1,
               \end{cases}  
\\ \\
J^{\Bp}(\e) &= \begin{cases}
                \CT_0^{\Bp}(\la)  &\quad\text{ if }\ve = 0, \\
                \CT_0(\la)        &\quad\text{ if } \ve = 1. 
              \end{cases}
\end{align*}
Note that 
$I^{\Bp}(\e)$ and $J^{\Bp}(\e)$ are not empty. 
If we put, for $\e \in \vS^{\Bp}$,
\begin{equation*}
\ZC^{\Bp}(\e) = \{\vf_{ST} \mid S \in I^{\Bp}(\e), 
                                     T \in J^{\Bp}(\e)\},
\end{equation*}
we see easily that 
\begin{equation*}
\ZC^{\Bp} = \coprod_{\e \in \vS^{\Bp}}\ZC^{\Bp}(\e).
\end{equation*}
For each $\e \in \vS^{\Bp}$, we define a submodule 
$(\CS^{\Bp})^{\vee\e}$ of $\CS^{\Bp}$ as the $R$-span of 
$\vf_{ST}$, where $S \in I^{\Bp}(\e'), T \in J^{\Bp}(\e')$
for some $\e' \in \vS^{\Bp}$ such that $\e' > \e$.
\par
By using the cellular structure of $\CS(\vL)$, the following result
can be proved in a similar way as in [Sa, Lemma 2.4].
\newpage
\begin{lem} %%% Lemma 2.5
Take $\la_i \in \vL^+$, $\mu_i, \nu_i \in \vL$ for 
$i = 1,2$ such that $\nu_1 = \mu_2$.  
Then for $\vf_{S_iT_i} \in \ZC^{\Bp}$ 
with $S_i \in \CT_0(\la_i,\mu_i), T_i \in \CT_0(\la_i, \nu_i)$,
the followings hold.
\begin{equation*}
\begin{split}
&\vf_{S_1T_1}\cdot \vf_{S_2T_2} =  \\
&= \begin{cases}
    \displaystyle
       \sum_{\vf_{ST} \in \ZC^{\Bp}(\la_1,0)}r_{ST}\vf_{ST}
           + \sum_{\la \triangleright \la_1}\sum_{\vf_{ST} \in
            \ZC^{\Bp}(\la)}r_{ST}\vf_{ST}
                 &\quad\text{ if } \vf_{S_1T_1} \in \ZC^{\Bp}(\la_1,0), \\ 
      \\  
    \displaystyle\sum_{\la \trianglerighteq \la_1}
          \sum_{\vf_{ST}\in\ZC^{\Bp}(\la,1)}r_{ST}\vf_{ST} 
                 &\quad\text{ if } \vf_{S_1T_1} \in \ZC^{\Bp}(\la_1,1), \\
       \\     
     \displaystyle
       \sum_{\vf_{ST} \in \ZC^{\Bp}(\la_2,0)}r_{ST}\vf_{ST}
           + \sum_{\la \triangleright \la_2}\sum_{\vf_{ST} \in
            \ZC^{\Bp}(\la)}r_{ST}\vf_{ST}
                 &\quad\text{ if } \vf_{S_2T_2} \in \ZC^{\Bp}(\la_2,0), \\     
         \\     
      \displaystyle\sum_{\la \trianglerighteq \la_2}
          \sum_{\vf_{ST}\in\ZC^{\Bp}(\la,1)}r_{ST}\vf_{ST} 
                 &\quad\text{ if } \vf_{S_2T_2} \in \ZC^{\Bp}(\la_2,1), \\  
\end{cases}
\end{split}
\end{equation*}
where $r_{ST} \in R$ and 
$\ZC^{\Bp}(\la) = \ZC^{\Bp}(\la,0) \cup \ZC^{\Bp}(\la,1)$.
\end{lem}
The following theorem is a generalization of [Sa, Theorem 2.6].
The proof is done similarly by using Lemma 2.5.
%%%
\begin{thm}  %%% Theorem 2.6.
$\CS^{\Bp}$ is a subalgebra of $\CS(\vL)$ containing the identity 
element of $\CS(\vL)$.  Moreover, $\CS^{\Bp}$ turns out to be 
a standardly based algebra with the standard basis $\ZC^{\Bp}$ 
in the sense of [DR], i.e., the 
following holds; for any 
$\vf \in \CS^{\Bp}, \vf_{ST} \in \ZC^{\Bp}(\e)$, we have
\begin{align*}
\vf\cdot\vf_{ST} &\equiv \sum_{S' \in I^{\Bp}(\e)}
          f_{S'}\vf_{S'T}
                  \quad \mod (\CS^{\Bp})^{\vee\e}, \\
\vf_{ST}\cdot\vf &\equiv \sum_{T' \in J^{\Bp}(\e)}
          f'_{T'}\vf_{ST'}        
                   \quad \mod (\CS^{\Bp})^{\vee\e}
\end{align*}
with $f_{S'}, f'_{T'} \in R$, 
where in the first formula $f_{S'}$ depends on 
$(\vf, S, S')$ but does not depend on $T$, and in the second 
formula $f'_{T'}$ depends on $(\vf, T,T')$ but does not
depend on $S$.  
\end{thm}
\para{2.7.}
For each $\e \in \vS^{\Bp}$, let 
${}^{\dia}Z_{\Bp}^{\e}$ be an $R$-module with 
a basis $\{ \vf_S^{\e} \mid S \in I^{\Bp}(\e)\}$, and
$Z_{\Bp}^{\e}$ be an $R$-module with a basis 
$\{ \vf_T^{\e} \mid T \in J^{\Bp}(\e)\}$.
In view of Theorem 2.6, one can define actions of $\CS^{\Bp}$
on ${}^{\dia}Z_{\Bp}^{\e}$ and on $Z_{\Bp}^{\e}$ by
\begin{align*}
\vf\cdot\vf_S^{\e} = 
   \sum_{S' \in I^{\Bp}(\e)}f_{S'}\vf_{S'}^{\e}
\qquad(S \in I^{\Bp}(\e), \vf \in \CS^{\Bp}), \\
\vf_T^{\e}\cdot\vf = 
   \sum_{T' \in J^{\Bp}(\e)}f'_{T'}\vf_{T'}^{\e}
\qquad (T \in J^{\Bp}(\e), \vf \in \CS^{\Bp}),  
\end{align*}
where $f_{S'}, f'_{T'}$ are as in the theorem.
Then ${}^{\dia}Z_{\Bp}^{\e}$ (resp. $Z_{\Bp}^{\e}$)
has a structure of the left $\CS^{\Bp}$-module (resp. the right 
$\CS^{\Bp}$-module).
Moreover the theorem implies, for any 
$\vf_{UT}, \vf_{SV} \in \ZC^{\Bp}(\e)$,  that there exists 
$f_{TS} \in R$ (independent of the choice of $U,V$ ) such that 
\begin{equation*}
\vf_{UT}\vf_{SV} \equiv f_{TS}\vf_{UV} \mod (\CS^{\Bp})^{\vee\e}.
\end{equation*} 
We define a bilinear form 
$\b_{\e}: {}^{\dia}Z_{\Bp}^{\e} 
        \times Z_{\Bp}^{\e} \to R$
by $\b_{\e}(\vf_S^{\e}, \vf_T^{\e}) = f_{TS}$ 
for $S \in I^{\Bp}(\e), T \in J^{\Bp}(\e)$.
Put 
\begin{equation*}
\rad Z_{\Bp}^{\e} = \{ y \in Z_{\Bp}^{\e} \mid 
       \b_{\e}(x,y) = 0 \text{ for any } x \in {}^{\dia}Z_{\Bp}^{\e}\}. 
\end{equation*}
Then $\rad Z_{\Bp}^{\e}$ is an $\CS^{\Bp}$-submodule of $Z_{\Bp}^{\e}$
and we define the quotient module 
$L_{\Bp}^{\e} = Z_{\Bp}^{\e}/\rad Z_{\Bp}^{\e}$.
By the general theory of standardly based algebras (see [DR]), 
we obtain the following corollary, which is a strengthened form 
of [Sa, Proposition 3.7].
%%% 
\begin{cor}  %%%  Cor. 2.8
Assume that $R$ is a field.  Then 
\begin{enumerate}
\item
$L_{\Bp}^{\e}$ is an absolutely irreducible $\CS^{\Bp}$-module 
if it is non-zero.
\item
The set $\{ L_{\Bp}^{\e}\ne 0 \mid \e \in \vS^{\Bp}\}$ gives a complete
set of non-isomorphic irreducible right $\CS^{\Bp}$-modules.
\end{enumerate}
\end{cor}
\remarks{2.9.} \ (i)
In [Sa], only the case $Z_{\Bp}^{(\la,0)}$
is discussed (for $\Bp = (1^r)$).
In that case (for arbitrary $\Bp$), we have 
the following description on the basis of $Z_{\Bp}^{(\la,0)}$ 
as in the case of the Weyl module $W^{\la}$.
For each $\la \in \vL^+$, $\vf_{\la} = \vf_{T^{\la}T^{\la}}$ is 
contained in $\ZC^{\Bp}(\la,0)$.  We consider the $\CS^{\Bp}$-submodule
$W_{\Bp}^{\la}$ of $\CS^{\Bp}/(\CS^{\Bp})^{\vee(\la,0)}$ 
generated by the image of $\vf_{\la}$.
Since $T^{\la} \in I^{\Bp}(\la,0)$, we see that 
$\vf_{T^{\la}T} \in \ZC^{\Bp}(\la,0)$ for any $T \in J^{\Bp}(\la,0)$.
We denote by $\vf'_T$ the image of $\vf_{T^{\la}T}$ on 
$\CS^{\Bp}/(\CS^{\Bp})^{\vee(\la,0)}$.  Then one can check that 
$\vf'_T \in W_{\Bp}^{\la}$ and that  
the map $\vf_T \to \vf'_T$ gives an isomorphism 
$Z_{\Bp}^{(\la,0)} \to W_{\Bp}^{\la}$ of $\CS^{\Bp}$-modules. 
In particular, we see that $Z_{\Bp}^{(\la,0)}$ is generated by
$\vf_{T^{\la}}^{(\la,0)}$ as an $\CS^{\Bp}$-module.
\par
However, the above argument can not be applied to $Z_{\Bp}^{(\la,1)}$
since $\vf_{T^{\la}T} \notin \CS^{\Bp}$ for 
$T \in J^{\Bp}(\la,1) \backslash J^{\Bp}(\la,0)$.
It is not known whether $Z_{\Bp}^{(\la,1)}$ is generated by one element 
as an $\CS^{\Bp}$-module.
\par 
(ii) \ For any $\la \in \vL^+$, we have $L_{\Bp}^{(\la,0)} \ne 0$.
In fact, since $T^{\la} \in I^{\Bp}(\la,0) \cap J^{\Bp}(\la,0)$, 
we have $f_{T^{\la}T^{\la}} = 1$.  This implies that 
$\b_{(\la,0)}(\vf_{T^{\la}}^{(\la,0)}, \vf_{T^{\la}}^{(\la,0)}) = 1$
and we see that $\rad Z_{\Bp}^{(\la,0)} \ne Z_{\Bp}^{(\la,0)}$.
\par
This argument cannot be applied to $Z_{\Bp}^{(\la,1)}$ since
$T^{\la} \notin I^{\Bp}(\la,1)$ and so 
$\vf_{T^{\la}}^{(\la,1)} \notin {}^{\dia}Z_{\Bp}^{(\la,1)}$.
It is not known when $L_{\Bp}^{(\la,1)} \ne 0$.
\para{2.10.}
Recall that $\vf \mapsto \vf^*$ be the anti-automorphism on 
$\CS(\vL)$ defined by $\vf_{ST} \mapsto \vf_{TS}$, related to 
the cellular structure.  Let  
$\CS^{\Bp*} = (\CS^{\Bp})^*$ be the image of $\CS^{\Bp}$ under 
the map $*$.
Then $\CS^{\Bp*}$ is a subalgebra of $\CS(\vL)$, and 
it is easy to check that $\CS^{\Bp*}$ is a standardly
based algebra with the standard basis 
$\ZC^{\Bp*} = \coprod_{\e \in \vS^{\Bp}}\ZC^{\Bp}(\e)^*$, 
where 
\begin{equation*}
\ZC^{\Bp}(\e)^* = \{ \vf_{ST}\in \ZC(\vL) 
             \mid S \in J^{\Bp}(\e), T \in I^{\Bp}(\e)\}.
\end{equation*}
In a similar way as in [Sa, Proposition 3.2], one can show 
the following result.
%%%
\begin{prop}  %%% Prop. 2.11
We have
$\CS(\vL) = \CS^{\Bp}\cdot \CS^{\Bp*}$.
\end{prop}
\para{2.12.}  
Let $\wh\CS^{\Bp}$ be the $R$-submodule of $\CS^{\Bp}$
spanned by 
\begin{equation*}
\wh\ZC^{\Bp} = \ZC^{\Bp} \backslash 
                       \{ \vf_{ST} \mid S, T \in \CT_0^{\Bp}(\la)
                            \text{ for some }\la \in \vL^+ \}.
\end{equation*}
Then by the second and the fourth formulas in Lemma 2.5, 
$\wh\CS^{\Bp}$ turns out to be a two-sided ideal of 
$\CS^{\Bp}$.  We denote by $\ol\CS^{\Bp} = \ol\CS^{\Bp}(\vL)$ the 
quotient algebra $\CS^{\Bp}/\wh\CS^{\Bp}$.
Let $\pi: \CS^{\Bp} \to \ol\CS^{\Bp}$ be the
natural projection, and put $\ol\vf = \pi(\vf)$ for 
$\vf \in \CS^{\Bp}$.  It is easy to see that $\ol\CS^{\Bp}$
is an $R$-free module with the basis
\begin{equation*}
\ol\ZC^{\Bp} = \{ \ol\vf_{ST} \mid S, T \in \CT_0^{\Bp}(\la)
                     \text{ for } \la \in \vL^+\}.
\end{equation*}
Similarly, one can define a quotient algebra 
$\ol\CS^{\Bp*} = \CS^{\Bp*}/\wh\CS^{\Bp*}$, where 
$\wh\CS^{\Bp*} = (\wh\CS^{\Bp})^*$ is a two-sided ideal of 
$\CS^{\Bp*}$.  Let $\pi'$ be the natural projection 
$\CS^{\Bp*} \to \ol\CS^{\Bp*}$, and put $\ol\vf' = \pi'(\vf)$
for $\vf \in \CS^{\Bp*}$.  Then $\ol\CS^{\Bp*}$ has an $R$-free
basis 
$\ol\ZC^{\Bp*} = \{ \ol\vf'_{ST} \mid 
    S, T \in \CT_0^{\Bp}(\la), \la \in \vL^+\}$.
It is clear that $\ol\vf_{ST} \mapsto \ol\vf'_{ST}$
gives an isomorphism $\ol\CS^{\Bp} \to \ol\CS^{\Bp*}$ of
$R$-algebras.  
On the other hand, the anti-algebra isomorphism 
$\CS^{\Bp} \to \CS^{\Bp*}$ induces an anti-algebra isomorphism
$\ol\CS^{\Bp} \to \ol\CS^{\Bp*}, \ol\vf_{ST} \mapsto \ol\vf'_{TS}$.
It follows that the map $\ol\vf_{ST} \mapsto \ol\vf_{TS}$ induces
an anti-algebra automorphism $*$ on $\ol\CS^{\Bp}$.
Thus we have the following theorem (cf. [Sa, Theorem 4.8]).
Note that the second assertion is obtained from the cellular structure of
$\CS(\vL)$.
%%%
\begin{thm}   %%% Theorem 2.13.
$\ol\CS^{\Bp}$ is a cellular algebra with a cellular basis 
$\ol\ZC^{\Bp}$, i.e., the following property holds;
\begin{enumerate}
\item
$\ol\vf_{ST} \mapsto (\ol\vf_{ST})^* = \ol\vf_{TS}$ gives 
an anti-algebra automorphism $*$ on $\ol\CS^{\Bp}$.
\item
Let $(\ol\CS^{\Bp})^{\vee\la}$ be the $R$-submodule of $\ol\CS^{\Bp}$
spanned by $\ol\vf_{ST}$ such that $S,T \in \CT_0^{\Bp}(\la')$ with 
$\la' \triangleright \la$. Then for any 
$\la \in \vL^+, S, T \in \CT_0^{\Bp}(\la), \ol\vf \in \ol\CS^{\Bp}$, 
\begin{equation*}
\ol\vf_{ST}\cdot\ol\vf \equiv \sum_{T' \in \CT_0^{\Bp}(\la)}
                r_{T'}\ol\vf_{ST'} \quad\mod (\ol\CS^{\Bp})^{\vee\la},
\end{equation*}
where $r_{T'} \in R$ depends on $\la, T, \ol\vf$, 
but does not depend on $S$.
\end{enumerate}
\end{thm}
\para{2.14.}
We apply the general theory of cellular algebras to $\ol\CS^{\Bp}$.
For each $\la \in \vL^+$, $(\ol\CS^{\Bp})^{\vee\la}$
is a two-sided ideal of $\ol\CS^{\Bp}$, and we define the Weyl module
$\ol Z_{\Bp}^{\la}$ as the $\ol\CS^{\Bp}$-submodule of the right 
$\ol\CS^{\Bp}$-module $\ol\CS^{\Bp}/(\ol\CS^{\Bp})^{\vee\la}$
generated by $\ol\vf_{T^{\la}T^{\la}} + (\ol\CS^{\Bp})^{\vee\la}$.
Let $\ol\vf_T$ be the image of $\ol\vf_{T^{\la}T}$ on 
$\ol\CS^{\Bp}/(\ol\CS^{\Bp})^{\vee\la}$.  Then the set 
$\{ \ol\vf_T \mid T \in \CT_0^{\Bp}(\la)\}$ gives a basis of 
$\ol Z_{\Bp}^{\la}$.  The symmetric bilinear form 
$\lp\ ,\ \rp_{\Bp}: \ol Z^{\la}_{\Bp} \times \ol Z^{\la}_{\Bp} \to R$
is defined by the equation
\begin{equation*}
\lp\ol\vf_S, \ol\vf_T\rp_{\Bp}\ol\vf_{T^{\la}T^{\la}} 
       \equiv \ol\vf_{T^{\la}S}\ol\vf_{TT^{\la}}
         \quad \mod (\ol\CS^{\Bp})^{\vee\la}.
\end{equation*}
Then the radical $\rad \ol Z_{\Bp}^{\la}$ of $\ol Z_{\Bp}^{\la}$
with respect to this form is an $\ol\CS^{\Bp}$-submodule of 
$\ol Z_{\Bp}^{\la}$, and we define an $\ol\CS^{\Bp}$-module 
$\ol L_{\Bp}^{\la}$ by 
$\ol L_{\Bp}^{\la} = \ol Z_{\Bp}^{\la}/\rad \ol Z_{\Bp}^{\la}$.
Since $\lp\ol\vf_{T^{\la}}, \ol\vf_{T^{\la}}\rp_{\Bp} = 1$, we 
see that $\ol L_{\Bp}^{\la} \ne 0$ for any $\la \in \vL^+$.
By the general theory of cellular algebras, we have
%%%
\begin{cor}  %%% Cor 2.15
Suppose that $R$ is a field.  Then,
for any $\la \in \vL^+$, $\ol L_{\Bp}^{\la}$ is an absolutely
irreducible $\ol\CS^{\Bp}$-module, and the set 
$\{ \ol L_{\Bp}^{\la} \mid \la \in \vL^+ \}$ gives a complete
set of non-isomorphic $\ol\CS^{\Bp}$-modules.
\end{cor}
\section{Decomposition numbers for $\CS(\vL), \CS^{\Bp}$ and
$\ol\CS^{\Bp}$}
 \para{3.1.}
By the discussion in the previous section, we have the following diagram.
\begin{equation*}
\begin{CD}
\CS^{\Bp} @> \io >> \CS(\vL)  \\
@V\pi VV                          \\
\ol\CS^{\Bp} 
\end{CD}
\end{equation*}
where $\io$ is the inclusion map, and $\pi$ is the natural 
surjective map.  
We have constructed the Weyl modules 
$W^{\la}, Z_{\Bp}^{\e}$ and $\ol Z^{\la}_{\Bp}$ for 
$\CS(\vL), \CS^{\Bp}$ and $\ol\CS^{\Bp}$,
and assuming that $R$ is a field, the irreducible modules
$L^{\mu}, L_{\Bp}^{\e'}, \ol L_{\Bp}^{\mu}$, respectively
for $\la, \mu \in \vL^+, \e, \e' \in \vS^{\Bp}$.
We consider the decomposition numbers 
\begin{equation*}
[W^{\la}: L^{\mu}]_{\CS(\vL)}, \quad 
[Z_{\Bp}^{\e} : L_{\Bp}^{\e'}]_{\CS^{\Bp}}, \quad
[\ol Z_{\Bp}^{\la} : \ol L_{\Bp}^{\mu}]_{\ol\CS^{\Bp}}
\end{equation*}
for $\CS(\vL), \CS^{\Bp}$ and $\ol\CS^{\Bp}$.
By using the above maps, we shall discuss the 
relationship among these decomposition numbers.
\par
First we consider the relation between $\CS^{\Bp}$ and 
$\ol\CS^{\Bp}$.  We regard an $\ol\CS^{\Bp}$-module
as an $\CS^{\Bp}$-module through the map $\pi$.   
The following lemma is easily verified if we notice that 
$\pi((\CS^{\Bp})^{\vee(\la,0)}) = \ol\CS^{\vee\la}$
and that 
$\b_{(\la,0)}(\vf_S^{(\la,0)}, \vf_T^{(\la,0)})
   = \lp\ol\vf_S, \ol\vf_T\rp_{\Bp}$ for $S, T \in \CT_0^{\Bp}(\la)$.
%%%
\begin{lem}  %%% Lemma 3.2
\begin{enumerate}
\item
For any $\la \in \vL^+$, the map $\vf_T^{(\la,0)} \mapsto \ol\vf_T$
($T \in \CT_0^{\Bp}(\la)$) gives an isomorphism 
$Z_{\Bp}^{(\la,0)} \isom \ol Z_{\Bp}^{\la}$ of $\CS^{\Bp}$-modules.
\item
Assume that $R$ is a field.  Then the above map induces 
an isomorphism $L_{\Bp}^{(\la,0)} \isom \ol L_{\Bp}^{\la}$ of
$\CS^{\Bp}$-modules.
\end{enumerate}
\end{lem}
The following proposition is proved in a similar way as 
in [Sa, Theorem 4.15] by taking the lemma into account.
%%%
\begin{prop}  %%% Prop 3.3
Assume that $R$ is a field.  Then
\begin{enumerate}
\item
The composition factors of $Z_{\Bp}^{(\la,0)}$ are isomorphic to
$L_{\Bp}^{(\mu,0)}$ for some $\mu \in \vL^+$  such that 
$\la \triangleright \mu$.
\item
For any $\la,\mu \in \vL^+$, we have
$[\ol Z_{\Bp}^{\la} : \ol L_{\Bp}^{\mu}]_{\ol\CS^{\Bp}}
  = [Z_{\Bp}^{(\la,0)} : L_{\Bp}^{(\mu,0)}]_{\CS^{\Bp}}$.
\item
For $\la, \mu \in \vL^+$ such that 
$\a_{\Bp}(\la) \ne \a_{\Bp}(\mu)$, we have
$[\ol Z_{\Bp}^{\la} : \ol L_{\Bp}^{\mu}]_{\ol\CS^{\Bp}} 
  = 0. $
\end{enumerate}
\end{prop}
\para{3.4.}
Next we consider the relation between $\CS(\vL)$ and $\CS^{\Bp}$.
Since $\CS^{\Bp}$ is a subalgebra of $\CS(\vL)$, we regard
an $\CS(\vL)$-module as an $\CS^{\Bp}$-module by restriction.
Recall that $J^{\Bp}(\la,0) = \CT_0^{\Bp}(\la), 
J^{\Bp}(\la,1) = \CT_0(\la)$ for $\la \in \vL^+$.  Thus the 
basis of $Z_{\Bp}^{(\la,0)}$ is 
$\{ \vf_T^{(\la,0)} \mid T \in \CT_0^{\Bp}(\la)\}$, the basis of 
$Z_{\Bp}^{(\la,1)}$ is $\{ \vf_T^{(\la,1)} \mid T \in \CT_0(\la)\}$, 
and the basis of 
$W^{\la}$ is $\{ \vf_T \mid T \in \CT_0(\la)\}$, respectively. 
The following result is implicit in [Sa].
%%%
\begin{lem} %%% Lemma 3.5
For each $\la \in \vL^+$, the followings hold.
\begin{enumerate}
\item
The map $\vf_T^{(\la,0)} \mapsto \vf_T^{(\la,1)} \
      (T \in \CT_0^{\Bp}(\la))$ gives an injective homomorphism 
$Z_{\Bp}^{(\la,0)} \to Z_{\Bp}^{(\la,1)}$ of $\CS^{\Bp}$-modules.
\item
The map $\vf_T^{(\la,1)} \mapsto \vf_T$ ($T \in \CT_0(\la)$) 
gives an isomorphism 
$Z_{\Bp}^{(\la,1)} \isom W^{\la}$ of $\CS^{\Bp}$-modules.
\end{enumerate}
\end{lem}
\begin{proof}
Take $\vf_{ST} \in \ZC(\vL)$ ($S,T \in \CT_0(\la)$), and 
$\vf \in \CS^{\Bp}$.  By the property of the cellular algebra 
$\CS(\vL)$, we have
\begin{equation*}
\vf_{ST}\cdot\vf \equiv \sum_{T' \in \CT_0(\la)}r_{T'}\vf_{ST'}
         \quad\mod \CS(\vL)^{\vee\la}.
\end{equation*}
Since $\vf_{ST} \in \CS^{\Bp}$ and 
$\CS(\vL)^{\vee\la} \cap \CS^{\Bp} = (\CS^{\Bp})^{\vee(\la,1)}$, 
the congruence relation by $\CS(\vL)^{\vee\la}$ in the above formula 
can be replaced by $(\CS^{\Bp})^{\vee(\la,1)}$.
In particular, for $\vf_{ST} \in \ZC^{\Bp}(\la,1), \vf \in \CS^{\Bp}$,
we have
\begin{equation*}
\tag{3.5.1}
\vf_{ST}\cdot\vf \equiv \sum_{T' \in \CT_0(\la)}r_{T'}\vf_{ST'}
\quad\mod (\CS^{\Bp})^{\vee(\la,1)}.
\end{equation*}
On the other hand by the second formula in Theorem 2.6 we have, 
for $\vf_{ST} \in \ZC^{\Bp}(\la,0), \vf \in \CS^{\Bp}$, 
\begin{equation*}
\vf_{ST}\cdot\vf \equiv \sum_{T' \in \CT_0^{\Bp}(\la)}f_{T'}\vf_{ST'}
   \quad\mod (\CS^{\Bp})^{\vee(\la,0)}.
\end{equation*}
 But the first formula in Lemma 2.5 shows that 
$\vf_{ST'} \in \ZC^{\Bp}(\la,1)$ does not appear in the 
expression of $\vf_{ST}\cdot\vf$ except 
$\vf_{ST'} \in \ZC^{\Bp}(\la,0)$.  It follows that 
the congruence relation $(\CS^{\Bp})^{\vee(\la,0)}$ in the above formula
can be replaced by $(\CS^{\Bp})^{\vee(\la,1)}$.  Thus we have, for
$\vf_{ST} \in \ZC^{\Bp}(\la,0), \vf \in \CS^{\Bp}$
\begin{equation*}
\tag{3.5.2}
\vf_{ST}\cdot\vf \equiv \sum_{T' \in \CT_0^{\Bp}(\la)}f_{T'}\vf_{ST'}
   \quad\mod (\CS^{\Bp})^{\vee(\la,1)}.
\end{equation*}
We now prove (i).  For $T \in \CT_0^{\Bp}(\la), \vf \in \CS^{\Bp}$, one can
write as 
\begin{align*}
\vf_T^{(\la,0)}\cdot\vf &= 
   \sum_{T' \in \CT_0^{\Bp}(\la)}g_{T'}\vf_{T'}^{(\la,0)}, \\
\vf_T^{(\la,1)}\cdot\vf &= 
   \sum_{T' \in \CT_0(\la)}g'_{T'}\vf_{T'}^{(\la,1)}.
\end{align*}
By the definition of Weyl modules, we see that $g_{T'} = f_{T'}$
and $g'_{T'} = r_{T'}$.  Thus by comparing (3.5.1) and (3.5.2), we have
\begin{equation*}
g'_{T'} = \begin{cases}
            g_{T'} &\quad\text{ if } T' \in \CT_0^{\Bp}(\la), \\
            0      &\quad\text{ otherwise.}
          \end{cases}
\end{equation*}
This proves (i).  The assertion (ii) is proved in a similar way.
\end{proof}
\par
The following proposition can be proved in a similar way as in
[Sa, Theorem 3.3].
%%%
\begin{prop}  %%% Prop. 3.6
For each $\la \in \vL^+$, there exists an isomorphism of 
$\CS(\vL)$-modules
\begin{equation*}
Z_{\Bp}^{(\la,0)}\otimes_{\CS^{\Bp}}\CS(\vL) \isom W^{\la}
\end{equation*}
which maps $\vf_{T^{\la}}^{(\la,0)}\p\otimes\vf$ to 
$\vf_{T^{\la}}\p\vf$ for $\p \in \CS^{\Bp}, \vf \in \CS(\vL)$.
\end{prop}
\para{3.7.}
By Lemma 3.5, the map $\vf_T^{(\la,0)} \mapsto \vf_T$ gives 
an injective homomorphism $f_{\la} : Z_{\Bp}^{(\la,0)} \to W^{\la}$ of 
$\CS^{\Bp}$-modules.  By this map we regard $Z_{\Bp}^{(\la,0)}$
as an $\CS^{\Bp}$-submodule of $W^{\la}$.  We have the following 
lemma.
%%%
\begin{lem}  %%% Lemma 3.8
Assume that $\la \in \vL^+$.
\begin{enumerate}
\item
Let $M$ be an $\CS^{\Bp}$-submodule of $Z_{\Bp}^{(\la,0)}$, 
and $\wt M$ be the $\CS(\vL)$-submodule of $W^{\la}$ generated by
$M$.
Then $\wt M \cap Z_{\Bp}^{(\la,0)} = M$.
\item
Let $M_1 \subsetneq M_2$ be $\CS^{\Bp}$-submodules of 
$Z_{\Bp}^{(\la,0)}$. Let $\io_i$ be the inclusion map 
$M_i \to Z_{\Bp}^{(\la,0)}$, and $\io_i\otimes\Id$ be the induced map
$M_i\otimes_{\CS^{\Bp}}\CS(\vL) \to 
    Z_{\Bp}^{(\la,0)}\otimes_{\CS^{\Bp}}\CS(\vL)$ for $i = 1,2$.
Then we have
$\Im (\io_1\otimes\Id) \subsetneq \Im (\io_2\otimes\Id)$.
\end{enumerate}
\end{lem}
\begin{proof}
We show (i).  Take $x \in \wt M \cap Z_{\Bp}^{(\la,0)}$.
We write $x = \sum_{T \in \CT_0^{\Bp}(\la)}r_T\vf_T^{(\la,0)}$.
Since $x \in \wt M$, one can write 
$x = \sum_{i}y_i\p_i$ with $\p_i \in \CS(\vL)$, 
$y_i = \sum_{T \in \CT_0^{\Bp}(\la)}r_{T,i}\vf_T^{(\la,0)} 
       \in Z_{\Bp}^{(\la,0)}$.
Hence we have a relation as elements in $W^{\la}$
\begin{equation*}
\sum_{T \in \CT_0^{\Bp}(\la)}r_T\vf_T = \sum_{i}
       \sum_{T \in \CT_0^{\Bp}(\la)}r_{T,i}\vf_T\p_i.
\end{equation*}
This means that 
\begin{equation*}
\sum_{T \in \CT_0^{\Bp}(\la)}r_T\vf_{T^{\la}T}
 \equiv \sum_i\sum_{T \in \CT_0^{\Bp}(\la)}r_{T,i}\vf_{T^{\la}T}\p_i
\quad \mod \CS(\vL)^{\vee\la}.
\end{equation*}
Put $\a = \a_{\Bp}(\la)$.  Take $\nu \in \vL$ such that 
$\a_{\Bp}(\nu) = \a$ and multiply $\vf_{\nu}$ on both sides of 
the above equation.
Note that $\vf_{\nu} \in \CS^{\Bp}$ is a projection from $M$ to $M^{\nu}$,  
and we have 
$\CS^{\Bp}\cap \CS(\vL)^{\vee\la} = (\CS^{\Bp})^{\vee(\la,1)} 
         \subset (\CS^{\Bp})^{\vee(\la,0)}$.
It follows that 
\begin{equation*}
\sum_{T \in \CT_0(\la,\nu)}r_T\vf_{T^{\la}T} \equiv
      \sum_i\sum_{T \in \CT_0^{\Bp}(\la)}r_{T,i}\vf_{T^{\la}T}\p_i\vf_{\nu}
   \quad\mod (\CS^{\Bp})^{\vee(\la,0)}.
\end{equation*}
Since this holds for any $\nu \in \vL$ such that $\a_{\Bp}(\nu) = \a$, 
we have
\begin{equation*}
\tag{3.8.1}
\sum_{T \in \CT_0^{\Bp}(\la)}r_T\vf_{T^{\la}T} \equiv
   \sum_{\substack{ \nu \in \vL \\ \a_{\Bp}(\nu) = \a}}
    \sum_i\sum_{T \in \CT_0^{\Bp}(\la)}
           r_{T,i}\vf_{T^{\la}T}\p_i\vf_{\nu}  
   \quad\mod (\CS^{\Bp})^{\vee(\la,0)}.
\end{equation*}
Put 
$\vf_{\a} = \sum_{\nu}\vf_{\nu}$, where $\nu$ runs over all the 
elements in $\vL$ such that $\a_{\Bp}(\nu) = \a$.
Since $\vf_{\a}$ is the projection from $M$ onto 
$M^{\a} = \bigoplus_{\nu}M^{\nu}$, 
we see that $\vf_{T^{\la}T}\vf_{\a} = \vf_{T^{\la}T}$ for any 
$T \in \CT_0^{\Bp}(\la)$.  Moreover we note that 
$\vf_{\a}\p_i\vf_{\nu} \in \CS^{\Bp}$ since it is contained 
in $\Hom_{\CH}(M^{\nu}, M^{\a})$.  
It follows that 
\begin{equation*}
\vf_{T^{\la}T}\p_i\vf_{\nu} = 
   \vf_{T^{\la}T}(\vf_{\a}\p_i\vf_{\nu}) \in \CS^{\Bp}
\end{equation*}
for $T \in \CT_0^{\Bp}(\la)$.
Thus one can rewrite (3.8.1) as a relation on $Z_{\Bp}^{(\la,0)}$
as 
\begin{equation*}
x = \sum_{T \in \CT_0^{\Bp}(\la)}r_T\vf_T^{(\la,0)} = 
   \sum_{\substack{ \nu \in \vL \\ \a_{\Bp}(\nu) = \a}}
    \sum_i\sum_{T \in \CT_0^{\Bp}(\la)}
           r_{T,i}\vf_T^{(\la,0)}(\vf_{\a}\p_i\vf_{\nu}).  
\end{equation*}
This shows that 
$x = \sum_{\nu}\sum_iy_i(\vf_{\a}\p_i\vf_{\nu}) \in M$ as asserted.
\par
Next we show (ii).  Under the embedding 
$Z_{\Bp}^{\la} \hra W^{\la}$, Proposition 3.6 implies that 
$\Im (\io_i\otimes\Id) = \wt M_i$.  Take $x \in M_2 \backslash M_1$.
Suppose that $\wt M_1 = \wt M_2$.  Then 
$x \in \wt M_1 \cap Z_{\Bp}^{(\la,0)} = M_1$ by (i).  This is 
a contradiction. 
\end{proof}
\par
By making use of Lemma 3.8, we show the following lemma.
%%%
\begin{lem} %%% Lemma 3.9
Assume that $R$ is a field.  Then for each $\la \in \vL^+$, 
there exists a unique maximal $\CS(\vL)$-submodule $N^{\la}$ of 
$L_{\Bp}^{(\la,0)}\otimes_{\CS^{\Bp}}\CS(\vL)$ such that
\begin{equation*}
L_{\Bp}^{(\la,0)}\otimes_{\CS^{\Bp}}\CS(\vL)/N^{\la} \simeq L^{\la}.
\end{equation*}
\begin{proof}
Since $L_{\Bp}^{(\la,0)} \simeq 
  Z_{\Bp}^{(\la,0)}/\rad Z_{\Bp}^{(\la,0)}$, 
we have a surjective homomorphism 
\begin{equation*}
Z_{\Bp}^{(\la,0)}\otimes_{\CS^{\Bp}}\CS(\vL) \to 
L_{\Bp}^{(\la,0)}\otimes_{\CS^{\Bp}}\CS(\vL)
\end{equation*}
as $\CS(\vL)$-modules.
Since $L_{\Bp}^{(\la,0)} \ne 0$, we have 
$L_{\Bp}^{(\la,0)}\otimes_{\CS^{\Bp}}\CS(\vL) \ne 0$ by 
Lemma 3.8, and the kernel of this map is a proper $\CS(\vL)$-submodule 
of $Z_{\Bp}^{(\la,0)}\otimes_{\CS^{\Bp}}\CS(\vL)$.  But 
$Z_{\Bp}^{(\la,0)}\otimes_{\CS^{\Bp}}\CS(\vL)$ 
is isomorphic to $W^{\la}$ by Proposition 3.6, and 
$L^{\la} \simeq W^{\la}/\rad W^{\la}$.  Since $\rad W^{\la}$ is 
the unique maximal submodule of $W^{\la}$, we see that there 
exists a surjective homomorphism 
$L_{\Bp}^{(\la,0)}\otimes_{\CS^{\Bp}}\CS(\vL) \to L^{\la}$ of 
$\CS(\vL)$-modules. It is clear that $N^{\la}$ is the unique maximal 
submodule since it is a quotient of $\rad W^{\la}$.
\end{proof}
\end{lem}
%%%
\begin{lem} %%% Lemma 3.10
Assume that $R$ is a field.  For each $\la \in \vL^+$, 
the $\CS^{\Bp}$-module $L^{\la}$ contains $L_{\Bp}^{(\la,0)}$ as a
submodule.
\end{lem}
\begin{proof}
By definition, we have 
$\b_{\la}(\vf_S^{(\la,0)},\vf_T^{(\la,0)}) = \lp \vf_S, \vf_T\rp$
for any $S, T \in \CT_0^{\Bp}(\la)$.
Moreover, one can check that  $\lp\vf_S, \vf_T\rp = 0$ for 
$S \in \CT_0(\la) \backslash \CT_0^{\Bp}(\la)$, $T \in \CT_0^{\Bp}(\la)$.
It follows that 
$f_{\la}(\rad Z_{\Bp}^{(\la,0)}) \subset \rad W^{\la}$, where 
$f_{\la}: Z_{\Bp}^{(\la,0)} \hra W^{\la}$ is the injective map 
given in 3.7.  Then $f_{\la}$ induces a homomorphism 
$\bar f_{\la} : L_{\Bp}^{(\la,0)} \to L^{\la}$ of $\CS^{\Bp}$-modules. 
Since $f_{\la}(\vf_{T^{\la}}^{(\la,0)}) = 
    \vf_{T^{\la}} \notin \rad W^{\la}$,
$\bar f_{\la}$ is a non-zero map.  Since $L_{\Bp}^{(\la,0)}$ 
is an irreducible $\CS^{\Bp}$-module, $\bar f_{\la}$ is injective.
This proves the lemma.  
\end{proof}
The following two results are generalizations of 
[Sa, Theorem 5.6, Theorem 5.7].
%%%%
\begin{prop} %%% Prop. 3.11
Assume that $R$ is a field.  Then for $\la, \mu \in \vL^+$, 
\begin{equation*}
[Z_{\Bp}^{(\la,0)}: L_{\Bp}^{(\mu,0)}]_{\CS^{\Bp}} 
              \le [W^{\la} : L^{\mu}]_{\CS(\vL)}. 
\end{equation*}
\end{prop}
\begin{proof}
We consider a composition series of $Z_{\Bp}^{(\la,0)}$ as 
an $\CS^{\Bp}$-module
\begin{equation*}
0 = M_0 \subsetneq M_1 \subsetneq \cdots 
               \subsetneq M_k = Z_{\Bp}^{(\la,0)}
\end{equation*}
such that $M_j/M_{j-1} \simeq L_{\Bp}^{(\mu_j,0)}$.
Let $i_j : M_j \hra Z_{\Bp}^{(\la,0)}$ be the inclusion map
and $\io_j\otimes\Id: M_j\otimes_{\CS^{\Bp}}\CS(\vL) \to 
Z_{\Bp}^{(\la,0)}\otimes_{\CS^{\Bp}}\CS(\vL)$ be the
induced map.  Put $\CM_j = \Im(\io_j\otimes\Id)$.
We have a filtration 
\begin{equation*}
0 = \CM_0 \subsetneq \CM_1 \subsetneq \cdots \subsetneq \CM_k = 
     Z_{\Bp}^{(\la,0)}\otimes_{\CS^{\Bp}}\CS(\vL) \simeq W^{\la}
\end{equation*}
of $\CS(\vL)$-submodules of $W^{\la}$ by Proposition 3.6 and Lemma 3.8.
In order to prove the proposition, it is enough to show that 
$L^{\mu_j}$ occurs in the composition series of $\CM_j/\CM_{j-1}$
for each $j$.
Since $M_j/M_{j-1} \simeq L_{\Bp}^{(\mu_j,0)}$, we have the following 
diagram of $\CS(\vL)$-modules
\begin{equation*}
\begin{CD}
 @.  M_{j-1}\otimes\CS(\vL) @>>> M_j\otimes\CS(\vL)
               @>>> L_{\Bp}^{(\mu_j,0)}\otimes\CS(\vL)
               @>>> 0  \\
@.   @VVV    @VVV  \\
0 @>>> \CM_{j-1} @>>> \CM_j @>>> \CM_j/\CM_{j-1} @>>> 0
\end{CD}
\end{equation*}
where the vertical maps are surjective.
Thus we obtain a surjective homomorphism
$L_{\Bp}^{(\la,0)}\otimes_{\CS^{\Bp}}\CS(\vL) \to \CM_j/\CM_{j-1}$.
On the other hand, by Lemma 3.9, we have a surjective homomorphism
$L_{\Bp}^{(\la,0)}\otimes_{\CS^{\Bp}}\CS(\vL) \to L^{\mu}$, 
whose kernel $N^{\la}$
is the unique maximal submodule. 
This implies that we have a surjectie homomorphism
$\CM_j/\CM_{j-1} \to L^{\mu_j}$. Hence $L^{\m_j}$ occurs 
in the composition series of $\CM_j/\CM_{j-1}$, and the proposition
is proved.  
\end{proof}
%%%
\begin{prop} %%% Prop. 3.12
Assume that $R$ is a field. Then for any $\la, \mu \in \vL^+$
such that $\a_{\Bp}(\la) = \a_{\Bp}(\mu)$, we have
\begin{equation*}
[Z_{\Bp}^{(\la,0)} : L_{\Bp}^{(\mu,0)}]_{\CS^{\Bp}}
   \ge [W^{\la} : L^{\mu}]_{\CS(\vL)}.
\end{equation*}
\end{prop}
\begin{proof}
Consider a composition series of $W^{\la}$ as an $\CS(\vL)$-module
\begin{equation*}
0 = W_0 \subsetneq W_1 \subsetneq \cdots \subsetneq W_k = W^{\la}
\end{equation*}
such that $W_j/W_{j-1}\simeq L^{\mu_j}$ for some $\mu_j \in \vL^+$.
We consider this as a filtration of $W^{\la}$ as $\CS^{\Bp}$-modules.
Since $L^{\mu_j}$ contains $L_{\Bp}^{(\mu_j,0)}$ as an 
$\CS^{\Bp}$-submode by Lemma 3.10, there exists an $\CS^{\Bp}$-submodule
$W_j'$ of $W_j$ containing $W_{j-1}$ 
such that $W_j'/W_{j-1} \simeq L_{\Bp}^{(\mu_j,0)}$. 
By 3.7, we identify $Z_{\Bp}^{(\la.0)}$ as an $\CS^{\Bp}$-submodule of 
$W^{\la}$, and put $M_j = Z_{\Bp}^{(\la,0)} \cap W_j$ and 
$M_j' = Z_{\Bp}^{(\la,0)} \cap W_j'$.  We have a filtration of 
$Z_{\Bp}^{(\la,0)}$ by $\CS^{\Bp}$-modules
\begin{equation*}
0 = M_0 \subset M_1' \subset M_1 \subset \cdots \subset 
M_{k-1} \subset M_k' \subset M_k = Z_{\Bp}^{(\la,0)}\cap W^{\la} 
                                 = Z_{\Bp}^{(\la,0)}.
\end{equation*}
 We claim that 
\begin{equation*}
\tag{3.12.1}
M_{j-1} \ne M'_j \quad\text{ if }\quad \a_{\Bp}(\mu_j) = \a_{\Bp}(\la). 
\end{equation*}
\par\medskip
Note that (3.12.1) implies the proposition.  In fact,
$M_j'/M_{j-1} \simeq L_{\Bp}^{(\mu_j,0)}$ since it is isomorphic to a 
non-zero submodule of $L_{\Bp}^{(\mu_j,0)}$.    
It follows that $L_{\Bp}^{(\mu_j,0)}$ occurs in the composition 
series of $M_j/M_{j-1}$ for each $j$, and the proposition follows. 
\par
We show (3.12.1).  Assume that $\a_{\Bp}(\mu_j) = \a_{\Bp}(\la)$.
Then the image of 
$\vf_{T^{\mu_j}}^{(\mu_j,0)} \in Z_{\Bp}^{(\mu_j,0)}$
to $L_{\Bp}^{(\mu_j,0)}$ gives a non-zero element $\bar\vf_j$in 
$L_{\Bp}^{(\mu_j,0)}$ by Remark 2.9 (ii).  We choose 
$x_j \in W'_j \backslash W_{j-1}$ corresponding to $\bar\vf_j$ 
under the isomorphism $W'_j/W_{j-1} \simeq L_{\Bp}^{(\mu_j,0)}$.
Since $\bar\vf_j\vf_{\mu_j} = \bar\vf_j$, we have
$x_j\vf_{\mu_j} \in W'_j\backslash W_{j-1}$.
Since $x_j \in W'_j \subset W^{\la}$, one can write 
$x_j = \sum_{T \in \CT_0(\la)}r_T\vf_T$.  Now 
$\vf_{\mu_j}$ is a projection from $M$ onto $M^{\mu_j}$. 
Hence 
\begin{equation*}
x_j\vf_{\mu_j} = \sum_{T \in \CT_0^{\Bp}(\la,\mu_j)}r_T\vf_T. 
\end{equation*}
Here  
$\CT_0^{\Bp}(\la,\mu_j) \subset \CT_0^{\Bp}(\la) = J_{\Bp}(\la,0)$
since $\a_{\Bp}(\mu_j) = \a_{\Bp}(\la)$.  It follows that the right 
hand side of the above equation is contained in $Z_{\Bp}^{(\la,0)}$, and 
so $x_j\vf_{\mu_j} \in M_j'\backslash M_{j-1}$.
This proves (3.12.1), and the proposition follows. 
\end{proof}
Combining Proposition 3.3, Proposition 3.11 and Proposition 3.12, 
we have the following theorem (cf. [SawS, Theorem 13.6]).
%%%
\begin{thm}  %%%% Theorem 3.13
Assume that $R$ is a field.  For any $\la, \mu \in \vL^+$
such that $\a_{\Bp}(\la) = \a_{\Bp}(\mu)$, 
we have
\begin{equation*}
[\ol Z_{\Bp}^{\la} : \ol L_{\Bp}^{\mu}]_{\ol\CS^{\Bp}}
= [Z_{\Bp}^{(\la,0)} : L_{\Bp}^{(\mu,0)}]_{\CS^{\Bp}}
= [W^{\la} : L^{\mu}]_{\CS(\vL)}.
\end{equation*}
\end{thm}

%%%
%%%
\section{Structure theorem for $\ol\CS^{\Bp}$}
\para{4.1.}
In this section, we assume that 
$\vL = \wt\CP_{n,r}(\Bm)$.
For each $\mu \in \vL$, let $\wh N^{\Ba_{\Bp}(\mu)}$ be 
the $R$-submodule of $\CH$ 
spanned by $m_{\Fs\Ft}$ such that $\Fs, \Ft \in \Std(\la)$
with $\Ba_{\Bp}(\la) > \Ba_{\Bp}(\mu)$.
Since $\la \trianglerighteq \mu$ implies 
$\Ba_{\Bp}(\la) \ge \Ba_{\Bp}(\mu)$, $\wh N^{\Ba_{\Bp}(\mu)}$
is a two sided ideal of $\CH$.  Put 
$\wh M^{\mu} = M^{\mu} \cap \wh N^{\Ba_{\Bp}(\mu)}$.
Then $\wh M^{\mu}$ is an $\CH$-module with the basis 
 $\{ m_{S\Ft} \mid S \in \CT_0(\la,\mu), 
        \Ft \in \Std(\la), \Ba_{\Bp}(\la) > \Ba_{\Bp}(\mu)\}$.
We define an $\CH$-module $\ol M^{\mu}$ by 
$\ol M^{\mu} = M^{\mu}/\wh M^{\mu}$ and let 
$f : M^{\mu} \to \ol M^{\mu}$ be the natural surjection. 
Put $\ol m_{S\Ft} = f(m_{S\Ft})$ for a basis $m_{S\Ft} \in M^{\mu}$.
Then 
\begin{equation*}
\tag{4.1.1}
\{ \ol m_{S\Ft} \mid S \in \CT^{\Bp}_0(\la,\mu), 
  \Ft \in \Std(\la) \text{ for } \la \in \vL^+\}
\end{equation*}
gives a basis of $\ol M^{\mu}$.
\para{4.2.}
We write $\Bm = (m_1, \dots, m_r)$ in the form 
$\Bm = (\Bm^{[1]}, \dots, \Bm^{[g]})$ where  
$\Bm^{[k]} = (m_{p_k+1}, \dots, m_{p_k+r_k})$.
For each $n_k \in \ZZ_{\ge 0}$, put 
$\vL_{n_k} = \wt\CP_{n_k,r_k}(\Bm^{[k]})$ and 
$\vL^+_{n_k} = \CP_{n_k,r_k}(\Bm^{[k]})$.
($\vL_{n_k}$ or $\vL_{n_k^+}$ is regarded as the empty set if 
$n_k = 0$.)
Let $\mu = (\mu^{(1)}, \dots, \mu^{(r)}) \in \vL$ 
be an $r$-composition and  
write it as $\mu = (\mu^{[1]}, \dots, \mu^{[g]})$.
Then 
an $\mu$-tableau $\Ft = (\Ft^{(1)}, \dots, \Ft^{(r)})$
can be expressed as $\Ft = (\Ft^{[1]}, \dots, \Ft^{[g]})$ with 
$\Ft^{[k]} = (\Ft^{(p_k+1)}, \dots, \Ft^{(p_k+r_k)})$, where
$\Ft^{[k]}$ is a $\mu^{[k]}$-tableau.  Take $\la \in \vL^+$,
$\mu \in \vL$ such that $\a_{\Bp}(\la) = \a_{\Bp}(\mu)$. 
Then an $\la$-tableau $T = (T^{(1)}, \dots, T^{(r)})$ 
of type $\mu$ can be expressed as $T = (T^{[1]}, \dots, T^{[g]})$
with $T^{[k]} = (T^{(p_k+1)}, \dots, T^{(p_k+r_k)})$, where
$T^{[k]}$ is a $\la^{[k]}$-tableau of type $\mu^{[k]}$. 
\par
The following lemma is easily verified.
%%%%
\begin{lem}  %%%% Lemma 4.3
Let $\a = (n_1, \dots, n_g) \in \ZZ_{> 0}^g$ be such that 
$n_1 + \cdots + n_g = n$.  Then 
\begin{enumerate}
\item
The map $\mu \mapsto (\mu^{[1]}, \dots, \mu^{[g]})$ gives a bijection 
between $\{ \mu \in \vL \mid \a_{\Bp}(\mu) = \a\}$ and 
$\vL_{n_1}\times\cdots\times \vL_{n_g}$.
\item
The map $\la \mapsto (\la^{[1]}, \dots, \la^{[g]})$ gives a bijection 
between $\{ \la \in \vL^+ \mid \a_{\Bp}(\la) = \a\}$ and 
$\vL^+_{n_1}\times\cdots\times \vL^+_{n_g}$.
\item
For each $\la \in \vL^+, \mu \in \vL$ such that 
$\a_{\Bp}(\la) = \a_{\Bp}(\mu)$, the map 
$T \mapsto (T^{[1]}, \dots, T^{[g]})$ gives a bijection 
$\CT_0^{\Bp}(\la,\mu) \simeq \CT_0(\la^{[1]},\mu^{[1]}) \times 
   \cdots \times \CT_0(\la^{[g]},\mu^{[g]})$.
\end{enumerate}
\end{lem}
\para{4.4.}
Let $\a = (n_1, \dots, n_g) \in \wt \CP_{n,1}$.  For each
$\la \in \vL^+$ such that $\a_{\Bp}(\la) = \a$, we define 
a subset $\Std(\la)_0$ of $\Std(\la)$ as the set of 
$\Ft = (\Ft^{[1]}, \dots, \Ft^{[g]})$
such that the letters contained in the tableau $\Ft^{[k]}$ are exactly
$\{ n_1 + \cdots + n_{k-1} + 1, \dots, n_1 + \cdots + n_k\}$.
Then the set $\Std(\la)_0$ is in bijection with the set
$\Std(\la^{[1]}) \times \cdots \times \Std(\la^{[g]})$ under the 
map $\Ft \mapsto (\Ft^{[1]}, \dots, \Ft^{[g]})$. 
For each $\mu \in \vL$ such that $\a_{\Bp}(\mu) = \a$, we define 
an $R$-submodule $\ol M_0^{\mu}$ of $\ol M^{\mu}$  as the $R$-span 
of $\ol m_{S\Ft}$ such that 
$S \in \CT_0^{\Bp}(\la, \mu), \Ft \in \Std(\la)_0$ for various 
$\la \in \vL^+$.  
We write $\mu = (\mu^{[1]}, \dots, \mu^{[g]})$ as before.
Take $\Fs \in \Std(\la)$ such that $\mu(\Fs) = S$ for 
$S \in \CT_0^{\Bp}(\la,\mu)$ with $\a_{\Bp}(\la) = \a_{\Bp}(\mu) = \a$. 
Then $\Fs \in \Std(\la)_0$ and $\Fs^{[k]} \in \Std(\la^{[k]})$
has the property that $\mu^{[k]}(\Fs^{[k]}) = S^{[k]}$. 
This gives a bijection between the set of $\Fs \in \Std(\la)$
such that $\mu(\Fs) = S$ and the set of 
$(\Fs^{[1]}, \dots, \Fs^{[g]}) \in \Std(\la^{[1]}) 
    \times \cdots \times \Std(\la^{[g]})$ such that 
    $\mu^{[k]}(\Fs^{[k]}) = S^{[k]}$ for each $k$.
Combined with (1.4.1), (4.1.1), this implies that 
\par\medskip\noindent
(4.4.1) \  The map $\ol m_{S\Ft} \mapsto 
    m_{S^{[1]}\Ft^{[1]}}\otimes\cdots\otimes m_{S^{[g]}\Ft^{[g]}}$
gives an isomorphism of $R$-modules 
$\f_{\mu} : \ol M_0^{\mu} \isom M^{\mu^{[1]}}
                 \otimes\cdots\otimes M^{\mu^{[g]}}$. 
\par\medskip
Put $\CH_{\a} = \CH_{n_1,r_1}\otimes\cdots\otimes \CH_{n_g,r_g}$.
Since $M^{\mu^{[k]}}$ is an $\CH_{n_k,r_k}$-module, 
$M^{\mu^{[1]}}\otimes\cdots\otimes M^{\mu^{[g]}}$ has a structure of
an $\CH_{\a}$-module.  We denote by  
$T_0^{[k]}, \dots, T_{n_k-1}^{[k]}$ the generators of 
$\CH_{n_k,r_k}$ corresponding to $T_0, \dots, T_{n-1}$ 
in the case of $\CH_{n,r}$, and more generally we denote by 
$T_w^{[k]}$ for $w \in \FS_{n_k}$ the element corresponding to 
$T_w \in \CH_n$.
Then $T_i^{[k]}$ acts on  
$M^{\mu^{[1]}}\otimes\cdots\otimes M^{\mu^{[g]}}$ for 
$i = 0, \dots, n_k-1$, through the action of 
$1^{\otimes (k-1)}\otimes T_i^{[k]}
        \otimes 1^{\otimes (g-k)} \in \CH_{\a}$ on it.
\par
Recall that $L_i = T_{i-1}T_{i-2}\cdots T_1T_0T_1\cdots 
                     T_{i-2}T_{i-1} \in \CH$
for $i = 0, \dots, n-1$.
The following lemma is crucial for later discussions.
%%%
\begin{lem} %%% Lemma 4.5
Let $\mu \in \vL$ be such that $\a_{\Bp}(\mu) = \a$.
For $\mu \in \vL$, put $\Ba_{\Bp}(\mu) = (a_1, \dots, a_g)$.
Then the action of $L_{a_k+1}$ on $\ol M^{\mu}$ stabilizes the 
submodule $\ol M_0^{\mu}$, and it gives rise
to the action of $T_0^{[k]}$ on 
$M^{\mu^{[1]}}\otimes\cdots\otimes M^{\mu^{[g]}}$ under the
isomorphism $\f_{\mu}$ in (4.4.1).
\end{lem}
\begin{proof}
Take $\la \in \vL^+$ such that $\Ba_{\Bp}(\la) = \Ba_{\Bp}(\mu)$
and consider $\CT_0^{\Bp}(\la,\mu)$.
Let $\Fs \in \Std(\la)$ such that $\mu(\Fs) = S$ for 
$S = (S^{[1]}, \dots, S^{[g]}) \in \CT_0^{\Bp}(\la,\mu)$. Then 
$\Fs = (\Fs^{[1]}, \dots, \Fs^{[g]}) \in \Std(\la)_0$ 
and $\mu^{[k]}(\Fs^{[k]}) = S^{[k]}$.
Take $\Fs$ as above, and take 
$\Ft = (\Ft^{[1]}, \dots, \Ft^{[g]}) \in \Std(\la)_0$.
We consider a basis 
$m_{\Fs\Ft} \in \CH$ and $m_{\Fs^{[k]}\Ft^{[k]}} \in \CH_{n_k,r_k}$.
We show that
\par\medskip\noindent
(4.5.1) \ $m_{\Fs\Ft}L_{a_k+1}$ is written as a linear 
combination of the basis elements $m_{\Fu\Fv}$ of $\CH$, where
$\Fu = (\Fu^{[1]}, \dots, \Fu^{[g]})$ 
is obtained from $\Fs$ by replacing $\Fs^{[k]}$ by some $\Fu^{[k]}$,
and $\Fv$ is obtained from $\Ft$ similarly.  Here $\Fu^{[k]}$ and
$\Fv^{[k]}$ has the same shape.
The coefficient of $m_{\Fu\Fv}$ in the expansion of 
$m_{\Fs\Ft}L_{a_k+1}$ coincides with the coefficient of 
$m_{\Fu^{[k]}\Fv^{[k]}}$ in the expansion of 
$m_{\Fs^{[k]}\Ft^{[k]}}T_0^{[k]}$ under the bijection 
$\Fu \lra \Fu^{[k]}, \Fv\lra \Fv^{[k]}$.  
\par\medskip
(4.5.1) implies the lemma since $\Fu, \Fv$ are standard tableau
of shape $\nu$ with $\Ba_{\Bp}(\nu) = \Ba_{\Bp}(\mu)$ and 
$\Fv \in \Std(\nu)_0$.  We shall show (4.5.1).
First we compute $m_{\Fs^{[k]}\Ft^{[k]}}T_0^{[k]}$ following 
the argument in the proof of [DJM, Proposition 3.20].
Recall that
\begin{equation*}
\la^{[k]} = (\la^{(p_k+1)}, \dots, \la^{(p_k + r_k)}), \qquad
\Ft^{[k]} = (\Ft^{(p_k+1)}, \dots, \Ft^{(p_k + r_k)}), 
\end{equation*}
and put 
\begin{align*}
\b &= (|\la^{(p_k+1)}|, \dots, |\la^{(p_k+r_k)}|) 
                   = (\b_1, \dots, \b_{r_k}), \\
\Bb &= (b_1, \dots, b_{r_k}) \text{ with } b_j = \sum_{i=1}^{j-1}\b_i.
\end{align*}
The letters contained in $\Ft^{[k]}$ consist of 
$\{ a_k+1, \dots, a_k+ n_k\}$.  By the shift by $-a_k$, we regard
$\Ft^{[k]}$ as the tableau consisting of letters $\{ 1, \dots, n_k\}$.
Assume that the letter $1$ is contained in $\Ft^{(p_k + f)}$.
One can write $d(\Ft^{[k]}) = yc$, where $y \in \FS_{\b}$ and
$c$ is a distinguished coset representative in 
$\FS_{\b}\backslash  \FS_{n_k}$.
Then $\Ft^{\la^{[k]}}y$ is a standard tableau, and 
$c$ is a permutation which maps the letters 
$\{ b_i +1, \dots, b_i + \b_i\}$ to the letters contained in 
$\Ft^{(p_k + i)}$ for $i = 1, \dots, r_k$.
Thus $y$ fixes the letter $b_f + 1$, and $c$ can be expressed as 
$c = (b_{f+1}, b_{f})(b_{f}, b_{f-1})\cdots (2,1)c'$, where 
$l(c) = b_f + l(c')$ and $c'$ fixes the letter 1.
It follows that 
$T^{[k]}_c = T^{[k]}_{b_f}T^{[k]}_{b_f-1}\cdots 
     T_2^{[k]}T_1^{[k]}T_{c'}^{[k]}$ 
and $T^{[k]}_{c'}T_0^{[k]} = T_0^{[k]}T^{[k]}_{c'}$.
Recall that $m_{\Fs^{[k]}\Ft^{[k]}} = 
T^{[k]*}_{d(\Fs^{[k]})}m_{\la^{[k]}}T^{[k]}_{d(\Ft^{[k]})}$
with $m_{\la^{[k]}} = u_{\Bb}^+x_{\la^{[k]}} = 
x_{\la^{[k]}}u_{\Bb}^+$.
Here $u_{\Bb}^+ = u_{\Bb,1}\cdots u_{\Bb,r_k}$, 
with 
\begin{equation*}
u_{\Bb, j} = \prod_{i=1}^{b_j}(L^{[k]}_i - Q^{[k]}_j),
\end{equation*}
where $L_i^{[k]}$ is the element in $\CH_{n_k,r_k}$ 
corresponding to $L_i \in \CH$, and $Q^{[k]}_j = Q_{p_k + j}$.
Then as in the computation in [DJM, Prop.3.20, Lemma 3.4], 
by noticing that $T_y^{[k]}$ commutes with $L_{b_f+1}^{[k]}$,  we have
\begin{align*}
u_{\Bb}^+T_y^{[k]}T^{[k]}_cT^{[k]}_0 
&= u_{\Bb}^+T_y^{[k]}T^{[k]}_{b_f}T^{[k]}_{b_f-1}
     \cdots T_1^{[k]}T_0^{[k]}T_{c'}^{[k]} \\
&= u_{\Bb}^+L^{[k]}_{b_f+1}T_y^{[k]}(T^{[k]}_{b_f})\iv\cdots 
        (T_1^{[k]})\iv T^{[k]}_{c'}  \\
&= (Q^{[k]}_fu_{\Bb}^+ + u_{\Bb'}^+)T_y^{[k]}(T^{[k]}_{b_f})\iv\cdots 
        (T_1^{[k]})\iv T^{[k]}_{c'}
\end{align*}
where $\Bb' = (b_1, \dots, b_{f-1}, b_{f}+1, b_{f+1}, \dots, b_{r_k})$.
It follows that 
\begin{equation*}
\tag{4.5.2}
m_{\Fs^{[k]}\Ft^{[k]}}T_0^{[k]} = T^{[k]*}_{d(\Fs^{[k]})}
                 x_{\la^{[k]}}(Q^{[k]}_fu_{\Bb}^+ + 
                u_{\Bb'}^+ )T^{[k]}_yh, 
\end{equation*}
with
\begin{equation*}
h = (T^{[k]}_{b_f})\iv\cdots (T_1^{[k]})\iv T_{c'}^{[k]}.
\end{equation*}
\par
Next we  shall compute 
$m_{\Fs\Ft}L_{a_k+1}$ for $\Fs, \Ft \in \Std(\la)_0$.
Recall that $m_{\Fs\Ft} = T^*_{d(\Fs)}m_{\la}T_{d(\Ft)}$
with $m_{\la} = x_{\la}u_{\Ba}^+$.
Since $\Ft \in \Std(\la)_0$, we have 
$d(\Ft) = d(\Ft^{[1]})\cdots d(\Ft^{[g]})$.  (Note that 
the letters contained in $\Ft^{[k]}$ consist of
$\{ a_k +1, \dots, a_k+ n_k\}$, and we compute 
$d(\Ft^{[k]})$ with respect to these letters.)
We note that
for $\Ft = (\Ft^{[1]}, \dots, \Ft^{[g]})$, 
the letters contained in 
$\Ft^{[1]}, \dots, \Ft^{[k-1]}$ consist of $\{ 1,2, \dots, a_k\}$, 
and the letters contained in $\Ft^{[k]}$ consist of 
$\{ a_k+1, \dots, a_k + n_k = a_{k+1}\}$, the letters 
contained in $\Ft^{[k+1]}, \dots, \Ft^{[g]}$ consist of
$\{ a_{k+1}+1, \dots, n\}$. 
It follows that 
\begin{equation*}
\tag{4.5.3}
\begin{split}
m_{\Fs\Ft}T_{a_k}T_{a_k-1}\cdots T_1T_0 
    = T^*_{d(\Fs)}&x_{\la}u_{\Ba}^+T_{d(\Ft^{[1]})}\cdots
    T_{d(\Ft^{[k]})} \times \\
    &\times T_{a_k}T_{a_k-1}\cdots 
              T_1T_0 T_{d(\Ft^{[k+1]})}\cdots T_{d(\Ft^{[g]})}. 
\end{split}
\end{equation*}
From the previous computation, we have  $d(\Ft^{[k]}) =  yc$
with $y \in \FS_{\b}$ and $c \in \FS_{\b}\backslash \FS_{n_k}$.
(Here we regard $\FS_{n_k}$ as the permutation group with respect 
to the letters $\{ a_k+1, \dots, a_k+n_k\}$.  In particular, 
$y$ fixes the letter $a_k + b_f+1$). 
Hence
\begin{equation*}
T_{d(\Ft^{[k]})} = T_yT_c 
          = T_yT_{a_k + b_f}T_{a_k + b_f -1}\cdots T_{a_k+1}T_{c'}.
\end{equation*}
Let $X$ be the left hand side of (4.5.3). Since $T_{c'}$ commutes with
$T_{a_k}, \dots, T_1, T_0$, we have
\begin{equation*}
\tag{4.5.4}
\begin{split}
X = T^*_{d(\Fs)}&x_{\la}u_{\Ba}^+T_{d(\Ft^{[1]})}\cdots
    T_{d(\Ft^{[k-1]})}T_y \times \\
    &\times T_{a_k + b_f}\cdots T_{a_k+1}T_{a_k}T_{a_k-1}\cdots 
              T_1T_0T_{c'} T_{d(\Ft^{[k+1]})}\cdots T_{d(\Ft^{[g]})}. 
\end{split}
\end{equation*}
Recall that $\Ba = \Ba(\la) = (a_1', \dots, a_r')$ is defined 
by $a_j' = \sum_{i=1}^{j-1}|\la^{(i)}|$, and $u_{\Ba}^+$
is given by  
$u_{\Ba}^+ = u_{\Ba,1}u_{\Ba,2}\cdots u_{\Ba,r}$, 
where 
$u_{\Ba,j} = \prod_{i=1}^{a_j'}(L_i - Q_j)$. 
Hence $\Ba_{\Bp} = (a_1, \dots, a_g)$ is given by 
$a_i = a'_{p_i +1}$ for $i = 1, \dots, g$.  
Put
\begin{equation*}
\tag{4.5.5}
u_{\Ba_{\Bp}, i} = u_{\Ba, p_i+1}\cdots u_{\Ba, p_i+ r_i} 
\end{equation*}
for $i = 1, \dots, g$.
Then we have $u_{\Ba}^+ = u_{\Ba_{\Bp},1}\cdots u_{\Ba_{\Bp}, g}$
and $u_{\Ba_{\Bp}, k}, \dots, u_{\Ba_{\Bp}, g}$ commutes with 
$T_{d(\Ft^{[1]})}, \dots, T_{d(\Ft^{[k-1]})}$, and 
$u_{\Ba_{\Bp}, k+1}, \dots, u_{\Ba_{\Bp}, g}$ commutes with $T_y$.
It follows that 
\begin{align*}
u_{\Ba}^+&T_{d(\Ft^{[1]})}\cdots T_{d(\Ft^{[k-1]})} 
   T_yT_{a_k + b_f}\cdots T_1T_0 \\
&= u_{\Ba_{\Bp},1}\cdots u_{\Ba_{\Bp},k-1}
        T_{d(\Ft^{[1]})}\cdots T_{d(\Ft^{[k-1]})}
  u_{\Ba_{\Bp},k}T_yu_{\Ba_{\Bp},k+1}
           \cdots u_{\Ba_{\Bp},g}T_{a_k + b_f}\cdots T_1T_0.
\end{align*}
Since $u_{\Ba_{\Bp}, k+1}, \dots, u_{\Ba_{\Bp},g}$ commutes with 
$T_{a_k+b_f}, \dots, T_1, T_0$, we have
\begin{equation*}
u_{\Ba_{\Bp},k}T_yu_{\Ba_{\Bp},k+1}
    \cdots u_{\Ba_{\Bp},g}T_{a_k + b_f}\cdots T_1T_0
   = u_{\Ba_{\Bp},k}T_yT_{a_k + b_f}\cdots T_1T_0
       u_{\Ba_{\Bp},k+1}\cdots u_{\Ba_{\Bp},g}.
\end{equation*}
Since $T_{a_k + b_f}\cdots T_1T_0 = L_{a_k + b_f+1}h'$ with 
$h' = T_{a_k+ b_f}\iv\cdots T_1\iv$, and $T_y$ commutes with
$L_{a_k + b_f+1}$,  we have by [DJM, Lemma 3.4],
\begin{equation*}
\begin{split}
u_{\Ba_{\Bp},k}T_y&T_{a_k + b_f}\cdots T_1T_0
       u_{\Ba_{\Bp},k+1}\cdots u_{\Ba_{\Bp},g} \\
&= (Q_{p_k + f}u_{\Ba_{\Bp},k} + u'_{\Ba_{\Bp},k})T_yh'
        u_{\Ba_{\Bp},k+1}\cdots u_{\Ba_{\Bp},g},
\end{split}
\end{equation*}
where $u'_{\Ba_{\Bp}, k}$ is defined as in (4.5.5) by 
replacing $\Ba$ by 
\begin{equation*}
\Ba' = (a_1', \dots, a'_{p_k + f-1}, a'_{p_k+f}+1, 
       a'_{p_k + f + 1}, \dots, a'_r).
\end{equation*}
Summing up the above computation, we have
\begin{align*}
X = T^*_{d(\Fs)}&x_{\la}u_{\Ba_{\Bp},1}\cdots u_{\Ba_{\Bp},k-1}
     T_{d(\Ft^{[1]})}\cdots T_{d(\Ft^{[k-1]})} \times \\
   & \times (Q_{p_k + f}u_{\Ba_{\Bp},k} + u'_{\Ba_{\Bp},k})T_yh'
        T_{c'}u_{\Ba_{\Bp},k+1}\cdots u_{\Ba_{\Bp},g}
        T_{d(\Ft^{[k+1]})}\cdots T_{d(\Ft^{[g]})}
\end{align*}
 It follows that 
\begin{align*}
\tag{4.5.6}
m_{\Fs\Ft}L_{a_k + 1} &= XT_1\cdots T_{a_k} \\
  &= T^*_{d(\Fs)}x_{\la}u_{\Ba_{\Bp},1}\cdots u_{\Ba_{\Bp},k-1}
     T_{d(\Ft^{[1]})}\cdots T_{d(\Ft^{[k-1]})} \times \\
   & \hspace{1cm}
      \times (Q_{p_k + f}u_{\Ba_{\Bp},k} + u'_{\Ba_{\Bp},k})T_yh''
        T_{c'}u_{\Ba_{\Bp},k+1}\cdots u_{\Ba_{\Bp},g}
        T_{d(\Ft^{[k+1]})}\cdots T_{d(\Ft^{[g]})}      
\end{align*}
where $h'' = T_{a_k + b_f}\iv\cdots T_{a_k+1}\iv$.
\par
We now compare (4.5.2) and (4.5.6).
The right hand side of (4.5.2) is written as 
$X_1 + X_2$, where 
$X_1 =  Q_f^{[k]}T_{d(\Fs^{[k]})}^{[k]*}
       x_{\la^{[k]}}u_{\Bb}^+T_y^{[k]}h$ 
and $X_2 = T_{d(\Fs^{[k]})}^{[k]*}
       x_{\la^{[k]}}u_{\Bb'}^+T_y^{[k]}h$.
Since $x_{\la^{[k]}}u_{\Bb}^+ = m_{\la^{[k]}}$, 
$X_1$ can be written, by Lemma 3.15 in [DJM],   
as a linear combination of the elements 
$m_{\Fs^{[k]}\Ft_i^{[k]}}$, where $\Ft_i^{[k]}$ are 
row-standard tableaux.  Then
they are converted 
to a linear combination of the basis elements 
$m_{\Fu^{[k]}\Fv^{[k]}}$ in $\CH_{n_k,r_k}$ 
by the procedure given in 
Proposition 3.18 in [loc. cit.], where 
$\Fu^{[k]}, \Fv^{[k]}$ are standard tableaux of 
shape $\mu^{[k]}$ for some $r_k$-partitions $\mu^{[k]}$.
On the other hand, for $X_2$, 
first we convert  $T_{d(\Fs^{[k]})}^{[k]*}
       x_{\la^{[k]}}u_{\Bb'}^+$ 
to a linear combination of the elements
$m_{\Fu_1^{[k]}\Fv_1^{[k]}}$ where $\Fu_1^{[k]}, \Fv_1^{[k]}$ are 
row-standard tableau of shape $\nu^{[k]}$ 
($\nu^{[k]}$ is determined from $u_{\Bb'}$), and 
then we follow the argument in the case $X_1$. 
Note that in these computations, the parts 
$u_{\Bb}^+$ and $u_{\Bb'}^+$ remain unchanged.
\par
Next we consider (4.5.6).  Since 
$T_{d(\Fs)} = T_{d(\Fs^{[1]})}\cdots T_{d(\Fs^{[g]})}$ and 
$x_{\la} = x_{\la^{[1]}}\cdots x_{\la^{[g]}}$, 
one can write the formula (4.5.6) in the form 
\begin{equation*}
m_{\Fs\Ft}L_{a_k+1} = Z\cdot T^*_{d(\Fs^{[k]})}x_{\la^{[k]}}
       (Q_{p_k + f}u_{\Ba_{\Bp},k} + u'_{\Ba_{\Bp},k})T_yh''
        T_{c'}\cdot Z'
\end{equation*}
where
\begin{align*}
Z &= T^*_{d(\Fs^{[1]})}x_{\la^{[1]}}u_{\Ba_{\Bp},1}T_{d(\Ft^{[1]})}
          \cdots T^*_{d(\Fs^{[k-1]})}x_{\la^{[k-1]}}u_{\Ba_{\Bp}, k-1}
                       T_{d(\Ft^{[k-1]})}  \\
Z' &= T^*_{d(\Fs^{[k+1]})}x_{\la^{[k+1]}}u_{\Ba_{\Bp},k+1}T_{d(\Ft^{[k+1]})}
          \cdots T^*_{d(\Fs^{[g]})}x_{\la^{[g]}}u_{\Ba_{\Bp}, g}
                       T_{d(\Ft^{[g]})}.
\end{align*}
Put 
\begin{align*}
Y_1 &= Q_{p_k + f}T^*_{d(\Fs^{[k]})}x_{\la^{[k]}}
                 u_{\Ba_{\Bp},k}T_yh''T_{c'}, \\ 
Y_2 &= T^*_{d(\Fs^{[k]})}x_{\la^{[k]}}
                   u'_{\Ba_{\Bp},k}T_yh''T_{c'}
\end{align*}
so that $m_{\Fs\Ft}L_{a_k+1} = Z(Y_1 + Y_2)Z'$.
Let $\CH'_{n_k}$ be the subalgebra of $\CH_{n}$ generated by
$T_{a_k + 1}, \cdots, T_{a_k + n_k-1}$.  Then 
$T_y, T_{c'}, h''$ belong to $\CH'_{n_k}$, and under the 
identification $\CH'_{n_k} \simeq \CH_{n_k}$, 
$T_y, T_{c'}$ coincide with $T_y^{[k]}, T_{c'}^{[k]}$, 
 and $h''T_{c'}$ coincides with $h$. 
We also note that $Q_{p_k + f} = Q_f^{[k]}$.
Now by applying Lemma 3.15 and Proposition 3.18 in [loc.cit], 
$Y_1$ can be expressed as a linear combination of 
the terms 
$T^*_{d(\Fu^{[k]})}x_{\mu^{[k]}}u_{\Ba_{\Bp},k}T_{d(\Fv^{[k]})}$,
where $\Fu^{[k]}, \Fv^{[k]}$ are standard tableaux of shape 
$\mu^{[k]}$ for some $r_k$-partitions $\mu^{[k]}$. 
Since this computation proceeds without referring $u_{\Ba_{\Bp},k}$, 
the coefficients of these elements in the expansion of $Y_1$ 
are exactly the same as the coefficients of 
$m_{\Fu^{[k]}\Fv^{[k]}}$ in the expansion of $X_1$.  
For $Y_2$, first we convert 
$T^*_{d(\Fs^{[k]})}x_{\la^{[k]}}u'_{\Ba_{\Bp},k}$
to a linear combination of the terms 
$T^*_{d(\Fu_1^{[k]})}x_{\nu^{[k]}}u'_{\Ba_{\Bp},k}T_{d(\Fv_1^{[k]})}$
by using Proposition 3.20.  By comparing $\Bb'$ and $\Ba'$, we see
that the coefficients in this expansion are exactly the same as 
the coefficients of $m_{\Fu^{[k]}_1\Fv^{[k]}_1}$ in the expansion of
$T^{[k]*}_{d(s^{[k]})}x_{\la^{[k]}}u_{\Bb'}^+$.  
Thus again by applying Lemma 3.15 and Proposition 3.18, 
we conclude that $Y_2$ can be written as a linear combination of the 
terms $T^*_{d(\Fu^{[k]})}x_{\nu^{[k]}}u'_{\Ba_{\Bp},k}T_{d(\Fv^{[k]})}$, 
where $\Fu^{[k]}, \Fv^{[k]}$ are standard tableaux of shape 
$\nu^{[k]}$, and that their coefficients in the expansion of $Y_2$ 
is the same as the coefficients of $m_{\Fu^{[k]}\Fv^{[k]}}$ in the
expansion of $X_2$.
\par
Now one sees easily that 
$Z\cdot T^*_{d(\Fu^{[k]})}x_{\mu^{[k]}}u_{\Ba_{\Bp},k}
       T_{d(\Fv^{[k]})}\cdot Z' = m_{\Fu\Fv}$, where 
$\mu$ is an $r$-partition obtained from $\la$ by replacing 
$\la^{[k]}$ by $\mu^{[k]}$, and $\Fu, \Fv$ are standard tableau
of shape $\mu$ obtained from $\Fs, \Ft$ by replacing 
$\Fs^{[k]}, \Ft^{[k]}$ by $\Fu^{[k]}, \Fv^{[k]}$.
A similar result holds also for 
$Z\cdot T^*_{d(\Fu^{[k]})}x_{\nu^{[k]}}u'_{\Ba_{\Bp},k}
       T_{d(\Fv^{[k]})}\cdot Z'$.
Summing up the above arguments, we see that (4.5.1) holds. 
Hence the lemma is proved.
\end{proof}
The following lemma is easily verified by using a similar
(but simpler) argument as in the proof of the previous lemma. 
\begin{lem}  %%% Lemma 4.6.
Let the notations be as in Lemma 4.5.  Then, for 
$i = 1, \dots, n_k-1$,  the action of
$T_{a_k+i}$ on $\ol M^{\mu}$ stabilizes $\ol M_0^{\mu}$, and 
it gives rise to the action of
$T_i^{[k]}$ on $M^{\mu^{[1]}}\otimes\cdots\otimes M^{\mu^{[g]}}$
under the identification $\f_{\mu}$ in (4.4.1).
\end{lem}
\para{4.7.}
We fix $\a = (n_1, \dots, n_g) \in \wt\CP_{n,1}$, 
and let $\CH_{\a}$ be as in 4.4. 
Assume that $\a_{\Bp}(\mu) = \a$ for $\mu \in \vL$.
Then $\CH_{\a}$ acts naturally on  
$M^{\mu^{[1]}}\otimes\cdots\otimes M^{\mu^{[g]}}$.  Let
$\Ba_{\Bp}(\mu) = (a_1, \dots, a_g)$ be as before, and let 
$\wt\CH_{\a}$ be the subalgebra of $\CH$ generated by
$T_{a_k +1}, \dots, T_{a_k + r_k-1}, L_{a_k+1}$  
for $k = 1, \dots, g$.
As a corollary to Lemma 4.5 and Lemma 4.6, 
we have the following.
\begin{cor}  %%% Cor. 4.8
For each $\mu \in \vL$ such that $\a_{\Bp}(\mu) = \a$, 
$\ol M_0^{\mu}$ is $\wt H_{\a}$-stable.  
The action of $\wt H_{\a}$ on $\ol M_0^{\mu}$ 
coincides with the action of
$\CH_{\a}$ on $M^{\mu^{[1]}}\otimes\cdots\otimes M^{\mu^{[g]}}$.
\end{cor}
\para{4.9.}
Recall that $\CS = 
   \bigoplus_{\mu, \nu \in \vL}\Hom_{\CH}(M^{\nu}, M^{\mu})$.  
It follows from the description of the basis of $\CS^{\Bp}$, 
we see that 
\begin{equation*}
\tag{4.9.1}
\CS^{\Bp} = \bigoplus_{\mu,\nu \in \vL}H_{\mu\nu}
\end{equation*}
where $H_{\mu\nu} = \CS^{\Bp} \cap \Hom_{\CH}(M^{\nu}, M^{\mu})$ 
is an $R$-submodule of $\Hom_{\CH}(M^{\nu}, M^{\mu})$
spanned by $\vf_{ST}$ with $S \in \CT_0(\la,\mu)$, 
$T \in \CT_0(\la,\nu)$ such that $\Ba_{\Bp}(\la) > \Ba_{\Bp}(\mu)$
if $\a_{\Bp}(\mu) \ne \a_{\Bp}(\nu)$.
Then we have 
\begin{equation*}
\tag{4.9.2}
\ol\CS^{\Bp} = \bigoplus_{\substack{\mu, \nu \in \vL \\
                   \a_{\Bp}(\mu) = \a_{\Bp}(\nu)}}\ol H_{\mu\nu}, 
\end{equation*}
where $\ol H_{\mu\nu} = \pi(H_{\mu\nu})$ is the $R$-span of the
elements $\bar\vf_{ST}$ such that $S \in \CT_0^{\Bp}(\la, \mu)$, 
$T \in \CT_0^{\Bp}(\la, \nu)$ for various $\la \in \vL^+$.
\par
Assume that $\a_{\Bp}(\mu) = \a_{\Bp}(\nu)$.
We claim that any $\vf \in H_{\mu\nu}$ maps $\wh M^{\nu}$ into 
$\wh M^{\mu}$. 
In fact take $\vf \in H_{\mu\nu}$. Then by the property of $\vf_{ST}$, 
there exists $h_{\vf} \in \CH$ such that 
$\vf(m_{\nu}h) = h_{\vf}m_{\nu}h$ for any $h \in \CH$.
Recall that $\wh M^{\nu}$ is a linear combination of $m_{S\Ft}$
with $S \in \CT_0(\la, \nu), \Ft \in \Std(\la)$ such that 
$\Ba_{\Bp}(\la) > \Ba_{\Bp}(\nu)$. Suppose that $m_{S\Ft}$ is 
written as $m_{S\Ft} = m_{\nu}h$ for some $h \in \CH$.  Then 
by the property of cellular basis, 
$\vf(m_{St}) = h_{\vf}m_{S\Ft}$ is a linear combination of 
$m_{\Fs'\Ft'}$, where $\Fs', \Ft' \in \Std(\la')$ with 
$\la' \trreq \la$.  Then we have 
$\Ba_{\Bp}(\la') \ge \Ba_{\Bp}(\la) > \Ba_{\Bp}(\nu)$.
Since $\Ba_{\Bp}(\nu) = \Ba_{\Bp}(\mu)$, we have 
$\Ba_{\Bp}(\la') > \Ba_{\Bp}(\mu)$, and so 
$\vf(m_{S\Ft}) \in \wh M^{\mu}$.  Thus the claim holds.  
\par
By the claim, $\vf$ induces a linear map 
$\bar\vf \in \Hom_{\CH}(\ol M^{\nu}, \ol M^{\mu})$ under the condition 
that $\a_{\Bp}(\mu) = \a_{\Bp}(\nu)$.  We note that $\bar\vf = 0$
if $\vf \in \wh \CS^{\Bp}$.  In fact, since 
$\Ba_{\Bp}(\mu) = \Ba_{\Bp}(\nu)$, we may consider the case where 
$\vf = \vf_{ST}$
for $S \in \CT_0(\la,\mu), T \in \CT_0(\la,\nu)$ with 
$\Ba_{\Bp}(\la) \ne \Ba_{\Bp}(\mu)$.  Since $\la \trreq \mu$, 
we have $\Ba_{\Bp}(\la) > \Ba_{\Bp}(\mu)$.  
It follows that 
$\vf_{ST}(m_{\nu}) = m_{ST} \in \wh M^{\mu}$, and the image of $\vf$
is contained in $\wh M^{\mu}$. Hence $\bar\vf = 0$ as asserted.  
\par
The above discussion allows us to define a linear map
$\th : H_{\mu\nu} \to \Hom_{\CH}(\ol M^{\nu}, \ol M^{\mu})$
by $\vf \mapsto \bar\vf$, which factors through the map 
$\bar \th : \ol H_{\mu\nu} \to \Hom_{\CH}(\ol M^{\nu}, \ol M^{\mu})$.
We show the following lemma. 
%%%
\begin{lem} % Lemma 4.10 
\begin{enumerate}
\item
For $\mu \in \vL$, let $\f_{\mu} : \ol M_0^{\mu} \to 
M^{\mu^{[1]}}\otimes \cdots\otimes M^{\mu^{[g]}}$ be the 
isomorphism given in (4.4.1).  Then we have
\begin{equation*}
\f_{\mu}\iv(m_{\mu^{[1]}}\otimes\cdots\otimes m_{\mu^{[g]}}) 
    = \ol m_{\mu}.
\end{equation*}
\item
Assume that $\a_{\Bp}(\mu) = \a_{\Bp}(\nu) = \a$.  
Then for any $\vf \in \ol H_{\mu\nu}$, 
$\bar\vf = \bar\th(\vf)$ maps $\ol M_0^{\nu}$ to $\ol M_0^{\mu}$.
In particular, 
$\bar\vf \in \Hom_{\wt\CH_{\a}}(\ol M_0^{\nu}, \ol M_0^{\mu})$.
\end{enumerate}
\end{lem}
\begin{proof}
First we show (i).
Put $\Ba = \Ba(\mu)$ and $\Ba_{\Bp} = \Ba_{\Bp}(\mu)$. 
Then $m_{\mu} = x_{\mu}u_{\Ba}^+$, and 
$x_{\mu} = x_{\mu^{[1]}}\cdots x_{\mu^{[g]}}$, 
$u_{\Ba}^+ = u_{\Ba_{\Bp},1}\cdots u_{\Ba_{\Bp},g}$, where
$u_{\Ba_{\Bp},i}$ is defined as in (4.5.5).
One can write $m_{\mu} = x_1x_2\cdots x_g$ with 
$x_k = x_{\mu^{[k]}}u_{\Ba_{\Bp},k}$.
On the other hand, $m_{\mu^{[k]}} = x_{\mu^{[k]}}u_{\Bb}^+$, 
where $\Bb = \Ba(\mu^{[k]})$ is defined with respect to 
$\mu^{[k]} \in \wt\CP_{n_k,r_k}(m^{[k]})$.  Then by 
Proposition 3.18 in [DJM], $m_{\mu^{[k]}}$ is written as a linear 
combination of the basis elements $m_{\Fu^{[k]}\Fv^{[k]}}$ of 
$\CH_{n_k, r_k}$, where $\Fu^{[k]}, \Fv^{[k]}$ are standard tableau of 
shape $\la^{[k]}$.  By the same procedure, $x_k$ is written as a 
linear combination of 
$x_{\Fu^{[k]}\Fv^{[k]}} = 
 T^*_{d(\Fu^{[k]})}x_{\la^{[k]}}u_{\Ba_{\Bp},k}T_{d(\Fv^{[k]})}$, 
and the corresponding coefficient coincides each other. 
Note that in the latter case $\Fu^{[k]}, d(\Fu^{[k]})$, etc. 
are referred   
with respect to the letters $\{ a_k +1, \dots, a_k + n_k\}$ as in 
the proof of Lemma 4.5. 
We see that $x_{\Fu^{[1]}\Fv^{[1]}}\cdots x_{\Fu^{[g]}\Fv^{[g]}}$ 
gives rise to a basis element $m_{\Fu\Fv}$ of $\CH$, where 
$\Fu = (\Fu^{[1]}, \dots, \Fv^{[g]})$ and 
$\Fv = (\Fv^{[1]}, \dots, \Fv^{[g]})$ are
in $\Std(\la)_0$ with $\la = (\la^{[1]}, \dots, \la^{[g]})$.
The assertion (i) follows from this. 
\par
Next we show (ii). Now we have $\ol m_{\nu} \in \ol M_0^{\nu}$.
Since $M^{\nu^{[1]}}\otimes \cdots \otimes M^{\nu^{[g]}}$ is 
generated by $m_{\nu^{[1]}}\otimes\cdots\otimes m_{\nu^{[g]}}$
as an $\CH_{\a}$-module, $\ol M_0^{\nu}$ is generated by $\ol m_{\nu}$
as an $\wt H_{\a}$-module.  We take $\bar\vf_{ST} \in \ol H_{\mu\nu}$.
Then any element in $\ol M_0^{\nu}$ is written as 
$\ol m_{\nu}h$ with $h \in \wt \CH_{\a}$, and 
$\vf_{ST}(m_{\nu}h) = \vf_{ST}(m_{\nu})h = m_{ST}h$.
Since $\ol m_{ST} \in \ol M_0^{\mu}$, we see that 
$\bar\vf_{ST}(\ol M_0^{\nu}) \subseteq \ol M_0^{\mu}$. 
This proves (ii), and the lemma follows.
\end{proof}
\para{4.11.}
We keep the previous setting.
By Lemma 4.10, one can define an $R$-linear map
$\vT : \ol H_{\mu\nu} \to 
    \Hom_{\wt\CH_{\a}}(\ol M_0^{\nu}, \ol M_0^{\mu})$
induced from $\bar\th$.
On the other hand, in view of the isomorphisms 
$\f_{\mu}, \f_{\nu}$ together with Corollary 4.8, we have 
a natural isomorphism of $R$-modules
\begin{equation*}
\tag{4.11.1}
\begin{split}
\Hom&_{\wt\CH_{\a}}(\ol M_0^{\nu}, \ol M_0^{\mu})  \\
   &\simeq \Hom_{\CH_{n_1,r_1}}(M^{\nu^{[1]}}, M^{\mu^{[1]}})\otimes
       \cdots\otimes \Hom_{\CH_{n_g,r_g}}(M^{\nu^{[g]}}, M^{\mu^{[g]}}).
\end{split}
\end{equation*}
We have the following lemma.
%%%
\begin{lem}  %%% Lemma 4.12
The map $\vT$ gives an isomorphism 
\begin{equation*}
\ol H_{\mu\nu} \simeq 
   \Hom_{\wt\CH_{\a}}(\ol M_0^{\nu}, \ol M_0^{\mu})
\end{equation*}
of $R$-modules. 
Let $\bar\vf_{ST}$ be a basis element of $\ol H_{\mu\nu}$, 
where $S = (S^{[1]}, \dots, S^{[g]}) \in \CT_0^{\Bp}(\la, \mu)$ and 
$T = (T^{[1]}, \dots, T^{[g]}) \in \CT_0^{\Bp}(\la, \nu)$ 
for some $\la \in \vL^+$. Then
under the identification in (4.11.1), $\vT$ maps  
$\bar\vf_{ST}$ to 
$\vf_{S^{[1]}T^{[1]}}\otimes\cdots\otimes \vf_{S^{[g]}T^{[g]}}$.
\end{lem}
\begin{proof}
It is enough to show the second assertion since 
$\vf_{S^{[1]}T^{[1]}}\otimes\cdots\otimes\vf_{S^{[g]}T^{[g]}}$
gives a basis of $\Hom_{\wt\CH_{\a}}(\ol M_0^{\nu}, \ol M_0^{\mu})$
under the identification in (4.11.1).
Take $\bar\vf_{ST} \in \ol H_{\mu\nu}$.  Then $\bar\vf_{ST}$ is defined by 
$\bar\vf_{ST}(\ol m_{\nu}) = \ol m_{ST}$. By Lemma 4.10 (i), 
$\ol m_{\nu}$ is mapped to 
$m_{\nu^{[1]}}\otimes\cdots\otimes m_{\nu^{[g]}}$ via $\f_{\nu}$.
$\ol m_{ST}$ is also mapped to 
$m_{S^{[1]}T^{[1]}}\otimes\cdots\otimes m_{S^{[g]}T^{[g]}}$ via
$\f_{\mu}$. Hence via the isomorphism (4.11.1), $\bar\vf_{ST}$ 
corresponds to the $\CH_{\a}$-linear map sending 
$m_{\nu^{[1]}}\otimes\cdots\otimes m_{\nu^{[g]}}$ to 
$m_{S^{[1]}T^{[1]}}\otimes\cdots\otimes m_{S^{[g]}T^{[g]}}$, 
which coincides with 
$\vf_{S^{[1]}T^{[1]}}\otimes\cdots\otimes \vf_{S^{[g]}T^{[g]}}$.
The lemma is proved.
\end{proof}
%%%
\remark{4.13.} \ There exists an $R$-linear map
$\p : \bar\th(\ol H_{\mu\nu}) \to 
     \Hom_{\wt\CH_{\a}}(\ol M_0^{\nu}, \ol M_0^{\mu})$ 
such that $\p\circ\bar\th = \vT$ by Lemma 4.10.  Hence 
$\bar\th$ is injective by Lemma 4.12. However $\bar\th$
is not necessarily surjective.  In Section 7, we describe 
$\Im \bar\th$ in terms of a modified Ariki-Koike algebra.
\para{4.14.}
Let $\vD_{n,g}$ be the set of 
$\a = (n_1, \dots, n_g) \in \ZZ^g_{\ge 0}$ such that 
$n_1 + \dots n_g = n$. 
For $\a \in \vD_{n,g}$, put 
\begin{equation*}
M^{\a} = \bigoplus_{\substack{\mu \in \vL \\
             \a_{\Bp}(\mu) = \a}}M^{\mu}, \qquad 
\ol M_0^{\a} = \bigoplus_{\substack{ \mu \in \vL \\
             \a_{\Bp}(\mu) = \a}}\ol M_0^{\mu}.
\end{equation*}
Then $\CS^{\Bp}_{\a} = \End_{\CH}(M^{\a})$ is a subalgebra 
of $\CS^{\Bp}$, and we have 
$\CS^{\Bp}_{\a} = \bigoplus_{\mu,\nu}H_{\mu\nu}$, where the sum 
is taken over all $\mu, \nu \in \vL$ such that 
$\a_{\Bp}(\mu) = \a_{\Bp}(\nu) = \a$.
Put $\ol\CS^{\Bp}_{\a} = \pi(\CS^{\Bp}_{\a})$.  Then 
$\ol\CS^{\Bp}_{\a}$ is a subalgebra of $\ol\CS^{\BP}$ such that
$\ol\CS^{\Bp}_{\a} = \bigoplus_{\mu,\nu}\ol H_{\mu\nu}$.
Hence we have
\begin{equation*}
\tag{4.14.1}
\ol\CS^{\Bp} = \bigoplus_{\a \in \vD_{n,g}}\ol\CS^{\Bp}_{\a}.
\end{equation*}
On the other hand, Lemma 4.12 implies that 
\begin{equation*}
\tag{4.14.2}
\ol\CS^{\Bp}_{\a} \simeq \End_{\wt\CH_{\a}}(\ol M_0^{\a}).
\end{equation*}
\par
We define an $\CH_{n_k,r_k}$-module $M^{[k]}$ by  
$M^{[k]} = \bigoplus_{\mu^{[k]} \in \vL_{n_k}} M^{\mu^{[k]}}$. 
Define a cyclotomic $q$-Schur algebra $\CS(\vL_{n_k})$ associated
to $\CH_{n_k,r_k}$ by $\CS(\vL_{n_k}) = \End_{\CH_{n_k,r_k}}M^{[k]}$. 
 Then we see that 
\begin{equation*}
\tag{4.14.3}
\End_{\CH_{\a}}(\bigoplus_{\substack{ \mu \in \vL \\
                \a_{\Bp}(\mu) = \a}} M^{\mu^{[1]}}
            \otimes\cdots\otimes M^{\mu^{[g]}})
\simeq \CS(\vL_{n_1})\otimes\cdots\otimes \CS(\vL_{n_g}).
\end{equation*}
\par
The following structure theorem follows from (4.14.1) $\sim$ (4.14.3)
together with (4.11.1).  Note that in the special case where 
$\Bp = (1^r)$, this result was proved in 
[SawS, Theorem 5.5 (i)]
under the assumption that $Q_i - Q_j$ are units in $R$ for any 
$i \ne j$, and that $\vL = \wt\CP_{n,r}(\Bm)$ with $m_i \ge n$ 
for $i = 1, \dots, r$. In our case, we don't need 
any assumption for parameters $Q_i$ nor $\Bm$.
%%%
\begin{thm}  %% Theorem 4.15
Assume that $\vL = \wt\CP_{n,r}(\Bm)$.  Then there
exists an isomorphism of $R$-algebras
\begin{equation*}
\ol\CS^{\Bp}(\vL) \simeq \bigoplus_{\substack{ (n_1, \dots, n_g) \\
     n_1 + \cdots + n_g = n}}\CS(\vL_{n_1})
               \otimes\cdots\otimes\CS(\vL_{n_g}),
\end{equation*}
where $\bar\vf_{ST}$ is mapped to 
$\vf_{S^{[1]T^{[1]}}}\otimes\cdots\otimes\vf_{S^{[g]}T^{[g]}}$.
\end{thm}
For $\la^{[k]}, \mu^{[k]} \in \vL_{n_k}^+$, let 
$W^{\la^{[k]}}$ be the Weyl module, and $L^{\mu^{[k]}}$ be the 
irreducible module with respect to $\CS(\vL_{n_k})$.  As a corollary to the
previous theorem, we have
\begin{cor}  %% Cor 4.16
Assume that $R$ is a field and $\vL$ is as above.  
Let $\la, \mu \in \vL^+$.  
Then under the isomorphism in Theorem 4.15, we have the following.
\begin{enumerate}
\item
$\ol Z^{\la}_{\Bp} \simeq W^{\la^{[1]}}\otimes\cdots\otimes
     W^{\la^{[g]}}$.
\item
$\ol L_{\Bp}^{\mu} \simeq L^{\mu^{[1]}}\otimes\cdots\otimes 
L^{\mu^{[g]}}$.
\item
$[\ol Z^{\la}_{\Bp} : \ol L^{\mu}_{\Bp}]_{\ol\CS^{\Bp}} = \begin{cases}
     \prod_{k=1}^g[ W^{\la^{[k]}}: L^{\mu^{[k]}}]_{\CS(\vL_{n_k})}
                   &\quad\text{ if } \a_{\Bp}(\la) = \a_{\Bp}(\mu), \\
     0             &\quad\text{ otherwise. }
                           \end{cases}$
\end{enumerate}
\end{cor}
Combining this with Theorem 3.13, we have the following product formula
for the decomposition numbers of $\CS(\vL)$, which is a generalization
of [Sa, Corollary 5.10].
%%%
\begin{thm}  %%% Theorem 4.17
Assume that $R$ is a field and that $\vL = \wt\CP_{n,r}(\Bm)$.  
For $\la, \mu \in \vL^+$ such that 
$\a_{\Bp}(\la) = \a_{\Bp}(\mu)$, we have
\begin{equation*}
[W^{\la} : L^{\mu}]_{\CS(\vL)} 
     = \prod_{k=1}^g[W^{\la^{[k]}} : L^{\mu{[k]}}]_{\CS(\vL_{n_k})}.
\end{equation*}
\end{thm}
%%%%%
%%%%%
\section{Modified Ariki-Koike algebra of type $\Bp$}
\para{5.1.}
Throughout this section we assume the following property for 
$\Bm = (m_1, \dots, m_r)$.
\par\medskip\noindent
(5.1.1) \ $m_i \ge n$ for $i = 1, \dots, r$ 
\par\medskip
We keep the assumption that $\vL = \wt\CP_{n,r}(\Bm)$.
Let $\Om = \Om^{\Bp}$ be a subset of $\vL$ consisting of 
$\w = (\w_i^{(j)}) = (\w^{[1]}, \dots, \w^{[g]})$  
satisfying the properties 
\begin{enumerate}
\item
$\w_i^{(j)} \in \{ 0,1\}$, 
\item
$\sum_{j=1}^r\w_i^{(j)} = 1$ for $1 \le i \le n, 1 \le j \le r$,
\item
$\w_i^{(j)} = 0$ unless $j = p_1+r_1, \dots, p_g + r_g$. 
Hence $\w^{[k]} = (-, \dots, -, \w^{(p_k + r_k)})$
for $k = 1, \dots, g$.
\end{enumerate}
Note that $\Om$ coincides with $\Om$ in [SawS, 7.1]
in the case where $\Bp = (1^r)$ (i.e., the case $g = r$). 
While in the case where $\Bp = (r)$ (i.e., the case
$g = 1$), $\Om = \{ \w\}$, where $\w$ is an $r$-partition 
$\w = (-, \dots, -, (1^n))$ which coincides with $\w$ in [M, \S 4]. 
\par
Let $I = \{ 1, \dots, n\}$.  
For $\w \in \Om$, we denote by $I_k$ the set of 
$i$ such that $\w_i^{(p_k+r_k)} = 1$ for $k = 1, \dots, g$. 
Then $I = \coprod_{k=1}^g I_k$ gives a partition of
$I$ into $g$ parts, and thanks to (5.1.1), the set $\Om$ is in bijection 
with the set of partitions of $I$ into $g$ parts.
For $\Ft = (\Ft^{[1]}, \dots, \Ft^{[t]}) \in \Std(\la)$, 
we denote by $I_k$ the letters contained in the standard tableau
$\Ft^{[k]}$.  Then $I = \coprod I_k$ determines 
$\w = \w_{\Ft} \in \Om$.
We associate to $\Ft$ a semi-standard tableau $T$ of 
shape $\la$ as follows; for each $k$ ($1 \le k \le g$), 
the first terms of the 
entries of $T^{(p_k + i)}$  
consist of the entries of $\Ft^{(p_k + i)}$, and the second term of 
them has the common value $p_k + r_k$ for $i = 1, \dots, r_k$.
Then $T \in \CT_0^{\Bp}(\la, \w)$, and any element of 
$\CT_0^{\Bp}(\la, \w)$ is obtained from $\Ft \in \Std(\la)$ such that 
$\w = \w_{\Ft}$ by the above procedure.
The correspondence $\Ft \mapsto T$ gives a bijective correspondence
\begin{equation*}
\tag{5.1.2}
\Std(\la) \simeq \bigcup_{\w \in \Om}\CT_0^{\Bp}(\la, \w).
\end{equation*}
We denote by $\Std(\la)_{\w}$ the subset of $\Std(\la)$ 
corresponding to $\CT^{\Bp}_0(\la, \w)$ under the bijection (5.1.2), 
i.e., $\Std(\la)_{\w} = \{ \Ft \in \Std(\la) \mid \w_{\Ft} = \w\}$.
\par
Assume that $\w \in \Om$ corresponds to the partition 
$I = \coprod_{k}I_k$, where $\Ba_{\Bp}(\w) = (a_1, \dots, a_g)$.  
We write $I_k$ as 
$I_k = \{ i_{k1} < i_{k2} < \dots < i_{kn_k}\}$.  
We define $d(\w) \in \FS_n$ as 
\begin{equation*}
d(\w) = \begin{pmatrix}
         \dots & a_k + 1& a_k+2 & \dots &a_k + n_k & \dots  \\ 
         \dots & i_{k1} & i_{k2} & \dots &i_{kn_k}   & \dots                
        \end{pmatrix}.
\end{equation*}
Suppose that $T \in \CT_0^{\Bp}(\la, \w)$ corresponds to 
$\Ft \in \Std(\la)$ via (5.1.2).   Let $\Ft_1 \in \Std(\la)$ 
be such that $\Ft = \Ft_1 d(\w)$.  Then the letters contained 
in $\Ft_1^{[k]}$ consist of $\{ a_k+1, \dots, a_k + n_k\}$,
and $\Ft_1$ is the unique element in $\Std(\la)$ such that 
$\w(\Ft_1) = T$.   
In particular, assume that  
$S \in \CT_0^{\Bp}(\la, \mu), T \in \CT_0^{\Bp}(\la, \w)$, 
for $\mu \in \vL, \w \in \Om$, and that $\Ft \in \Std(\la)$
corresponds to $T$ via (5.1.2).  Then we have
\begin{equation*}
\tag{5.1.3}
m_{ST}T_{d(\w)} = m_{S\Ft}.
\end{equation*}
\para{5.2.}
For each $\mu \in \vL$, 
let $\vf_{\mu}$ be the identity map on $M^{\mu}$. By 2.4, 
$\vf_{\mu} \in H_{\mu\mu}$, and we put 
$\bar\vf_{\mu} = \pi(\vf_{\mu})\in \ol H_{\mu\mu}$. 
If we put $\bar\vf_{\Om} = \sum_{\w \in \Om}\bar\vf_{\w}$, 
$\bar\vf_{\Om}$ is an idempotent in $\ol\CS^{\Bp}$, and we define
a subalgebra $\ol\CH^{\Bp}$ of $\ol\CS^{\Bp}$ by  
$\ol\CH^{\Bp} = \bar\vf_{\Om}\ol\CS^{\Bp}\bar\vf_{\Om}$.
We call $\ol\CH^{\Bp}$ the modified Ariki-Koike algebra of type 
$\Bp$.  In the case where $\Bp = (1^r)$, $\ol\CH^{\Bp}$
can be identified with the modified Ariki-Koike algebra given in [SawS]
(see 7.1 in [loc. cit.]).
One can write $\ol\CH^{\Bp} = \bigoplus_{\w,\w' \in \Om}\ol H_{\w\w'}$.
In particular, $\ol\CH^{\Bp}$ has an $R$-free basis 
\begin{equation*}
\tag{5.2.1}
\CB^{\Bp} = \{ \ol\vf_{ST} \mid S \in \CT_0^{\Bp}(\la, \w),  
          T \in \CT_0^{\Bp}(\la, \w') \text{ for } 
               \w, \w' \in \Om, \la \in \vL^+\}.
\end{equation*}
Note that each $\ol\vf_{ST} \in \CB^{\Bp}$ determines uniquely the pair
$\Fs, \Ft$ of standard tableau of shape $\la$ by (5.1.2).  We denote 
$\ol\vf_{ST}$ by $m^{\Bp}_{\Fs\Ft}$ if $S, T$ correspond to 
$\Fs,\Ft \in \Std(\la)$.  Thus we see that 
\begin{equation*}
\tag{5.2.2}
\CB^{\Bp} = \{ m_{\Fs\Ft}^{\Bp} \mid \Fs, \Ft \in \Std(\la) 
                  \text{ for some } \la \in \vL^+\}.
\end{equation*}
Note that $\ol\CS^{\Bp}$ has a structure of the cellular algebra 
with the cellular basis $\ol\ZC^{\Bp}$.  Since the involution 
$*$ on $\ol\CS^{\Bp}$ stabilizes the set $\CB^{\Bp}$, we see that 
\par\medskip\noindent
(5.2.3) \ $\ol\CH^{\Bp}$ is a cellular algebra with the cellular
basis $\CB^{\Bp}$.
\par\medskip
More generally, we consider for each $\mu \in \vL$ an $R$-submodule
$\bar\vf_{\mu}\ol\CS^{\Bp}\bar\vf_{\Om}$ of $\ol\CS^{\Bp}$.
Then $\bar\vf_{\mu}\ol\CS^{\Bp}\bar\vf_{\Om}$ has an $R$-basis 
\begin{equation*}
\{ \ol\vf_{ST} \mid S \in \CT_0^{\Bp}(\la,\mu), 
              T \in \CT_0^{\Bp}(\la, \w) 
     \text{ for } \w \in \Om, \la \in \vL^+ \}.
\end{equation*}
\par
Let $\ol M^{\Om} = \bigoplus_{\w \in \Om}\ol M^{\w}$, 
and put $\ol m_{\Om} = \sum_{\w \in \Om}\ol m_{\w}T_{d(\w)} 
          \in \ol M^{\Om}$.   
Then for $S \in \CT_0^{\Bp}(\la, \mu), T \in \CT_0^{\Bp}(\la,\w)$, 
we have
\begin{equation*}
\ol\vf_{ST}(\ol m_{\Om}) = \ol\vf_{ST}(\ol m_{\w}T_{d(\w)}) 
                         = \ol m_{ST}T_{d(\w)} = \ol m_{S\Ft} 
\end{equation*}
by (5.1.3), where $\Ft \in \Std(\la)$ corresponds to $T$ via 
(5.1.2).  Since 
$\{ \ol m_{S\Ft} \mid 
            S \in \CT_0^{\Bp}(\la, \mu), \Ft \in \Std(\la)\}$
gives a basis of $\ol M^{\mu}$, we see that the map 
$\vf \mapsto \vf(\ol m_{\Om})$ gives an isomorphism of 
$R$-modules 
\begin{equation*}
\tag{5.2.4}
\bar\vf_{\mu}\ol\CS^{\Bp}\bar\vf_{\Om} \simeq \ol M^{\mu}, \quad 
         \ol\vf_{ST} \lra \ol m_{S\Ft}.
\end{equation*}
Since $\bar\vf_{\Om}\ol\CS^{\Bp}\bar\vf_{\Om} = \ol\CH^{\Bp}$ acts naturally 
on $\bar\vf_{\mu}\ol\CS^{\Bp}\bar\vf_{\Om}$ from the right, one can define a
right action of $\ol\CH^{\Bp}$ on $\ol M^{\mu}$ through (5.2.4). 
Let $\mu, \nu \in \vL$.  By 4.9, we know that 
$\vf \in \ol H_{\mu\nu}$ gives a map $\ol\th(\vf)$ 
from $\ol M^{\nu}$ to $\ol M^{\mu}$. It is clear by definition, that 
$\bar\th(\vf)$ commutes with the action of $\ol\CH^{\Bp}$. Hence 
we have an $R$-linear map 
$\th': \ol H_{\mu\nu} \to \Hom_{\ol\CH^{\Bp}}
                         (\ol M^{\nu}, \ol M^{\mu})$, 
which induces an $R$-algebra homomorphism 
$\th': \ol\CS^{\Bp} \to \End_{\ol\CH^{\Bp}}(\ol M)$,
where $\ol M = \bigoplus_{\mu \in \vL}\ol M^{\mu}$.
\par
The following result is a generalization of Proposition 7.5 
in [SawS]. 
%%%
\begin{prop}  %%% Prop. 5.3.  
For each $\a = (n_1, \dots, n_g) \in \vD_{n,g}$, 
put $n_{\a} = n!/n_1!\cdots n_g!$.  Then 
we have an isomorphism of $R$-algebras
\begin{equation*}
\ol\CH^{\Bp} \simeq \bigoplus_{\a \in \vD_{n,g}}
     M_{n_{\a}}(\CH_{n_1,r_1}\otimes\cdots\otimes \CH_{n_g,r_g}).
\end{equation*}
\end{prop}
\begin{proof}
By (4.14.1), one can write 
\begin{equation*}
\ol\CH^{\Bp} = 
\bar\vf_{\Om}\ol\CS^{\Bp}\bar\vf_{\Om} = \bigoplus_{\a \in \vL_{n,g}}
           \bar\vf_{\Om,\a}\ol\CS^{\Bp}_{\a}\bar\vf_{\Om,\a}.
\end{equation*}
Here $\bar\vf_{\Om,\a} = \sum_{\w}\bar\vf_{\w}$ is an idempotent of 
$\ol\CS^{\Bp}_{\a}$, where the sum is taken over all $\w \in \Om$
such that $\a_{\Bp}(\w) = \a$.
We define a subalgebra $\ol\CH^{\Bp}_{\a}$ of $\ol\CH^{\Bp}$ by 
$\ol\CH^{\Bp}_{\a} = \bar\vf_{\Om,\a}\ol\CS^{\Bp}_{\a}\bar\vf_{\Om,\a}$.
Put $\ol M_0^{\Om,\a} = \ol M^{\Om} \cap \ol M_0^{\a}$.  Then 
by (4.14.2) we have
\begin{equation*}
\tag{5.3.1}
\ol\CH^{\Bp}_{\a} \simeq 
                \End_{\wt\CH_{\a}}(\ol M^{\Om,\a}_0) = 
                \bigoplus_{\substack{\w,\w' \in \Om \\
    \a_{\Bp}(\w) = \a_{\Bp}(\w') = \a}}
            \Hom_{\wt\CH_{\a}}(\ol M_0^{\w}, \ol M_0^{\w'}).
\end{equation*}
Now the $\wt\CH_{\a}$-module $\ol M^{\w}_0$ is isomorphic to  
   the $\CH_{\a}$-module 
$M^{\w^{[1]}}\otimes\cdots\otimes M^{\w^{[g]}}$ by 
Corollary 4.8.  In our case $M^{\w^{[k]}} = \CH_{n_k,r_k}$
(see 5.1).  Hence for any $\w, \w' \in  \Om$ such that 
$\a_{\Bp}(\w) = \a_{\Bp}(\w') = \a$,  we have
\begin{align*}
\tag{5.3.2}
\Hom_{\wt\CH_{\a}}(\ol M_0^{\w}, \ol M_0^{\w'}) &\simeq
\End_{\CH_{\a}}(\CH_{n_1,r_1}\otimes\cdots\otimes \CH_{n_g,r_g}) \\
 &\simeq \CH_{n_1,r_1}\otimes\cdots\otimes \CH_{n_g,r_g}.
\end{align*}
The proposition follows from this by noticing that 
$\sharp\{ \w \in \Om \mid \a_{\Bp}(\w) = \a\} = n_{\a}$.
\end{proof}
\para{5.4.}
By $\bar \th$, $\ol\CS^{\Bp}$ acts on $\ol M$ from the left, and which
commutes with the right action of $\CH$.  
Hence we have a homomorphism
$\r : \CH \to \End^0_{\ol\CS^{\Bp}}\ol M$ (see Notation).
Since $\sum_{\mu \in \vL}\bar\vf_{\mu} = \Id_{\ol M}$, we have
$\ol\CS^{\Bp}\bar\vf_{\Om} \simeq \ol M$ by (5.2.4).  This implies 
a natural isomorphism of $R$-algebras
\begin{equation*}
\tag{5.4.1}
\End^0_{\ol\CS^{\Bp}}\ol M \simeq 
         \End^0_{\ol\CS^{\Bp}}(\ol\CS^{\Bp}\bar\vf_{\Om}) 
       \simeq \bar\vf_{\Om}\ol\CS^{\Bp}\bar\vf_{\Om} = \ol\CH^{\Bp},
\end{equation*}
where the second isomorphism is given by $f \mapsto f(\bar\vf_{\Om})$
for $f \in \End_{\ol\CS^{\Bp}}(\ol\CS^{\Bp}\bar\vf_{\Om})$.
It follows that we have a homomorphism 
$\r_0 : \CH \to \ol\CH^{\Bp}$ of $R$-algebras thorough 
$\CH \to \End^0_{\ol\CS^{\Bp}} \ol M$. 
The homomorphism $\r_0$ is explicitly given as follows; 
we have $\ol\CH^{\Bp} = \bar\vf_{\Om}\ol\CS^{\Bp}\bar\vf_{\Om} 
        \simeq \ol M^{\Om}$ via 
$\vf \mapsto \vf(\ol m_{\Om})$. Then for each $h \in \CH$, there 
exists a unique $\vf_{h} \in \ol\CH^{\Bp}$ such that 
$\vf_h(\ol m_{\Om}) = \ol m_{\Om}h \in \ol M^{\Om}$. 
The map $h \mapsto \vf_h$ gives $\r_0$.
\par 
Now $\ol\CH^{\Bp}$-module $\ol M$ is regarded as an $\CH$-module
via $\r_0$, which coincides with the original $\CH$-module $\ol M$.
It follows that we have an injection
\begin{equation*}
 \Hom_{\ol\CH^{\Bp}}(\ol M^{\nu}, \ol M^{\mu})
    \hra \Hom_{\CH}(\ol M^{\nu}, \ol M^{\mu}),
\end{equation*}
and $\th$ factors through $\th'$ via this injection.
Since $\bar\th$ is injective by Remark 4.13, we see that 
\par\medskip\noindent
(5.4.2) \ The map $\th': \ol H_{\mu\nu} \to 
   \Hom_{\ol\CH^{\Bp}}(\ol M^{\nu}, \ol M^{\mu})$ is injective.
\par\medskip
Since $\ol M^{\mu}$ is generated by $\ol m_{\mu}$ as an $\CH$-module,
it is generated by $\ol m_{\mu}$ as an $\ol\CH^{\Bp}$-module, i.e.,
we have $\ol M^{\mu} = \ol m_{\mu}\ol \CH^{\Bp}$.
The following lemma is also clear from the fact 
that $\ol\CH^{\Bp} \simeq \ol M^{\Om}$ via 
$ \vf \mapsto \vf(\ol m_{\Om})$ as noticed above. 
%%%%
\begin{lem}  %%% Lemma 5.5
We have $\ol M^{\Om} = \ol m_{\Om}\ol\CH^{\Bp}$.  The map
$h \mapsto \ol m_{\Om}h$ gives an isomorphism of $R$-modules
$\ol\CH^{\Bp} \to \ol M^{\Om}$, namely $\ol M^{\Om}$ is the 
regular representation of $\ol \CH^{\Bp}$.
\end{lem}

\section{Presentation for $\ol\CH^{\Bp}$}
\para{6.1}
We shall define several elements in $\ol\CH^{\Bp}$, and show 
that they generate $\ol\CH^{\Bp}$. 
For each $\w \in \Om$ let $I = \coprod I_k$ be the corresponding 
partition of $I$.  Define a map $b_{\w} : I \to \ZZ_{>0}$ by 
$b_{\w}(i) = k$ if $i \in I_k$.  We put $Q^{\Bp}_k = Q_{p_k + r_k}$ 
for $k = 1, \dots, g$.  Under this notation, we define elements 
$\xi_i \in \ol\CS^{\Bp}$, for $i = 1, \dots, n$, by 
\begin{equation*}
\tag{6.1.1}
\xi_i = \sum_{\w \in \Om}Q^{\Bp}_{b_{\w}(i)}\bar\vf_{\w}.
\end{equation*}
Clearly, $\bar\vf_{\Om}\xi_i\bar\vf_{\Om} = \xi_i$, and so
$\xi_1, \dots, \xi_n$ are elements in $\ol\CH^{\Bp}$. They 
commute each other.   Moreover, they satisfy the relation
\begin{equation*}
\tag{6.1.2}
(\xi_j - Q^{\Bp}_1)(\xi_j - Q^{\Bp}_2)\cdots (\xi_j - Q^{\Bp}_g) = 0
\end{equation*}
for $j = 1, \dots, n$.
\par
Under the isomorphism in (5.2.4), the action of $\xi_i$ on 
the basis element $\ol m_{S\Ft}$ in $\ol M^{\mu}$ is given as follows. 
\begin{equation*}
\tag{6.1.3}
\ol m_{S\Ft}\xi_i = Q^{\Bp}_{b_{\w}(i)}\ol m_{S\Ft} \quad\text{ if }
                         \Ft \in \Std(\la)_{\w},
\end{equation*}
where $\Std(\la)_{\w}$ is as in 5.1.  
Note that in this case $b_{\w}(i)$ coincides with $k$ such that
the letter $i$ is contained in $\Ft^{[k]}$.
By [DJM, Proposition 3.18], $\ol m_{\mu}$ is written, for
$\mu \in \vL$,  as a linear
combination of $\ol m_{\Fs\Ft}$ such that the letters contained 
in the $k$ component of $\Ft$ is the same as that of $\Ft^{\mu}$.
It follows from this, by making use of (6.1.3), that 
\begin{equation*}
\tag{6.1.4}
\ol m_{\mu}\xi_i = Q^{\Bp}_{b(i)}\ol m_{\mu},
\end{equation*}
where $b(i) = k$ if $a_k+1 \le i \le a_k + n_k$
under the notation $\Ba_{\Bp}(\mu) = (a_1, \dots, a_g)$ and 
$\a_{\Bp}(\mu) = (n_1, \dots, n_g)$.
\par
Let $\r_0 : \CH \to \ol \CH^{\Bp}$ be the homomorphism 
defined in 5.4.  We note that 
\par\medskip\noindent
(6.1.5) \ The restriction of $\r_0$ on $\CH_n$ is injective.
\par\medskip
In fact, it is enough to show that $\r_0(T_w)$ ($w \in \FS_n$) are
linearly independent as operators on $\ol M$.
Now $\ol M = \bigoplus_{\a \in \vD_{n,g}}\ol M^{\a}$, 
and $T_w$ preserves the 
subspaces $\ol M^{\a}$.  We choose $\a$ such that 
$\a = (n, 0, \dots, 0)$.  Then $\CH_n$ is contained in 
$\wt\CH_{\a} = \CH$, and $\r_0(T_w)$ induces an operator on 
$\ol M^{\a}_0$.  By our choice of $\a$, Corollary 4.8 implies that 
$\ol M_0^{\a}$ can be identified with $M'$, 
the $\CH_{n, r_1}$-module corresponding to $M$ for $\CH$, and the 
action of $\wt\CH_{\a}$ on $\ol M_0^{\a}$ coincides with the action of 
$\CH_{n, r_1}$ on $M'$. In particular, the action of $\r_0(T_w)$ on
$\ol M_0^{\a}$ corresponds to the action of $T_w$ on $M'$ 
(we regard $T_w \in \CH_n \subset \CH_{n,r_1}$).  Since $T_w$ 
($w \in \FS_n$) are linearly independent as operators on $M'$, 
we see that $\r_0(T_w)$ are linearly independent as asserted.   
\par
By (6.1.5), we regard $\CH_n$ as a subalgebra of $\ol\CH^{\Bp}$, 
and define the elements $T_1, \dots, T_{n-1} \in \ol\CH^{\Bp}$ 
by the generators of $\CH_n$. 
\para{6.2.}
We shall determine the commutation relations between 
$T_j$ and $\xi_k$.  In view of Lemma 5.5, we compare the 
elements $\ol m_{\Om}T_j\xi_k$ and $\ol m_{\Om}\xi_kT_j$.
First 
we compute the element $\ol m_{\Om}T_j$ 
for $T_j \in \CH_n$.
Since $\ol m_{\Om}T_j = \sum_{\w \in \Om}\ol m_{\w}T_{d(\w)}T_j$, 
we compute $m_{\w}T_{d(\w)}T_j$.
Let $I = \coprod I_k$ be the partition corresponding to $\w$.
Assume that $j \in I_k$ and $j+1 \in I_{k'}$. Then we see that 
\begin{equation*}
T_{d(\w)}T_j = \begin{cases}
                 T_{d(\w)s_j}  &\quad\text{ if } k \le k',  \\
                 T_{d(\w)s_j} + (q-q\iv)T_{d(\w)}
                               &\quad\text{ if }  k > k', 
              \end{cases}
 \end{equation*}
where $s_j$ is the element in $\FS_n$ corresponding to $T_j$.
Note that $m_{\w} = u_{\Ba}^+ = m_{\la}$, where $\la$ is 
the multi-partition 
obtained from $\w$ by rearranging the rows.
Put $\Ft_{\w} = \Ft^{\la}d(\w) \in \Std(\la)$.  Put 
$\Fv_{\w} = \Ft_{\w}s_j$.  If $k \ne k'$, then $\Fv_{\w} \in \Std(\la)$ and
it is expressed as $\Ft_{\w'}$, where $\w' \in \Om$ is obtained from 
$\w$ by exchanging $j$ and $j+1$ in $I_k$ and $I_{k'}$.
One can write $m_{\w}T_{d(\w)} = m_{S_{\w}\Ft_{\w}}$ and 
$m_{\w}T_{d(\w)s_j} = m_{S_{\w}\Fv_{\w}}$, where 
$S_{\w} = \w(\Ft^{\la}) \in \CT_0^{\Bp}(\la, \w)$.  
Hence we have 
\begin{equation*}
\tag{6.2.1}
\ol m_{\w}T_{d(\w)}T_j = \begin{cases}
              \ol m_{S_{\w}\Fv_{\w}}
                                  &\quad\text{ if } k = k', \\
              \ol m_{S_{\w}\Ft_{\w'}}
                                  &\quad\text{ if } k < k', \\
      \ol m_{S_{\w}\Ft_{\w'}} + (q-q\iv)\ol m_{S_{\w}\Ft_{\w}}
                                  &\quad\text{ if } k > k'.
                         \end{cases}            
\end{equation*}
Note that in the first case, by [DJM, Proposition 3.18], 
$\ol m_{S_{\w}\Fv_{\w}}$ is expressed 
as a linear combination of basis elements $\ol m_{S'\Fv}$ such that 
$\w_{\Fv} = \w$. 
It follows from (6.2.1) that 
\begin{equation*}
\begin{split}
\ol m_{\Om}T_j &= \sum_{\substack{\w \in \Om \\
                        b_{\w}(j) < b_{\w}(j+1)}}\ol m_{S_{\w}\Ft_{\w'}} \\
    &+ \sum_{\substack{\w \in \Om \\ 
                        b_{\w}(j) = b_{\w}(j+1)}}\ol m_{S_{\w}\Fv_{\w}}
    + \sum_{\substack{\w \in \Om \\
                        b_{\w}(j) > b_{\w}(j+1)}}
        (\ol m_{S_{\w}\Ft_{\w'}} + (q-q\iv)\ol m_{S_{\w}\Ft_{\w}}),
\end{split}
\end{equation*}
where $\w' \in \Om$ is obtained from $\w$ by $s_j$ as above, and 
$\Fv_{\w} = \Ft_{\w}s_j$.  Thus by (6.1.3) and (6.1.4),  
we have
\begin{equation*}
\begin{split}
\ol m_{\Om}T_j\xi_k &= \sum_{\substack{\w \in \Om \\
                              b_{\w}(j) \ne b_{\w}(j+1)}}
                 Q^{\Bp}_{b_{\w'}(k)}\ol m_{S_{\w}\Ft_{\w'}} \\
&+ \sum_{\substack{ \w \in \Om \\  b_{\w}(j) = b_{\w}(j+1)}}
                Q^{\Bp}_{b_{\w}(k)}\ol m_{S_{\w}\Fv_{\w}}
 + \sum_{\substack{ \w \in \Om \\ b_{\w}(j) > b_{\w}(j+1)}}
                 Q^{\Bp}_{b_{\w}(k)} (q-q\iv)\ol m_{S_{\w}\Ft_{\w}}.  
\end{split}
\end{equation*}
On the other hand, we have
\begin{equation*}
\begin{split}
\ol m_{\Om}\xi_kT_j &= \sum_{\substack{\w \in \Om \\
                              b_{\w}(j) \ne b_{\w}(j+1)}}
                 Q^{\Bp}_{b_{\w}(k)}\ol m_{S_{\w}\Ft_{\w'}} \\
&+ \sum_{\substack{ \w \in \Om \\  b_{\w}(j) = b_{\w}(j+1)}}
                Q^{\Bp}_{b_{\w}(k)}\ol m_{S_{\w}\Fv_{\w}}
 + \sum_{\substack{ \w \in \Om \\ b_{\w}(j) > b_{\w}(j+1)}}
                 Q^{\Bp}_{b_{\w}(k)} (q-q\iv)\ol m_{S_{\w}\Ft_{\w}}.  
\end{split}
\end{equation*}
It follows that 
\begin{equation*}
\tag{6.2.2}
\ol m_{\Om}(T_j\xi_k - \xi_kT_j)
 = \sum_{\substack{\w \in \Om \\ 
        b_{\w}(j) \ne b_{\w}(j+1)}}
   (Q^{\Bp}_{b_{\w'}(k)} - 
            Q^{\Bp}_{b_{\w}(k)})\ol m_{S_{\w}\Ft_{\w'}}.
\end{equation*}
Note that if $k \ne j, j+1$, then $b_{\w}(k) = b_{\w'}(k)$
for any $\w$. It follows that 
\par\medskip\noindent
(6.2.3) \ $T_j\xi_k = \xi_kT_j$ if $k \ne j, j+1$.
\para{6.3.}
Let $A$ be a square matrix of degree $g$ whose $ij$-entry 
is given by $(Q^{\Bp}_j)^{i-1}$ for $1 \le i,j \le g$.
Thus $A$ is the Vandermonde matrix, and 
$\vD = \det A = \prod_{i>j}(Q^{\Bp}_i - Q^{\Bp}_j)$.
We pose the following assumption so that $\vD\iv \in R$.
\par\medskip\noindent
(6.3.1) \ $Q^{\Bp}_i - Q^{\Bp}_j$ are units in $R$ for any $i \ne j$.
\par\medskip
We express $A\iv = \vD\iv B$ with
$B = (h_{ij})$ for $h_{ij} \in R$.  We define a polynomial
$F_i(X) \in R[X]$,  for $1 \le i \le g$, by
\begin{equation*}
F_i(X) = \sum_{j = 1}^g h_{ij}X^{j-1}.
\end{equation*}
\par
\par\medskip
We denote by $\Om_j^{[c]}$ the set of $ \w \in \Om$ such that
$b_{\w}(j) = c$ for $1 \le j \le n, 1 \le c \le g$.
As in 6.2, one can write 
$\ol m_{\Om} = \sum_{\w \in \Om}\ol m_{S_{\w}\Ft_{\w}}$, and so
\begin{equation*}
\tag{6.3.2}
\ol m_{\Om}\xi_j^b = 
   \sum_{\w \in \Om}(Q^{\Bp}_{b_{\w}(j)})^b\ol m_{S_{\w}\Ft_{\w}}
  = \sum_{c = 1}^g(Q_c^{\Bp})^b
            \sum_{\w \in \Om_j^{[c]}}\ol m_{S_{\w}\Ft_{\w}}
\end{equation*}
for $b = 0, \dots, g-1$.
We regard (6.3.2) as a system of linear equations with 
unknown variables $\sum_{\w \in \Om^{[c]}_j}\ol m_{S_{\w}\Ft_{\w}}$.
Since $\vD\iv \in R$, we see that 
\begin{equation*}
\sum_{\w \in \Om^{[c]}_j}\ol m_{S_{\w}\Ft_{\w}} = 
     \ol m_{\Om}\cdot\vD\iv\sum_{b = 1}^gh_{cb}\xi_j^{b-1} 
       = \ol m_{\Om}\cdot\vD\iv F_c(\xi_j).
\end{equation*}
Repeating a similar procedure, we have 
\begin{equation*}
\tag{6.3.3}
\sum_{\w \in \Om_{j}^{[c_1]} \cap \Om_{j+1}^{[c_2]}}\ol m_{S_{\w}\Ft_{\w}}
   = \ol m_{\Om}\cdot\vD^{-2}F_{c_1}(\xi_j)F_{c_2}(\xi_{j+1}).
\end{equation*}
By applying $T_j$ on both side of (6.3.3), and by using (6.2.1), we
have
\begin{equation*}
\tag{6.3.4}
\begin{split}
\sum_{\w \in \Om_j^{[c_1]} \cap \Om_{j+1}^{[c_2]}}\ol m_{S_{\w}\Ft_{\w'}}
  = \begin{cases}
        \ol m_{\Om}\cdot \vD^{-2}F_{c_1}(\xi_j)F_{c_2}(\xi_{j+1})T_j
              &\quad\text{ if }c_1 < c_2, \\
        \ol m_{\Om}\cdot \vD^{-2}F_{c_1}(\xi_j)F_{c_2}(\xi_{j+1})
                       (T_j - (q - q\iv))
              &\quad\text{ if } c_1 > c_2. 
     \end{cases}
\end{split}
\end{equation*}
We show the following lemma, which is analogous to [Sh, Lemma 3.4].
%%%
\begin{lem}  %%% Lemma 6.4
For $j = 1, \dots, n-1$, we have
\begin{align*}
T_j\xi_{j+1} &= \xi_jT_j +  \vD^{-2}\sum_{c_1 > c_2}
           (Q^{\Bp}_{c_2} - Q^{\Bp}_{c_1})(q - q\iv)
                               F_{c_1}(\xi_j)F_{c_2}(\xi_{j+1}), \\
T_j\xi_j &= \xi_{j+1}T_j - \vD^{-2}\sum_{c_1 > c_2}
           (Q^{\Bp}_{c_2} - Q^{\Bp}_{c_1})(q - q\iv)
                               F_{c_1}(\xi_j)F_{c_2}(\xi_{j+1}), \\
T_j\xi_k &= \xi_kT_j \quad (k \ne j, j+1).
\end{align*}
\end{lem}
\begin{proof}
The third formula is already shown in (6.2.3).  So 
assume that $k = j$ or $j+1$. 
Substituting (6.3.4) into (6.2.2), 
and by using Lemma 5.5, we have
\begin{equation*}
\tag{6.4.1}
\begin{split}
T_j\xi_k - \xi_kT_j = 
  &\ve\vD^{-2}\biggl\{\sum_{c _1 < c_2}(Q^{\Bp}_{c_2} - Q^{\Bp}_{c_1})
        F_{c_1}(\xi_j)F_{c_2}(\xi_{j+1})T_j \\
  &+ \sum_{c_1 > c_2}(Q^{\Bp}_{c_2} - Q^{\Bp}_{c_1})
    F_{c_1}(\xi_j)F_{c_2}(\xi_{j+1})(T_j - (q - q\iv))\biggr\},
\end{split}
\end{equation*}
where $\ve = 1$ (resp. $\ve = -1$) if $k = j$ (resp. $k = j+1$).
\par
We note that the following formula holds.
\begin{equation*}
\tag{6.4.2}
\xi_{j+1} - \xi_j = 
             \vD^{-2}\sum_{c_1 < c_2}(Q^{\Bp}_{c_2} - Q^{\Bp}_{c_1})
                \biggl\{ F_{c_1}(\xi_j)F_{c_2}(\xi_{j+1})
                       - F_{c_2}(\xi_j)F_{c_1}(\xi_{j+1})\biggr\}
\end{equation*}
In fact it is enough to compare the values at 
$\ol m_{S_{\w}\Ft_{\w}} \in \ol M^{\Om}$.  This is essentially the 
same as the case where $\Bp = (1^r)$, and in that case 
the formula is proved in [Sh, (3.4.2)].
\par
Now (6.4.1) can be written, by making use of (6.4.2), as 
\begin{equation*}
\begin{split}
T_j\xi_k - \xi_kT_j = \ve&(\xi_{j+1} - \xi_j)T_j \\
    & - \ve\sum_{c_1 > c2 }(Q^{\Bp}_{c_2} - Q^{\Bp}_{c_1})(q - q\iv)
               F_{c_1}(\xi_j)F_{c_2}(\xi_{j+1}). 
\end{split}
\end{equation*}
The first and the second equalities in the lemma follow from this.
\end{proof}
\para{6.5.}
For each $\a \in \vD_{n,g}$ and for $k = 1, \dots, g$, 
we define $T^{[k]}_{\a, 0} \in \ol\CH^{\Bp}$ as follows.
We regard $T^{[k]}_0 \in \CH_{n_k,r_k}$ as an element in 
$\CH_{n_1,r_1}\otimes\cdots\otimes \CH_{n_g,r_g}$, and  
we denote by $T^{[k]}_{\a,0} \in \ol\CH^{\Bp}_{\a}$ 
the diagonal matrix consisting of 
$T_0^{[k]}$ in the diagonal entries under the isomorphism in 
Proposition 5.3.
In particular, one can write 
\begin{equation*}
\ol m_{\Om} T^{[k]}_{\a,0} = 
     \sum_{ \w \in \Om^{\a}}\ol m_{\w}T_{d(\w)}L_{a_k+1},
\end{equation*}
where $\Om^{\a} = \{ \w \in \Om \mid \a_{\Bp}(\w) = \a\}$.
Thus we see that $T^{[k]}_{\a,0}$ acts on 
$\ol M^{\Om,\a} = \ol M^{\Om}\cap \ol M^{\a}$ 
as $L_{a_k+1}$,  and annihilates $\ol M^{\Om,\a'}$ for 
any $\a' \ne \a$.  
\par
For a given $\w \in \Om$, put $c_i = b_{\w}(i)$ for 
$i = 1, \dots, n$.  We define $F_{\w}(\xi) \in \ol\CH^{\Bp}$ by 
\begin{equation*}
\tag{6.5.1}
F_{\w}(\xi) = F_{c_1}(\xi_1)F_{c_2}(\xi_2)\cdots F_{c_n}(\xi_n). 
\end{equation*}
We have the following lemma.
%%%
\begin{lem}   %%% Lemma 6.6
Under the assumption of (6.3.1), the elements
\begin{equation*}
\xi_i\  (1 \le i \le n),  \quad T_j\  (1 \le j \le n-1), 
\quad  T^{[k]}_{\a, 0}\ (\a \in \vD_{n,g}, 1 \le k \le g)     
\end{equation*}
generate $\ol\CH^{\Bp}$.
\end{lem}
\begin{proof}
Let $\CK$ be the subalgebra of $\ol\CH^{\Bp}$ generated 
by elements in the lemma. In view of Lemma 5.5, it is enough 
to show that $\ol M^{\Om} = \ol m_{\Om}\CK$.  
First we show that 
\begin{equation*}
\tag{6.6.1}
\ol m_{\w} \in \ol m_{\Om}\CK
\end{equation*}
for any $\w \in \Om$.
In fact, 
we have 
$\bigcap_{i=1}^n\Om^{[c_i]}_i = \{ \w \}$ with $c_i = b_{\w(i)}$. 
Hence by repeating the argument used to prove (6.3.3), we see that
\begin{equation*}
\tag{6.6.2}
\ol m_{\w}T_{d(\w)} = \ol m_{S_{\w}\Ft_w}
     = \ol m_{\Om}\cdot \vD^{-n}F_{\w}(\xi).
\end{equation*}
This implies that $\ol m_{\w}T_{d(\w)} \in \ol m_{\Om}\CK$.
Since $T_{d(\w)}$ is an invertible element in $\CK$, 
we obtain (6.6.1). 
\par
Now take $\ol m_{\w}$ and put $\a = \a_{\Bp}(\w)$.
We know that $\ol m_{\w} \in \ol M_0^{\w}$, and that 
$\ol M_0^{\w} = \ol m_{\w}\wt\CH_{\a}$ (see the proof of Lemma 4.10).
Note that $\wt H_{\a}$ is generated by $L_{a_k +1}$ and 
$\wt H_{\a} \cap \CH_n$, and the action of $L_{a_k+1}$ on 
$\ol M_0^{\a}$ coincides with that of $T^{[k]}_{\a,0}$. 
It follows that $\ol M_0^{\w} = \ol m_{\w}\wt\CH_{\a} 
                 \subset \ol m_{\Om}\CK$.
Here $\ol M_0^{\w}$ has the basis $\{\ol m_{S\Ft}\}$ with 
$S \in \CT_0^{\Bp}(\la,\w)$ and $\Ft \in \Std(\la)_0$.
While the basis of $\ol M^{\w}$ is given by $\{\ol m_{S\Ft'}\}$
for $S \in \CT^{\Bp}(\la, \w)$ and $\Ft' \in \Std(\la)$.
If we take $\Ft = \Ft^{\la} \in \Std(\la)_0$, any $\Ft'$ is 
obtained as $\Ft' = \Ft d(\Ft')$, and we have 
$\ol m_{S\Ft'} = \ol m_{S\Ft}T_{d(\Ft')}$. 
It follows that 
$\ol M^{\w} \subseteq \ol M_0^{\w}\CH_n 
                     \subseteq \ol m_{\Om}\CK$ for any 
$\w \in \Om$, and so $\ol M^{\Om} = \ol m_{\Om}\CK$.
The lemma is proved.
\end{proof}
\para{6.7.}
Recall that 
$\ol\CH^{\Bp} = \bigoplus_{\a \in \vD_{n,g}}\ol\CH^{\Bp}_{\a}$.
For each $\a \in \vD_{n,g}$, we denote by 
$T_{\a,j}$ the projection of $T_j$ onto $\ol\CH^{\Bp}_{\a}$. 
Also we denote by $\xi_{\a,i}$ the projection of $\xi_i$ onto
$\ol\CH_{\a}^{\Bp}$.
Hence we have $T_j = \sum_{\a}T_{\a,j}$ and 
$\xi_i = \sum_{\a}\xi_{\a,i}$.
It follows from the construction that under the isomorphism 
$\ol M_0^{\w} \simeq M^{\w^{[1]}}\otimes\cdots\otimes M^{\w^[g]}$, 
the action of $T_{\a, a_k + i}$ corresponds to the action of 
$T^{[k]}_i$ on $M^{\w^{[k]}}$.
\par
We note the following relation.
\begin{equation*}
\tag{6.7.1}
\xi_{\a,i}T^{[k]}_{\a,0} = T^{[k]}_{\a,0}\xi_{\a,i}
\end{equation*}
for any $i$ and any $k$.  In fact by (5.3.1), it is enough to show the
formula regarding $\xi_{\a, i}$ and $T^{[k]}_{\a,0}$ as operators 
on $\ol M_0^{\Om, \a}$.  Under the isomorphism 
$\ol M_0^{\w} \simeq M^{\w^{[1]}}\otimes\cdots\otimes M^{\w^{[g]}}$ 
for $\w \in \Om$ such that $\a_{\Bp}(\w) = \a$, 
$\xi_{\a,a_h + i}$ corresponds to the operator $\xi_i^{[h]}$ on 
$M^{\w^{[h]}}$, where $\xi_i^{[h]}$ is an element of $\CH_{n_h,r_h}$
defined similar to $\xi_i$ for $\ol\CH^{\Bp}$ (i.e., the special case
where $n = n_h, r = r_h, g = 1, \Bp = (r_h)$). 
But it is easy to see that in this case $\xi_i^{[h]}$ is a scalar 
multiplication on $M^{\w^{[h]}}$ by $Q^{\Bp}_h$.
Hence $\xi_{\a,i}$ is a scalars operator on $\ol M^{\w}$, and so 
commutes with $T^{[k]}_{\a,0}$. 
(6.7.1) follows from this.
\par
For each $\w \in \Om$ and $\a \in \vD_{n,g}$, let 
\begin{equation*}
F_{\w}(\xi_{\a}) = F_{c_1}(\xi_{\a,1})F_{c_2}(\xi_{\a,2})
      \cdots F_{c_n}(\xi_{\a,n})
\end{equation*}
with $c_i = b_{\w}(i)$.  We claim that 
\begin{equation*}
\tag{6.7.2}
F_{\w}(\xi_{\a}) = 0 \quad\text{ unless }  \a_{\Bp}(\w) =  \a.
\end{equation*}
\par
In fact, we have $\ol m_{\Om}F_{\w}(\xi_{\a}) \in \ol M^{\a'}$
by (6.6.2), where $\a' = \a_{\Bp}(\w)$.  
But since 
$F_{\w}(\xi_{\a}) \in \ol\CH^{\Bp}_{\a} = 
          \vf_{\Om,\a}\ol\CS^{\Bp}_{\a}\vf_{\Om,\a}$, 
we have $\ol m_{\Om}F_{\w}(\xi_{\a}) \in \ol M^{\a}$.
It follows that $\ol m_{\Om}F_{\w}(\xi_{\a}) = 0$ unless 
$\a_{\Bp}(\w) = \a$, and the claim
follows. 
\par
The following theorem gives a presentation of $\ol\CH^{\Bp}$. 
%%%%
\begin{thm}  %%%% Theorem 6.8
Assume that (6.3.1) holds.  
Recall that $Q^{\Bp}_k = Q_{p_k + r_k}$.
Then for each $\a \in \vD_{n,g}$, 
the algebra $\ol\CH^{\Bp}_{\a}$ is generated by 
\begin{equation*}
\xi_{\a,i}\  (1 \le i \le n),  \quad T_{\a,j}\  
(1 \le j \le n-1), 
\quad  T^{[k]}_{\a, 0}\ (1 \le k \le g)     
\end{equation*}
with  relations
\begin{align*}
\tag{A1}
&(T_{\a,i}-q)(T_{\a,i} + q\iv) = 0 \quad 
        (1 \le i \le n-1),  \\
\tag{A2}
&T_{\a,i}T_{\a,i+1}T_{\a,i} = T_{\a,i+1}T_{\a,i}T_{\a,i+1} 
        \quad  (1 \le i \le n-2), \\
\tag{A3}
&T_{\a, i}T_{\a,j} = T_{\a,j}T_{\a,i} \quad
              (1 \le i,j \le n-1, |i-j| \ge 2), \\
\tag{A4}
&(T^{[k]}_{\a,0} - Q_{p_k +1})\cdots 
        (T^{[k]}_{\a,0} - Q_{p_k + r_k}) = 0 \quad 
    (1 \le k \le g),  \\
\tag{A5}
&T^{[k]}_{\a,0}T_{\a,a_k+1}T^{[k]}_{\a,0}T_{\a,a_k+1}
   = T_{\a,a_k+1}T^{[k]}_{\a,0}T_{\a,a_k+1}T^{[k]}_{\a,0}
     \quad (1 \le k \le g),  \\
\tag{A6}
&T_{\a,0}^{[k]} = T_{\a, a_k}\cdots 
     T_{\a,a_{k-1}+ 2}T_{\a, a_{k-1}+1}
           T_{\a,0}^{[k-1]}T_{\a,a_{k-1}+1}T_{\a,a_{k-1}+2}\cdots
                           T_{\a, a_k}, \\
\tag{A7}
&T_{\a,0}^{[k]}T_{\a,j} = T_{\a,j}T_{\a,0}^{[k]} \quad
      (j \ne a_k, a_k +1), \\
\tag{A8}
&(\xi_{\a,i} - Q^{\Bp}_1)(\xi_{\a,i} - Q^{\Bp}_2)
               \cdots (\xi_{\a,i} - Q^{\Bp}_g) = 0
    \quad (1 \le i \le n), \\ 
\tag{A9}
&\xi_{\a,i}\xi_{\a,j} = \xi_{\a,j}\xi_{\a,i} \quad (1 \le i,j \le n), \\
\tag{A10}
&F_{\w}(\xi_{\a}) = 0 \quad\text{ if } \a_{\Bp}(\w) \ne \a,  \\
\tag{A11}
&T_{\a,j}\xi_{\a,j+1} = \xi_{\a,j}T_{\a,j} + \vD^{-2}
      \sum_{c_1 < c_2}(Q^{\Bp}_{c_2}- Q^{\Bp}_{c_1})(q-q\iv)
            F_{c_1}(\xi_{\a,j})F_{c_2}(\xi_{\a,j+1}), \\
\tag{A12}
&T_{\a,j}\xi_{\a,j} = \xi_{\a,j+1}T_{\a,j} - \vD^{-2}
      \sum_{c_1 < c_2}(Q^{\Bp}_{c_2}- Q^{\Bp}_{c_1})(q-q\iv)
            F_{c_1}(\xi_{\a,j})F_{c_2}(\xi_{\a,j+1}), \\
\tag{A13}
&T_{\a,j}\xi_{\a,k} = \xi_{\a,k}T_{\a,j} \quad (k \ne j, j+1), \\
\tag{A14}
&T_{\a,0}^{[k]}\xi_i = \xi_iT_{\a,0}^{[k]} \quad
      (1\le i \le n, 1 \le k \le g).
\end{align*}  
\end{thm}
\begin{proof}
One sees that these elements generate $\ol\CH^{\Bp}_{\a}$ 
by Lemma 6.6. We show that these generators satisfy the relations 
(A1) $\sim$ (A14).  (A1) $\sim$ (A3) follows from the relations for 
$\CH_n$.  (A8) follows from (6.1.2).  (A9) is also clear from 
6.1. (A10) follows from (6.7.2).  
(A11) $\sim$ (A13) follows from Lemma 6.4.  (A14) is given 
in (6.7.1).  We show the remaining relations (A4) $\sim$ (A7).
We may prove the formulas by regarding 
 $T_{\a,0}^{[k]}$ and $T_{\a,j}$ as operators on $\ol M_0^{\w}$
for $\w \in \Om$ such that $\a_{\Bp}(\w) = \a$ 
by (5.3.1).  Since $T^{[k]}_{\a,0}$ corresponds to the action of  
$T_0^{[k]} \in \CH_{n_k,r_k}$ on $M^{\w^{[k]}}$, 
and $T_{\a,a_k+i}$ corresponds to the action of $T^{[k]}_i$,
(A4), (A5) and (A7) follows from the relations for $\CH_{n_k,r_k}$.
While (A6) follows from the property that $T^{[k]}_{\a,0}$ is the 
restriction of $L_{a_k+1}$ on $\ol M^{\Om,\a}$.  
Thus those generators satisfy the relations (A1) $\sim$ (A14). 
\par
Next we show that (A1) $\sim$ (A14) gives a fundamental relation 
for $\ol\CH^{\Bp}_{\a}$.  Let $\wh\CH_{\a}$ be the algebra with 
generators $\wh\xi_{\a,i}, \wh T_{\a,j}$ and $\wh T_{\a,0}^{[k]}$,
and relations as in the theorem.  (We denote by $\wh X$ the generator
in $\wh\CH_{\a}$ corresponding to the generator $X$ in 
$\ol\CH^{\Bp}_{\a}$.)  
Let $\wh\CH^0_{\a}$ be the subalgebra of $\wh\CH_{\a}$ 
generated by $\wh T^{[k]}_{\a,i}$ for 
$k = 1, \dots, g, i = 0, \dots, n_k-1$.
Recall that $\wh T^{[k]}_{\a,i} = \wh T_{\a, a_k+i}$ for 
$1 \le i \le n_k-1$.  Then by the relations in the theorem, 
$\wh \CH_{\a}^0$ is isomorphic to the quotient of the algebra
$\CH_{n_1,r_1}\otimes\cdots\otimes \CH_{n_g,r_g}$.
Also we note that the subalgebra $\wh \CH_n$ of $\wh\CH_{\a}$ generated by 
$\wh T_{\a,j}$ is the quotient of $\CH_n$. We denote by $\wh T_{\a,w}$
the image of $T_w \in \CH_n$ to $\wh\CH_n$ for $w \in \FS_n$.
Let $\FS_{\a}$ be the Young subgroup of $\FS_n$ corresponding 
to the composition $\a$ of $n$.
Let $\wh\Xi_{\a}$ be the subalgebra of $\wh \CH_{\a}$ generated by 
$\wh\xi_{\a,1}, \dots, \wh\xi_{\a,n}$.   
For each $\w \in \Om^{\a}$, we define 
$F_{\w}(\wh\xi_{\a}) \in \wh\Xi_{\a}$ 
in a similar way as $F_{\w}(\xi_{\a})$, but replacing $\xi_{\a,i}$ by 
$\wh\xi_{\a,i}$.
We show that
\par\medskip\noindent
(6.8.1)\  Any element of $\wh\CH_{\a}$ can be written as a linear
combination of elements in  
\begin{equation*}
\ZC = \{ F_{\w}(\wh\xi_{\a})\wh\CH_{\a}^0\wh T_{\a, w} 
      \mid \w \in \Om^{\a}, w \in \FS_{\a}\backslash \FS_n\}.
\end{equation*}
\par
In fact, let $\wh \CH_{\a}\nat$ be the subalgebra of 
$\wh \CH_{\a}$ generated by $\wh T_{\a,j}$ and $\wh T^{[k]}_{\a,0}$.
Then by the commuting relations in the theorem, $\wh\CH_{\a}$ can be 
written as 
\begin{equation*}
\wh \CH_{\a} = \sum_{c_1, \dots, c_n}
    \wh \xi^{c_1}_{\a,1}\wh\xi_{\a,2}^{c_2}\cdots 
           \wh\xi_{\a,n}^{c_n}\wh\CH_{\a}\nat,
\end{equation*}
where $c_i$ are integers such that $0 \le c_i \le g-1$. 
It is easy to see that any element in $\wh\Xi_{\a}$ can be written 
as a linear combination of $F_{\w}(\wh\xi_{\a})$ 
for various $\w \in \Om$.  Thus by (A10), any element in $\wh\CH_{\a}$
is written as a linear combination of 
$F_{\w}(\wh\xi_{\a})\wh\CH\nat_{\a}$ with $\w \in \Om^{\a}$.
We now concentrate on $\wh\CH\nat_{\a}$.  
Define $\wh L_i^{[k]} \in \wh\CH^0_{\a}$
by 
\begin{equation*}
\wh L_i^{[k]} = 
    \wh T_{\a, a_k + i-1}\cdots \wh T_{\a,a_k + 2}\wh T_{\a,a_k + 1}
      \wh T^{[k]}_{\a,0}
         \wh T_{\a,a_k+1}\wh T_{\a,a_k+2}\cdots \wh T_{\a, a_k+i-1}
\end{equation*}
for $i = 1, \dots, n_k$.
Then $L_i^{[k]}$ commutes with $\wh T_{\a, j}$ for 
$j \ne a_k+ i-1, a_k + i$ 
and we have 
\begin{equation*}
\wh T_{\a, a_k+i}\wh L_{i}^{[k]}\wh T_{\a,a_k+i}
      = \begin{cases}
             \wh L_{i+1}^{[k]}  &\quad\text{ if } i \ne n_k, \\
             \wh L_1^{[k+1]}  &\quad\text{ if }  i = n_k.
        \end{cases}
\end{equation*}
by (A6). (Note that in the latter case, 
$\wh T_{\a, a_k+n_k} \notin \wh\CH^0_{\a}$.)
It follows that any element in $\wh\CH_{\a}\nat$ can be written 
as a linear combination of the elements in $\wh L \wh\CH_n$, where 
$\wh L$ is the subalgebra of $\wh\CH^0_{\a}$ generated by 
$\wh L_i^{[k]}$. 
Let $\CH_{n, \a}$ be the subalgebra of $\CH_n$ corresponding to the
Young subgroup $\FS_{\a}$ of $\FS_n$, and let $\wh\CH_{n,\a}$ be 
the corresponding subalgebra of $\wh\CH_n$.
Since $\wh\CH_n$ is the quotient of $\CH_n$, it is written as 
a sum of $\wh\CH_{n,\a}\wh T_{\a,w}$ with $w \in \FS_{\a}\backslash \FS_n$.
Since $\wh\CH_{n,\a} \subset \wh\CH_{\a}^0$, one sees that 
$\wh\CH\nat_{\a} = 
 \sum_{w \in \FS_{\a}\backslash \FS_n}\wh\CH^0_{\a}\wh T_w$.
Hence (6.8.1) holds.
\par
Since $\ol\CH^{\Bp}_{\a}$ satisfies the same relations, 
we have a surjective homomorphism 
$\p: \wh\CH_{\a} \to \ol\CH^{\Bp}_{\a}$.
In order to show that $\p$ is injective, it is enough to 
see that  
the set of elements in (6.8.1) gives an $R$-free basis of 
$\wh\CH_{\a}$ and that the image under $\p$ of this basis
gives a basis of $\ol\CH^{\Bp}_{\a}$. 
We denote by $\ZC'$ the image of $\ZC$ under $\p$.  By a similar
argument as above, we see that $\ZC'$ spans $\ol\CH^{\Bp}_{\a}$ 
as an $R$-module.  We show that $\ZC'$ gives an $R$-free basis of
$\ol\CH^{\Bp}$.  For this, it is enough to see that the elements 
in $\ZC'$ are linearly independent over $R$, or equivalently, 
they are linearly independent over $K$, where $K$ is the quotient 
field of $R$.
It is easy to see that the cardinality of the set $\ZC'$
is equal to 
\begin{equation*}
|\Om^{\a}| \times \dim \CH_{\a}\times n_{\a}
= n_{\a}^2 \times \dim \CH_{\a} = \dim \ol \CH^{\Bp}_{\a}
\end{equation*}
by Proposition 5.3.
Hence the elements of $\ZC'$ are linearly independent, and 
$\ZC'$ gives an $R$-free basis of $\ol\CH^{\Bp}_{\a}$.
This shows that the elements in $\ZC$ are also linearly independent, 
and so $\ZC$ is an $R$-free basis of $\wh \CH_{\a}$. 
Therefore $\p$ is an isomorphism, and the theorem is proved.
\end{proof}
\remark{6.9.}
In the case where $\Bp = (r)$, $\ol\CH^{\Bp}_{\a} = \ol\CH^{\Bp}$ 
coincides with  
$\CH$, and the fundamental relation (A1) $\sim (A14)$ is reduced 
to the fundamental relation for $\CH$.
On the other hand, in the case where $\Bp = (1^r)$,  $\CH_{\a}$ is 
a subalgebra of $\CH_n$ for each $\a \in \vD_{n,r}$. Then 
$T_{\a,0}^{[k]}$ turns out to be scalar operators, and the relations 
(A4) $\sim$ (A7), (A14) can be ignored. The remaining relations give
the fundamental relation for $\ol\CH^{\Bp}_{\a}$.  Note that a similar
argument as in the proof shows that the relations (A1) $\sim$ (A3), 
(A8), (A9), (A11) $\sim$ (A13) gives a fundamental relation for 
$\ol\CH^{\Bp}$, which is nothing but the fundamental relation for 
the modified Ariki-Koike algebra given in [SawS]. 
%%%
%%%
\section{Schur-Weyl duality}
\para{7.1.}
It is known by [M, \S 5] that the Schur-Weyl duality i.e., 
the double centralizer property holds between 
$\CH$ and $\CS = \End_{\CH}M$.
A similar duality also holds by [SawS, Theorem 8.2] for 
the modified Ariki-Koike algebra $\ol\CH$ on the action of the tensor space 
$V^{\otimes n}$.  In our setting, $\ol\CH$ coincides with 
$\ol\CH^{\Bp}$ with $\Bp = (1^r)$, 
and $V^{\otimes n} \simeq \ol M$ as $\ol\CH^{\Bp}$-modules.
In what follows we shall give a generalization of this property 
for the arbitrary $\Bp$, i.e., we show the Schur-Weyl duality 
between $\ol\CS^{\Bp}$ and $\ol\CH^{\Bp}$ acting on $\ol M$.  
Although the proof is carried out 
for the action on $\ol M$, we formulate the theorem for
$\ol\CH^{\Bp}$-module $M_{\Bp} = \bigoplus M_{\Bp}^{\mu}$ 
which is isomorphic to $\ol M$, where $M_{\Bp}^{\mu}$ is a right 
ideal of $\ol\CH^{\Bp}$, so that it fits to the situation above.
%%%
\para{7.2.}
In order to give an expression of $\ol M^{\mu}$ as a right ideal
of $\ol\CH^{\Bp}$, we describe the cellular basis $m_{\Fs\Ft}^{\Bp}$ of 
$\ol\CH^{\Bp}$ more explicitly.
For each $\a = (n_1, \dots, n_g)\in \vD_{n,g}$ 
we define $F_{\a} \in \ol\CH^{\Bp}$ by 
\begin{equation*}
\tag{7.2.1}
F_{\a} = \vD^{-n}F_{c_1}(\xi_1)\cdots F_{c_n}(\xi_n),
\end{equation*}
where 
\begin{equation*}
(c_1, \dots, c_n) = (\underbrace{1,\dots,1}_{n_1\text{-times}},
                     \underbrace{2,\dots,2}_{n_2\text{-times}},
                     \dots, 
                     \underbrace{g,\dots,g}_{n_g\text{-times}}).
\end{equation*}
If we define $\w = \w_{\a}$ as the unique element in $\Om^{\a}$ such that
$d(\w) = 1$, we see that $F_{\a} = \vD^{-n}F_{\w}(\xi)$ in the
notation of (6.5.1).  It follows from (6.6.2) that 
\begin{equation*}
\tag{7.2.2}
\ol m_{\Om}F_{\a} = \ol m_{\w}. 
\end{equation*}
Take $\la \in \vL^+$ such that $\a_{\Bp}(\la) = \a$.  
Then $\Ft^{\la} \in \Std(\la)_{\w}$.  Let 
$S_{\w} = \w(\Ft^{\la}) \in \CT_0^{\Bp}(\la, \w)$.
Then 
$\ol m_{S_{\w}\Ft^{\la}} = 
\ol m_{\Ft^{\la}\Ft^{\la}} = \ol m_{\la} \in \ol M^{\w}$.
Since $\ol M^{\w} = \ol m_{\w}\ol\CH^{\Bp}$, there exists 
$y_{\la} \in \ol\CH^{\Bp}$ such that 
$\ol m_{\la} = \ol m_{\w}y_{\la}$.
One can choose $y_{\la}$ in the following way. 
Let $\ol\CH^0_{\a}$ be the subalgebra of $\ol\CH^{\Bp}$
consisting of scalar matrices with entries in 
$\CH_{\a} = \CH_{n_1,r_1}\otimes\cdots\otimes\CH_{n_g,r_g}$ under the 
isomorphism in Proposition 5.3.  Thus $\ol\CH_{\a}^0 \simeq \CH_{\a}$.
We have $\ol m_{\w}, \ol m_{\la} \in \ol M_0^{\w}$, 
and under the isomorphism 
$\ol M_0^{\w} \simeq M^{\w^{[1]}}\otimes\cdots\otimes M^{\w^{[g]}}
= \CH_{\a}$,
$\ol m_{\w}$ corresponds to $1\otimes\cdots\otimes 1$, and 
$\ol m_{\la}$ corresponds to 
$m_{\la^{[1]}}\otimes\cdots\otimes m_{\la^{[g]}}$ in 
$\CH_{\a}$.
Then we choose $y_{\la} \in \ol \CH_{\a}^{0}$ as the scalar
matrix consisting of $m_{\la^{[1]}}\otimes\cdots\otimes m_{\la^{[g]}}$ in 
$\CH_{\a}$ under the isomorphism in Proposition 5.3.  
\par
Note that $F_{\a}$ commutes with any element in $\ol\CH_{\a}^0$. 
In fact by (7.2.1) $F_{\a} \in R[\xi_1, \dots, \xi_n]^{\FS_{\a}}$ 
with $\FS_{\a} = \FS_{n_1}\times\cdots\times \FS_{n_g}$, and a similar
argument as in [SawS, Lemma 2.8] can be applied. 
In particular, $y_{\la}$ commutes with $F_{\a}$.
Let $* : \ol\CH^{\Bp} \to \ol\CH^{\Bp}$ be the anti-automorphism.
Since $\xi_i$  are fixed by $*$, $F_{\a}$ is fixed by $*$.
Also $y_{\la}$ is fixed by $*$ since the corresponding elements 
in $\CH_{n_k,r_k}$ are fixed by $*$.
We have the following 
lemma.
%%%
\begin{lem}  %%% Lemma 7.3
For each $\Ft, \Fs \in \Std(\la)$, we have 
\begin{equation*}
m_{\Fs\Ft}^{\Bp} = T_{d(\Fs)}^*F_{\a}y_{\la}T_{d(\Ft)}.
\end{equation*}
\end{lem}
\begin{proof}
By the construction in 7.2, we see that 
$\ol m_{\Om}F_{\a}y_{\la} = \ol m_{\Ft^{\la}\Ft^{\la}}$ for 
$\la$ such that $\a_{\Bp}(\la) = \a$. 
Thus $\ol m_{\Om}F_{\a}y_{\la}T_{d(\Ft)} = \ol m_{\Ft^{\la}\Ft}$
for any $\Ft \in \Std(\la)$.
If $T \in \CT_0^{\Bp}(\la, \w')$ corresponds to 
$\Ft \in \Std(\la)$ under (5.1.2), and 
$S_{\w} \in \CT_0^{\Bp}(\la, \w)$ with $\w = \w_{\a}$, 
then we have $\vf_{S_{\w}T}(\ol m_{\Om}) = \ol m_{\Ft^{\la}\Ft}$.
It follows that $\vf_{S_{\w}T} = F_{\a}y_{\la}T_{d(\Ft)}$.
This shows that $\vf_{SS_{\w}} = T_{d(\Fs)}^*F_{\a}y_{\la}$.
Take $T \in \CT_0^{\Bp}(\la, \w')$ corresponding to 
$\Ft \in \Std(\la)$. Since 
$\vf_{SS_{\w}}(\ol m_{\Om}) = \ol m_{\Fs\Ft^{\la}}$, we have 
\begin{equation*}
\ol m_{\Om}\cdot T_{d(\Fs)}^*F_{\a}y_{\la}T_{d(\Ft)}
    = \ol m_{\Fs\Ft^{\la}}\cdot T_{d(\Ft)} = \ol m_{\Fs\Ft} 
    = \vf_{ST}(\ol m_{\Om}).
\end{equation*}
Thus we have 
$m_{\Fs\Ft}^{\Bp} = \vf_{ST} = T_{d(\Fs)}^*F_{\a}y_{\la}T_{d(\Ft)}$. 
\end{proof}
\para{7.4.}
For each $\mu \in \vL$ such that $\a_{\Bp}(\mu) = \a$, 
we define $y_{\mu} \in \ol\CH^{\Bp}_{\a}$ similarly as before, 
by extending the definition of $y_{\la}$ for $\la \in \vL^+$.
We define a right ideal $M_{\Bp}^{\mu}$ of $\ol\CH^{\Bp}$
by $M_{\Bp}^{\mu} = F_{\a}y_{\mu}\ol\CH^{\Bp}$ and put
$M_{\Bp} = \bigoplus_{\mu \in \vL}M_{\Bp}^{\mu}$.
By Lemma 4.10, we have $\ol m_{\Om}F_{\a}y_{\mu} = \ol m_{\mu}$
and so $\ol m_{\Om}F_{\a}y_{\mu}\ol\CH^{\Bp} = 
   \ol m_{\mu}\ol\CH^{\Bp} = \ol M^{\mu}$.
This shows that there exists an isomorphism
$\f : M_{\Bp}^{\mu} \to \ol M^{\mu}$ of $\ol\CH^{\Bp}$-modules by 
$F_{\a}y_{\mu}h \mapsto \ol m_{\Om}F_{\a}y_{\mu}h = \ol m_{\mu}h$.
\par
Recall that $\{ \ol m_{S\Ft} \mid S \in \CT_0^{\Bp}(\la,\mu), 
   \Ft \in \Std(\la) \text{ for } \la \in \vL^+\}$ gives a basis 
of $\ol M^{\mu}$.
In connection with this, we define,
for each $S \in \CT_0^{\Bp}(\la,\mu), \Ft \in \Std(\la)$ with 
$\la \in \vL^+$, 
\begin{equation*}
m_{S\Ft}^{\Bp} = \sum_{\substack{\Fs \in \Std(\la)\\
                         \mu(\Fs) = S}}
                         q^{l(d(\Fs)) + l(d(\Ft))}m_{\Fs\Ft}^{\Bp}.
\end{equation*}
The following lemma holds.
%%%
\begin{lem} %%% Lemma 7.5
The set $\{ m_{S\Ft}^{\Bp}\}$ gives rise to a basis of $M^{\mu}_{\Bp}$, 
and we have $\f(m^{\Bp}_{S\Ft}) = \ol m_{S\Ft}$ for each basis element.  
\end{lem}
\begin{proof}
By the proof of Lemma 7.3, we know that 
$\ol m_{\Om}m_{\Fs\Ft}^{\Bp} = \ol m_{\Fs\Ft}$ for any 
$\Ft, \Fs \in \Std(\la)$.  It follows that 
$\ol m_{\Om}m_{S\Ft}^{\Bp} = \ol m_{S\Ft} \in \ol M^{\mu}$ 
for any $S \in \CT_0^{\Bp}(\la,\mu)$
and $\Ft \in \Std(\la)$.  In particular, we see that 
$m_{S\Ft}^{\Bp} \in M_{\Bp}^{\mu}$, and the lemma follows. 
\end{proof}
The following result gives the Schur-Weyl duality, i.e., the 
double centralizer property  between $\ol\CH^{\Bp}$
and $\ol\CS^{\Bp}$.
\begin{thm}  %%%% Theorem 7.6
Under the assumptions (5.1.1) and (6.3.1), there exist isomorphisms 
of $R$-algebras
\begin{equation*}
\ol\CS^{\Bp} \simeq \End_{\ol\CH^{\Bp}}M_{\Bp}, \qquad
\ol\CH^{\Bp} \simeq \End^0_{\ol\CS^{\Bp}}M_{\Bp}. 
\end{equation*}
\end{thm}
\begin{proof}
We argue on $\ol M$ instead of $M_{\Bp}$.  The second isomorphism
is already shown in (5.4.1).  So we prove the first isomorphism.
Let $\mu, \nu \in \vL$ be such that 
$\a_{\Bp}(\mu) = \a_{\Bp}(\nu) = \a$, 
and take $\vf \in \Hom_{\ol\CH^{\Bp}}(\ol M^{\nu}, \ol M^{\mu})$.
Since $\ol M^{\nu} = \ol m_{\nu}\ol\CH^{\Bp}$, the map 
$\vf$ is determined by $\vf(\ol m_{\nu})$.  We show that
\begin{equation*}
\tag{7.6.1}
\vf(\ol m_{\nu}) \in \ol M^{\mu}_0.
\end{equation*}
\par
In fact, since $\ol m_{S\Ft}$  
($S \in \CT_0^{\Bp}(\la, \mu), \Ft \in \Std(\la)$) gives a basis of 
$\ol M^{\mu}$, one can write
\begin{equation*}
\vf(\ol m_{\nu}) = \sum_{S,\Ft}c_{S\Ft}\ol m_{S\Ft}
\end{equation*}
with $c_{S\Ft} \in R$.  By (6.1.4), we have
\begin{equation*}
\tag{7.6.2}
\vf(\ol m_{\nu}\xi_i) = Q^{\Bp}_{b(i)}\sum_{S, \Ft}c_{S\Ft}\ol m_{S\Ft}
\end{equation*}
for $i = 1, \dots, n$, where $b(i) = k$ if $a_k+1 \le i \le a_k+n_k$.
On the other hand, by (6.1.3), we have
\begin{equation*}
\tag{7.6.3}
\vf(\ol m_{\nu}\xi_i) = \vf(\ol m_{\nu})\xi_i = 
    \sum_{S, \Ft}c_{S\Ft}Q^{\Bp}_{\Ft(i)}\ol m_{S\Ft},
\end{equation*}
where $\Ft(i) = k$ if the letter $i$ is contained in 
$\Ft^{[k]}$ (see the remark after (6.1.3)). 
Comparing (7.6.2) and (7.6.3), we see that $\Ft \in \Std(\la)_0$.
Since $\ol M^{\mu}_0$ is spanned by $\ol m_{S\Ft}$ with
$\Ft \in \Std(\la)_0$, we obtain (7.6.1). 
\par
Let $\ol\CH_{\a}^0$ be the subalgebra of $\ol\CH^{\Bp}$ as before.
Since $\ol M_0^{\nu} = \ol m_{\nu}\ol\CH_{\a}^0$, and similarly for
$\ol M^{\mu}_0$, it follows from (7.6.1) that any 
$\vf \in \Hom_{\ol\CH^{\Bp}}(\ol M^{\nu}, \ol M^{\mu})$  
has the property that $\vf(\ol M^{\nu}_0) \subset \ol M^{\mu}_0$.
Thus we have a natural $R$-linear map 
\begin{equation*}
\th'': \Hom_{\ol\CH^{\Bp}}(\ol M^{\nu}, \ol M^{\mu}) 
      \to \Hom_{\ol\CH^0_{\a}}(\ol M^{\nu}_0, \ol M^{\mu}_0),
\end{equation*}
which is clearly injective.
Let $\ol H_{\mu\nu}$ be the $\mu\nu$-part of $\ol S^{\Bp}$ as in 
(4.9.2).  Since the action of $\wt H_{\a}$ on $\ol M^{\nu}_0$
coincides with the action of $\ol\CH^0_{\a}$, we see that 
there exists an $R$-linear isomorphism
\begin{equation*}
\vT: \ol H_{\mu\nu} \to \
    \Hom_{\ol \CH_{\a}^0}(\ol M^{\nu}_0, \ol M^{\mu}_0)
\end{equation*}
by Lemma 4.12.  On the other hand by (5.4.2), we know that there exists 
an injective map 
$\th': \ol H_{\mu\nu} \to \Hom_{\ol\CH^{\Bp}}(\ol M^{\nu}, \ol M^{\mu})$.
It is clear that the composite of $\th'$ and $\th''$ coincides with 
$\vT$.  Hence $\th'$ is an isomorphism.  This shows that
$\ol\CS^{\Bp} \simeq \End_{\ol\CH^{\Bp}}\ol M$, and the theorem follows.
\end{proof}
\remark{7.7.} 
The assumption (6.3.1) is used to give an expression of $\ol M$
as an ideal of $\ol\CH^{\Bp}$.  But the Schur-Weyl duality holds
for $\ol M$ without referring the ideal $M_{\Bp}$.  In that case, 
(6.3.1) can be replaced by a weaker assumption 
``the parameters $Q_1^{\Bp}, \dots, Q_g^{\Bp}$ are all distinct''. 
\par\medskip
By making use of Theorem 7.6, we obtain the following additional
information on the space 
$\ol H_{\mu\nu} = \Hom_{\ol\CH^{\Bp}}(M_{\Bp}^{\nu}, M_{\Bp}^{\mu})$.
Put $m_{\nu}^{\Bp} = F_{\a}y_{\nu}$ so that 
$M^{\nu}_{\Bp} = m_{\nu}^{\Bp}\ol \CH^{\Bp}$.
%%%
\begin{prop}  %%% Prop. 7.8.
Let $\mu, \nu \in \vL$ such that 
$\a_{\Bp}(\mu) = \a_{\Bp}(\nu) = \a$.
\begin{enumerate}
\item
The map $\vf \mapsto \vf(m_{\nu}^{\Bp})$ gives an isomorphism
of $R$-modules, 
\begin{equation*}
\Hom_{\ol\CH^{\Bp}}(M_{\Bp}^{\nu}, M_{\Bp}^{\mu}) \to 
M_{\Bp}^{\nu*} \cap M_{\Bp}^{\mu},
\end{equation*}
where $M_{\Bp}^{\nu*} = \ol\CH^{\Bp}m_{\nu}^{\Bp}$ is the 
image of $M_{\Bp}^{\nu}$ under the operation $*$.
\item
We have
\begin{equation*}
M_{\Bp}^{\nu*} \cap M_{\Bp}^{\mu} = 
   F_{\a}(\ol\CH^0_{\a}y_{\nu} \cap y_{\mu}\ol\CH^0_{\a}). 
\end{equation*}
\end{enumerate}
\end{prop}
\begin{proof}
For each $m \in M_{\Bp}^{\nu*} \cap M_{\Bp}^{\mu}$, 
the map $m_{\nu}^{\Bp}h \mapsto mh$ ($h \in \ol\CH^{\Bp}$)
gives a well-defined map 
$\vf_m \in \Hom_{\ol\CH^{\Bp}}(M_{\Bp}^{\nu}, M_{\Bp}^{\mu})$, 
and the map $m \mapsto \vf_m$ gives an $R$-linear 
map $M_{\Bp}^{\nu*}\cap M_{\Bp}^{\mu} \to 
      \Hom_{\ol\CH^{\Bp}}(M_{\Bp}^{\nu}, M_{\Bp}^{\mu})$, 
which is clearly injective.
\par
On the other hand, we have
\begin{align*}
\Hom_{\ol\CH^{\Bp}}(M_{\Bp}^{\nu}, M_{\Bp}^{\mu})
  &= \Hom_{\ol\CH^{\Bp}_{\a}}(M_{\Bp}^{\nu}, M_{\Bp}^{\mu})  \\
  &\simeq \Hom_{\ol\CH_{\a}^0}
   (m^{\Bp}_{\nu}\ol\CH_{\a}^0, m^{\Bp}_{\mu}\CH_{\a}^0)  \\
  &\simeq \Hom_{\CH_{\a}}
   (M^{\nu^{[1]}}\otimes\cdots\otimes M^{\nu^{[g]}},
    M^{\mu^{[1]}}\otimes\cdots\otimes M^{\mu^{[g]}}),
\end{align*}
where $\CH_{\a} = \CH_{n_1,r_1}\otimes\cdots\otimes \CH_{n_g,r_g}$.
(Note that $m^{\Bp}_{\mu}\ol\CH^0_{\a}$ corresponds to 
$\ol M_0^{\mu}$ under the isomorphism 
   $M_{\Bp}^{\mu} \simeq \ol M^{\mu}$.)
It is known, by [DJM] that the last set is isomorphic to
$\ol\CH_{\a}^0y_{\nu}\cap y_{\mu}\ol\CH_{\a}^0$ as $R$-modules.  Hence
\begin{equation*}
\Hom_{\ol\CH_{\a}^0}
   (m^{\Bp}_{\nu}\ol\CH_{\a}^0, m^{\Bp}_{\mu}\CH_{\a}^0) \simeq
(m^{\Bp}_{\nu}\ol\CH_{\a}^0)^* \cap 
                   m^{\Bp}_{\mu}\ol\CH_{\a}^0
\end{equation*}
via the map $\vf \mapsto \vf(m^{\Bp}_{\nu})$.
But since $m^{\Bp}_{\nu}\ol\CH^0_{\a} = F_{\a}y_{\nu}\ol\CH^0_{\a}$
and $F_{\a}$ commutes with $y_{\nu}$ and $\ol\CH^0_{\a}$, we see that
$(m^{\Bp}_{\nu}\ol\CH^0_{\a})^* = F_a\ol\CH^0_{\a}y_{\nu}$.  This shows 
that
\begin{equation*}
\Hom_{\ol\CH^{\Bp}}(M^{\Bp}_{\nu}, M^{\Bp}_{\mu})
\simeq F_{\a}(\ol\CH^0_{\a}y_{\nu} \cap y_{\mu}\ol\CH^0_{\a}) 
\subseteq M_{\Bp}^{\nu*} \cap M_{\Bp}^{\mu}, 
\end{equation*}
where the first isomorphism is given by the map 
$\vf \mapsto \vf(m^{\Bp}_{\nu})$.
Both statements of the proposition follow from this, by 
combined with the remark in the first paragraph. 
\end{proof}

%%%
%%%
\section{Comparison of $\ol\CH^{\Bp}$ for various $\Bp$}
\para{8.1.}
We shall consider the relationship among $\CS^{\Bp}$ and 
$\ol\CH^{\Bp}$ for various types $\Bp$.  First consider the case
of $\CS^{\Bp}$.
Let $\Bp = (r_1, \dots, r_g)$ and $\Bp' = (r'_1, \dots, r'_{g'})$
be two compositions of $r$.  We define $(p_1', \dots, p'_{g'})$ 
for $\Bp'$ similar to $\Bp$.   
We write $\Bp' \preceq \Bp$ if $\Bp'$ is obtained as a refinement of 
$\Bp$, namely if $p_j$ coincides with some $p'_{k_j}$ for each $j$.
In particular, we have $(1^r) \preceq \Bp \preceq (r)$ for any $\Bp$.
Assume that $\Bp' \preceq \Bp$.  Then we see that 
$\Ba_{\Bp}(\la) \ge \Ba_{\Bp}(\mu)$ if 
$\Ba_{\Bp'}(\la) \ge \Ba_{\Bp'}(\mu)$ for $\la, \mu \in \vL$.
Moreover, $\a_{\Bp}(\la) = \a_{\Bp}(\mu)$ if 
$\a_{\Bp'}(\la) = \a_{\Bp'}(\mu)$.
This implies that 
\begin{equation*}
\tag{8.1.1}
\CS^{\Bp'} \subseteq \CS^{\Bp} \quad\text{ if } \Bp' \preceq \Bp.
\end{equation*}
Concerning the modified Ariki-Koike algebras $\ol\CH^{\Bp}$, we have
the following.
%%%
\begin{prop}  %%% Proposition 8.2.
There exists an algebra homomorphism 
$\r_{\Bp'\Bp}: \ol\CH^{\Bp} \to \ol\CH^{\Bp'}$
for any pair $\Bp,\Bp'$ such that $\Bp' \preceq \Bp$ satisfying the 
following property;  for $\Bp'' \preceq \Bp' \preceq \Bp$, we have
$\r_{\Bp''\Bp} = \r_{\Bp''\Bp'}\circ \r_{\Bp'\Bp}$.
\end{prop}
\begin{proof}
Let $\ol M = \bigoplus_{\mu}\ol M^{\mu}$ be the $\CH$-module defined by 
$\wh N^{\Ba_{\Bp}(\mu)}$ as before.  We denote $\ol M$ by 
$\ol M_{\Bp}$ to indicate its dependence on $\Bp$. 
Assume that $\Bp' \preceq \Bp$.
Then we have a natural surjection 
$\ol M_{\Bp} \to \ol M_{\Bp'}$ of $\CH$-modules.
If we regard $\ol M_{\Bp}$ as a left $\CS^{\Bp}$-module, and 
$\ol M_{\Bp'}$ as a left $\CS^{\Bp'}$-module, then the map
$\ol M_{\Bp} \to \ol M_{\Bp'}$ is compatible with the actions 
of $\CS^{\Bp}$ and $\CS^{\Bp'}$  via the inclusion 
$\CS^{\Bp'} \hra \CS^{\Bp}$.  
Let $M_{\Bp} = \bigoplus_{\mu} M_{\Bp}^{\mu}$ be as in 7.4.  
By Theorem 7.6, $\ol\CH^{\Bp}$ is realized as 
$\ol\CH^{\Bp} = \End^0_{\ol\CS^{\Bp}}M_{\Bp}$.
By using the property of the cellular structure of $\ol\CH^{\Bp}$
described in the beginning of Section 7, together with Lemma 7.5, 
the above property of $\ol M_{\Bp}$  
can be made more precise for $M_{\Bp}$ as follows (which is a generalization
of the argument in 4.1).
Let $\wh N_{\Bp}^{\Ba_{\Bp'}(\mu)}$ be the $R$-submodule of $\ol\CH^{\Bp}$
spanned by $m^{\Bp}_{\Fs\Ft}$ such that $\Fs,\Ft \in \Std(\la)$ with
$\Ba_{\Bp'}(\la) > \Ba_{\Bp'}(\mu)$.  Then 
$\wh N_{\Bp}^{\Ba_{\Bp'}(\mu)}$ is a two-sided ideal of $\ol\CH^{\Bp}$.
Put $\wh M_{\Bp}^{\mu} = M_{\Bp}^{\mu} \cap N_{\Bp}^{\Ba_{\Bp'}(\mu)}$.
Then $\wh M_{\Bp}^{\mu}$ is an $\ol\CH^{\Bp}$-submodule of 
$M_{\Bp}^{\mu}$ with the basis 
$\{ m^{\Bp}_{S\Ft} \mid S \in \CT^{\Bp}_0(\la,\mu), \Ft \in \Std(\la),
                \Ba_{\Bp'}(\la) > \Ba_{\Bp'}(\mu)\}$, and we have 
an isomorphism of $R$-modules
\begin{equation*}
M_{\Bp}^{\mu}/\wh M_{\Bp}^{\mu} \simeq M_{\Bp'}^{\mu}.
\end{equation*}
A similar argument as in 4.9 shows that any $\vf \in \ol H_{\mu\nu}$
maps $\wh M^{\nu}$ to $\wh M^{\mu}$, and so 
$\vf \in \ol\CS^{\Bp} = \End_{\ol\CH^{\Bp}}M_{\Bp}$ induces an action 
on $M_{\Bp}/\wh M_{\Bp}$, where 
$\wh M_{\Bp} = \bigoplus_{\mu}\wh M_{\Bp}^{\mu}$. 
This gives an isomorphism $M_{\Bp}/\wh M_{\Bp} \simeq M_{\Bp'}$ as
$\ol\CS^{\Bp'}$-modules.
\par
Now the action of $\ol\CH^{\Bp}$ on $M_{\Bp}$ 
induces an action on $M_{\Bp}/\wh M_{\Bp}$, which is compatible
with the action of $\ol\CS^{\Bp}$.  Hence this induces an action of
$\ol\CH^{\Bp}$ on $M_{\Bp'}$ compatible with the action of 
$\ol\CS^{\Bp'}$.  Thus we have an $R$-algebra homomorphism
\begin{equation*}
\r_{\Bp'\Bp}: 
    \ol\CH^{\Bp} \to \End^0_{\ol\CS^{\Bp'}}M_{\Bp'} \simeq \ol\CH^{\Bp'}.
\end{equation*}
It is clear that this map $\r_{\Bp'\Bp}$ satisfies the required 
property.
\end{proof}
\para{8.3.}
In the case where $\Bp = (r)$, we have $\ol\CH^{\Bp} \simeq \CH$, and 
in the case where $\Bp' = (1^r)$, we have 
$\ol\CH^{\Bp'} \simeq \ol\CH$, the modified Ariki-Koike algebra 
introduced in [SawS].
We have $\Bp' \preceq \Bp$, and the map 
$\r_{\Bp'\Bp} : \CH \to \ol\CH$ coincides with the map $\r_0$ given 
in [SawS, Lemma 1.5].  We consider the following separation condition 
on parameters of $\CH$, which was first introduced in [A].
\par\medskip\noindent
(8.3.1) \ $q^{2k}Q_i-Q_j \in R$ are invertible in $R$ for 
$|k| < n, i \ne j$.
\par\medskip
Note that the condition (8.3.1) is stronger than the condition 
(6.3.1) for any $\Bp$.
It is shown by [SawS, 8.3.2], based on the result in [HS], 
 that $\r_0: \CH \to \ol\CH$ gives 
an isomorphism if the separation condition (8.3.1) holds.  
We have the following corollary.
%%%
\begin{cor}   %%% Corollary 8.4
Suppose that the condition (8.3.1) holds for $\CH$.  
Then $\CH \simeq \ol\CH^{\Bp}$ for any $\Bp$.  
In particular, Theorem 6.8 gives a new presentation for 
the Ariki-Koike algebra $\CH$.
\end{cor}
\begin{proof}
We have $(1^r) \preceq \Bp \preceq (r)$ for any $\Bp$.
Since $\r_0: \CH \to \ol\CH$ is an isomorphism, the map 
$\r_{\Bp'\Bp} : \ol\CH^{\Bp} \to \ol\CH^{\Bp'} \simeq \CH$ is surjective 
by Proposition 8.2 for $\Bp' = (1^r)$.  Since 
both of $\CH$ and $\ol\CH^{\Bp}$ are free $R$-modules of the same 
rank, we obtain $\ol\CH^{\Bp} \simeq \CH$ as asserted.
\end{proof}

\bigskip
%%%%%
%%%%%

 \end{document}